\documentclass[11pt]{article}
\usepackage{sustyle}
\usepackage[round]{natbib}
\usepackage{url}
\usepackage{bm}
\usepackage{multicol}
\usepackage[colorlinks=true, allcolors=black,c itecolor=blue, urlcolor=black]{hyperref}
\usepackage{fullpage}
\usepackage{cleveref}
\usepackage{setspace}
\usepackage{ulem}
\usepackage{booktabs}
\usepackage{color}
\usepackage{xcolor}
\usepackage{enumerate}
\usepackage{ulem}
\usepackage{fancyhdr}
\usepackage{verbatim}
\usepackage{latexsym}
\usepackage{float}
\numberwithin{equation}{section}

\newcommand{\Normal}{{\mathcal{N}}}
\newcommand\indep{\protect\mathpalette{\protect\independenT}{\perp}}
\def\independenT#1#2{\mathrel{\rlap{$#1#2$}\mkern2mu{#1#2}}}

\newcommand{\PP}{{\mathbb{P}}}
\newcommand{\EE}{{\mathbb{E}}}

\newcommand{\X}{{\bX}}

\newcommand{\FDP}{{\textnormal{FDP}}}
\newcommand{\TPP}{{\textnormal{TPP}}}
\newcommand{\fdp}{\textnormal{fdp}^{\infty}}
\newcommand{\tpp}{\textnormal{tpp}^{\infty}}
\newcommand{\tppt}{\textnormal{tpp}^{\infty (t)}}
\newcommand{\fdpt}{\textnormal{fdp}^{\infty (t)}}

\newcommand{\fd}{\textnormal{fd}^{\infty}}
\newcommand{\td}{\textnormal{td}^{\infty}}

\newcommand{\VAR}{{\textnormal{Var}}}
\newcommand{\Cov}{{\textnormal{Cov}}}

\def\bM{\bm{M}}
\newcommand{\bth}{{\widehat{\bbeta}}}
\newcommand{\tuj}{{\widetilde{u}_j}}
\newcommand{\tyi}{{\widetilde{Y}_i}}
\newcommand{\bs}{{\bbeta_{S_0}}}
\newcommand{\bss}{{\bbeta_{S_1}}}
\newcommand{\bw}{\bm w}
\newcommand{\zh}{{\widehat{\bw}}}
\newcommand{\Xs}{{\bX_{S_{0}}}}
\newcommand{\Xj}{{X_j}}
\newcommand{\Xss}{{\bX_{S_{1}}}}
\newcommand{\XSc}{{\bX_{S^C}}}
\newcommand{\Ps}{{\mathcal{P}_{S_0^{\perp}}}}

\title{The Price of Competition: Effect Size Heterogeneity\\ Matters in High Dimensions}
\author{Hua Wang{$^\ast$} \and Yachong Yang\footnote{The first two authors are listed in alphabetical order.} \and Weijie J.~Su}
\date{}

\begin{document}

\maketitle

{\centering
\vspace*{-0.5cm}
Department of Statistics and Data Science, University of Pennsylvania
\par\bigskip
\date{June 2020; Revised March 2022}\par
}

\begin{abstract}

In high-dimensional sparse regression, would increasing the signal-to-noise ratio while fixing the sparsity level always lead to better model selection? For high-dimensional sparse regression problems, surprisingly, in this paper we answer this question in the \textit{negative} in the regime of linear sparsity for the Lasso method, relying on a new concept we term \textit{effect size heterogeneity}. Roughly speaking, a regression coefficient vector has high effect size heterogeneity if its nonzero entries have significantly different magnitudes. From the viewpoint of this new measure, we prove that the false and true positive rates achieve the optimal trade-off \textit{uniformly} along the Lasso path when this measure is maximal in a certain sense, and the worst trade-off is achieved when it is minimal in the sense that all nonzero effect sizes are roughly equal. Moreover, we demonstrate that the first false selection occurs much earlier when effect size heterogeneity is minimal than when it is maximal. The underlying cause of these two phenomena is, metaphorically speaking, the ``competition'' among variables with effect sizes of the same magnitude in entering the model. Taken together, our findings suggest that effect size heterogeneity shall serve as an important complementary measure to the sparsity of regression coefficients in the analysis of high-dimensional regression problems. Our proofs use techniques from approximate message passing theory as well as a novel technique for estimating the rank of the first false variable.

\end{abstract}

\section{Introduction}

We consider high-dimensional sparse regression problems in which we observe an $n$-dimensional response vector $\by$ that is generated by a linear model
\begin{equation}\label{eq:basic_model}
\bm{y} = \bm{X\beta} + \bm{z},
\end{equation}
where $\bX$ is an $n\times p$ design matrix of features, $\bbeta \in \R^p$ denotes an unknown vector of regression coefficients, and $\bz \in \R^n$ is a noise term. In the big data era, this model has been increasingly applied to high-dimensional settings where the number of variables $p$ is comparable to or even much larger than the number of observational units $n$. While this reality poses challenges to the regression problem, in many scientific problems there are good reasons to suspect that truly relevant variables account for a small fraction of all the observed variables or, equivalently, $\bbeta$ is sparse in the sense that many of its components are zero or nearly zero. Indeed, a very impressive body of \textit{theoretical} work shows that the difficulty of variable selection in the high-dimensional setting relies crucially on how sparse the regression coefficients are~\citep{DS,bickel2009}.

\begin{figure}[!htb]
\centering
\includegraphics[width = 0.8\linewidth]{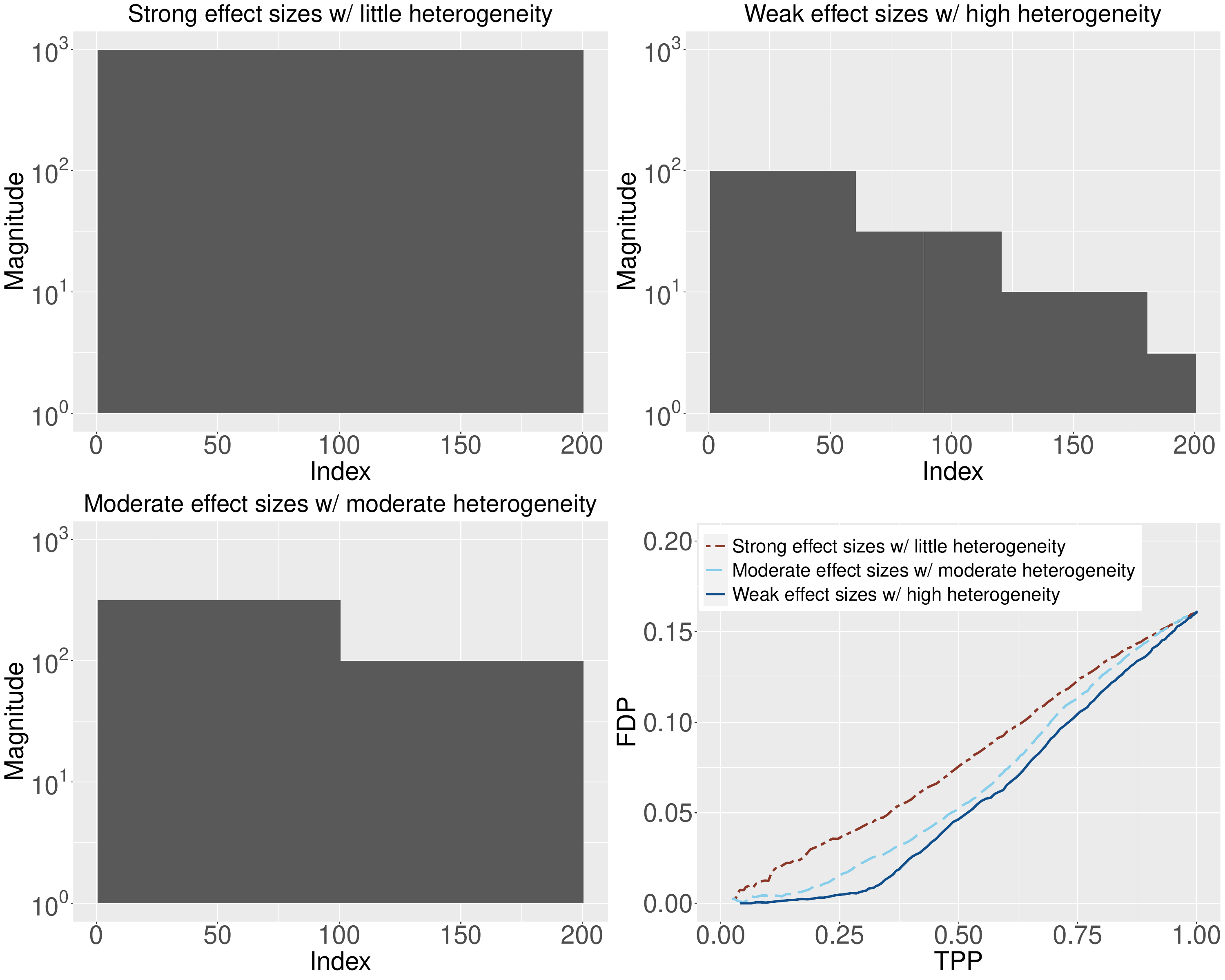}
\caption{The TPP--FDP trade-off along the entire Lasso path, with three different sets of regression coefficients. Note that the TPP--FDP trade-off is equivalent to the receiver operating characteristic curve. The sparsity of $\bbeta$ is fixed to $k = 200$ (throughout this paper, we use $k$ to denote the sparsity level) and the $200$ true effects are plotted in the logarithmic scale in the three panels. For example, in the ``Strong effect sizes w/ little heterogeneity'' setting, $\beta_1 = \cdots = \beta_{200} = 10^3$, and $\beta_{201} = \cdots = \beta_{1000} = 0$. The design matrix $\bX\in \R^{n\times p}$ has independent $\N(0, 1/n)$ entries, where $n = p = 1000$, and the noise term $\bz$ has independent $\N(0, \sigma^2)$ entries with $\sigma = 0.01$. The bottom-right panel shows the plot of FDP as a function of TPP, averaged over $100$ independent runs. For completeness, we remark that effect size heterogeneity influences model selection in a more complex manner at a higher noise level (see \Cref{fig:noisy}).}
\label{fig:3simulations}
\end{figure}

This paper, however, asks whether there are other measures concerning the regression coefficients that have a \textit{practical} impact on variable selection for the linear model~\eqref{eq:basic_model}. To address this question, we present a simulation study in \Cref{fig:3simulations}. Notably, the sparsity---or, put differently, the number of nonzero components---of the regression coefficients $\bbeta$ is \textit{fixed} to $200$ across three experimental settings, but with varying magnitudes of the $200$ true effect sizes. The method we use for variable selection is the Lasso~\citep{lasso}, which is perhaps the most popular model selector in the high-dimensional setting. Given a penalty parameter $\lambda > 0$, this method finds the solution to the convex optimization program
\begin{equation}\label{eq:Lasso}
\widehat{\bm{\beta}}(\lambda) = \argmin_{\bm{b} \in \R^p} \frac{1}{2}\|\bm{y} - \bm{Xb} \|^2+\lambda \|\bm{b}\|_1,
\end{equation}
where $\|\cdot\|$ and $\|\cdot\|_1$ denote the $\ell_2$ and the $\ell_1$ norms, respectively. A variable $j$ is selected by this method at $\lambda$ if $\widehat{\beta}_j(\lambda) \ne 0$, and a false selection occurs if it is a noise variable in the sense that $\beta_j = 0$. Formally, we use the false discovery proportion (FDP) and true positive proportion (TPP) as measures of the type I error and power, respectively, to assess the quality of the selected model $\{1 \le j \le p: \widehat{\beta}_j(\lambda)\ne 0\}$:
\begin{align}\label{eq:fdp_tpp}
\FDP_\lambda &= \frac{\#\{j:\beta_j = 0 ~\text{and}~ \widehat{\beta}_j(\lambda)\ne 0\}}{\#\{j: \widehat{\beta}_j(\lambda)\ne 0\}}, \\ \TPP_\lambda &= \frac{\#\{j:\beta_j \ne 0 ~\text{and}~ \widehat{\beta}_j(\lambda)\ne 0\}}{\#\{j: {\beta}_j\ne 0\}}.
\end{align}
As is clear, we wish to select a model with a small FDP and large TPP.

Despite weaker effect sizes, strikingly, \Cref{fig:3simulations} shows that the Lasso can achieve \textit{better} model selection in terms of the TPP--FDP trade-off and, in particular, this counterintuitive behavior holds uniformly along the entire Lasso path or, equivalently, over all values of $\lambda$. Existing theory often analyzes how the worst-case performance of the Lasso and other related procedures depends on the regression coefficients through the sparsity of the regression coefficients (see, for example, \cite{wainwright2019high}). However, the sparsity level is fixed across the experimental settings of \Cref{fig:3simulations}. In light of this, therefore, one would expect that the strong signals and weak signals would yield the best and worst model selection results, respectively. \Cref{fig:3simulations} shows that this is not necessarily the case.



Thus, a \textit{finer-grained} structural analysis of the effect sizes is needed to better understand the Lasso in some settings. In this paper, we address this important question by proposing a concept that we term \textit{effect size heterogeneity} concerning the regression coefficients in high dimensions. Roughly speaking, a regression coefficient vector has higher effect size heterogeneity than another vector (of the same sparsity) if the nonzero entries of the former are more heterogeneous than those of the latter in terms of magnitude. As a complement to sparsity, effect size heterogeneity will be shown to have a significant impact on how the Lasso performs in terms of the false and true positive rates trade-off: while the sparsity level of the regression coefficients is fixed, the higher the effect size heterogeneity is, the better the Lasso performs. Turning back to \Cref{fig:3simulations}, we note that the strong effect sizes are the least heterogeneous in magnitude, and the weak effect sizes
are the most heterogeneous. Therefore, the comparisons made in \Cref{fig:3simulations} match well the implication of effect size heterogeneity.

\begin{figure}[!htb]
	\centering	
	\includegraphics[width = 0.52\linewidth]{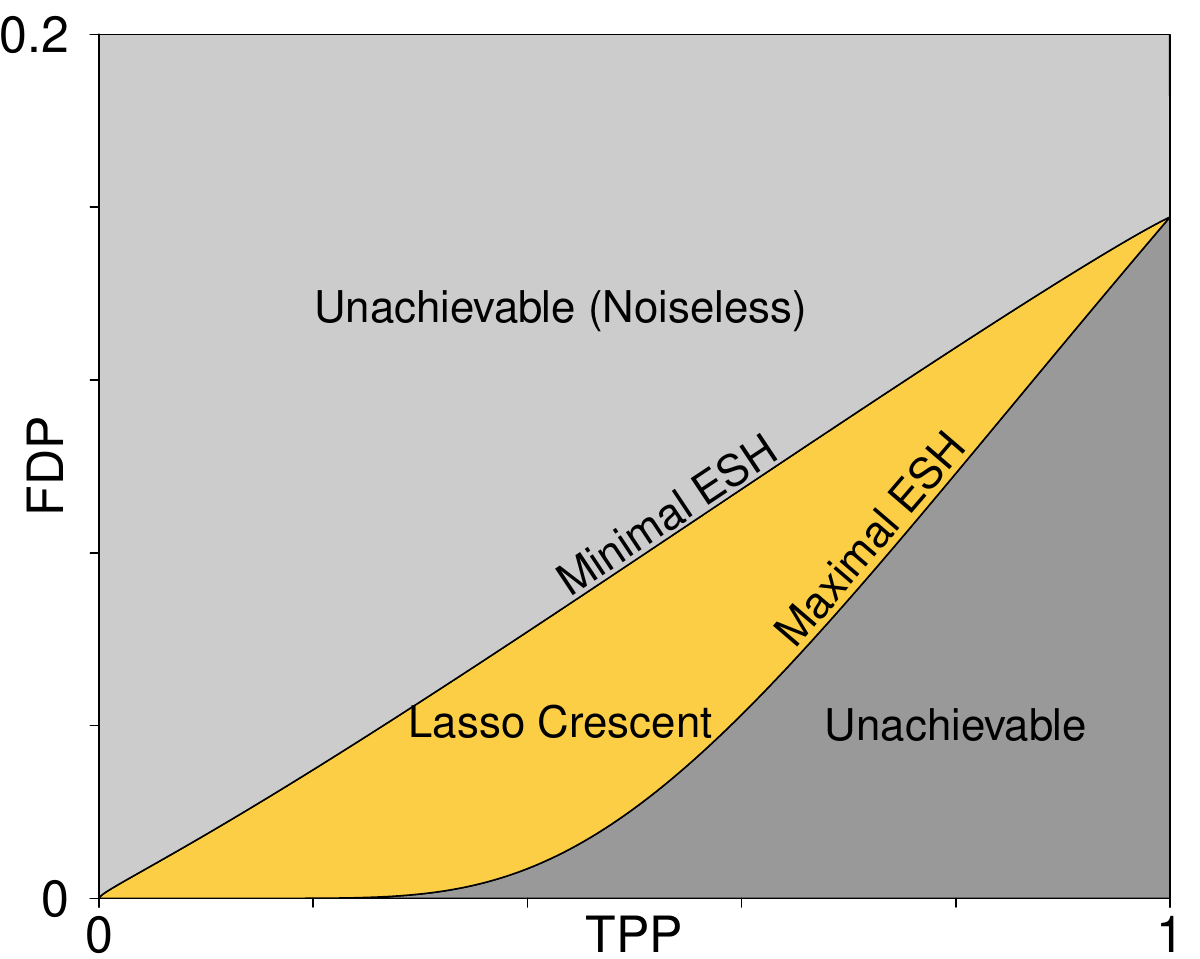}
	\caption{The Lasso Crescent diagram specified by the parameters $n/p = 1$ and $k/p = 0.2$, following the setting in \Cref{fig:3simulations}. The lower/upper smooth curve is asymptotically achieved with maximal/minimal effect size heterogeneity (ESH) in the regime of infinite signal-to-noise ratio. The explicit expressions of the curves are given in \Cref{sec:lasso-crescent}. Our \Cref{thm:lower} implies that nowhere on the Lasso path we can find any (TPP, FDP) pairs in the region below the Lasso Crescent (also see \cite{su2017false}). In the noiseless setting, moreover, this impossibility result continues to hold in the region above the Lasso Crescent (\Cref{thm:upper}), which is labeled ``Unachievable (Noiseless).''}
	\label{fig:PhaseDiagram}
\end{figure}

Concretely, the main thrust of this paper lies in the development of two complementary perspectives to precisely quantify the impact of effect size heterogeneity. First, following the setup of \Cref{fig:3simulations}, we consider the full possible range of the asymptotic trade-off between the TPP and FDP along the Lasso path, while varying the level of effect size heterogeneity. Assuming a random design with independent Gaussian entries and working in the regime of linear sparsity---meaning that the fraction of true effect sizes tends to a constant---we formally show that the full possible range is enclosed by two smooth curves in the (TPP, FDP) plane, which we referred to as the \textit{Lasso Crescent}. \Cref{fig:PhaseDiagram} presents an instance of the Lasso Crescent. More precisely, having excluded the impact of noise by taking $\bm z = \bm 0$ in the linear model~\eqref{eq:basic_model}, the lower curve is asymptotically achieved when effect
size heterogeneity is maximal in the sense that all true effect sizes are widely different from each other, while the upper curve is asymptotically achieved when the heterogeneity is minimal in the sense that all true effects are of the same size. In general, the $(\TPP, \FDP)$ pairs computed from the entire Lasso path must be asymptotically sandwiched between the two curves in the noiseless setting or, equivalently, in the regime of the infinite signal-to-noise ratio. The gap between the two curves is fundamental in the sense that it persists no matter how strong the effects are.

While the TPP--FDP trade-off essentially examines the ``bulk'' of the Lasso solution path, the second perspective we take extends to the ``edge'': when does the first noise variable enter the model along the Lasso path? More precisely, we decrease $\lambda$ from $\infty$ to $0$ and find the first ``time'' a false selection occurs. To indicate the difficulty of consistent model selection, formally, we consider the rank of the first noise variable or, put concretely, one plus the number of the true variables before the Lasso selects the first noise variable. Intuitively, a large rank is desirable. As with the first perspective, assuming a Gaussian random design and regression coefficients with linear sparsity, we prove that the rank of the first false selection is bounded above by $(1 + o(1)) n/(2 \log p)$. This upper bound, which approximately equals 72 in the setting of \Cref{fig:3simulations}, holds no matter how strong the effect sizes are. Interestingly, this upper bound is
\textit{exactly} achieved when effect size heterogeneity is maximal and the noise level tends to zero. On the other hand, \cite{su2018first} has obtained a sharp prediction of the rank of the first false variable in the case of minimal effect size heterogeneity, which, together with our new result, shows that the first noise variable occurs much earlier with minimal effect size heterogeneity than with maximal effect size heterogeneity. Although not entirely related, the two perspectives consistently demonstrate that effect size heterogeneity is an important and useful concept for understanding the performance of the Lasso.


The fact that effect size heterogeneity matters, as shown above, is due to the bias introduced by the shrinkage nature of the Lasso. This bias in turn makes the residuals absorb many of the true effects that act as what we may want to call ``shrinkage noise''. Metaphorically, variables yet to be selected tend to ``compete'' with each other in entering the Lasso path and contribute to the shrinkage noise. The ``competition'' is particularly intensive among variables having about the same effect sizes, which is the case when effect size heterogeneity is low. As a price, the shrinkage noise gets inflated and some noise variables may be selected early due to their high correlations with the shrinkage noise. This is why false selections occur with a good chance and early. In contrast, when the heterogeneity is high, the largest true effect yet to be selected tends to have a significant correlation with the residuals, thereby having a better chance to be
selected sooner. To appreciate this heuristic explanation, it is instructive to note that the least-squares estimator, if available, does not exhibit this price-of-competition phenomenon, as it is unbiased for the regression coefficients.\footnote{If the effect sizes are sufficiently strong, variable selection using the $t$ values of the least-squares estimator can lead to full power without any type I errors.} An alternative but less direct way to appreciate effect size heterogeneity is to relate it to the restricted eigenvalue condition~\citep{bickel2009}. Roughly speaking, this condition is concerned with a vector of regression coefficients such that its $\ell_1$ norm is largely contributed by a few components, and such approximate sparse regression coefficients can be well estimated by the Lasso and the Dantzig selector under certain designs~\citep{bickel2009,raskutti2010restricted}. From the viewpoint of this condition, therefore, regression coefficients with high effect size heterogeneity can be thought of as having a smaller effective sparsity level, which is favored by the Lasso.

As a final remark, the price-of-competition phenomenon does not appear if the sparsity is sub-linear in the ambient dimension $p$, which is often assumed in the copious body of literature on high-dimensional regression. In this regime of sparsity, effect size heterogeneity has a vanishing impact on the performance of the Lasso if the signal-to-noise ratio is sufficiently strong or the beta-min condition is satisfied. Our paper also departs from this line of literature from a technical standpoint. Indeed, the proofs of the results in this paper make heavy use of approximate message passing (AMP) theory~\citep{DMM09,ampdense, BM12}, with nontrivial extensions.

\subsection{Organization}
\label{sec:organ-notat}
The remainder of this paper is organized as follows. In \Cref{sec:lasso-crescent}, we formalize the Lasso Crescent diagram by presenting our theoretical results that predict the TPP--FDP trade-off with respect to effect size heterogeneity. \Cref{sec:heuristics} extends the investigation of effect size heterogeneity to the problem of the first false variable along the Lasso path. \Cref{sec:proofs} is devoted to proving the results in \Cref{sec:lasso-crescent}, whereas technical details of the proofs are deferred to the appendix. In \Cref{sec:illustrations}, we provide numerical studies to demonstrate the impact of effect size heterogeneity in general settings.  We conclude the paper in \Cref{sec:discussion} with a few directions for future research.


\section{The Lasso Crescent}
\label{sec:lasso-crescent}

In this section, we derive the full possible range of the asymptotic trade-off between the TPP and FDP along the Lasso path, with a focus on its dependence on effect size heterogeneity. Specifically, our results can be pictorially presented by the Lasso Crescent as in \Cref{fig:PhaseDiagram}, hence the title of this section. The proofs are deferred to \Cref{sec:proofs}.

Throughout this section, and indeed the entire paper, we assume the following \textit{working hypotheses} to specify the linear model~\eqref{eq:basic_model}. For ease of reading, we use boldface letters to denote vectors and matrices.
\smallskip\\
\noindent\textit{Gaussian design matrix.} We consider a sequence of designs $\bX \in \R^{n_l \times p_l}$ consisting of i.i.d.~$\N(0,1/n_l)$ entries so that each column has an approximate unit $\ell_2$ norm. As the index $l \goto \infty$, we assume $p_l, n_l \goto \infty$ with $n_l/p_l \rightarrow \delta$ for a constant $\delta > 0$. The index $l$ is often omitted for the sake of simplicity.
\smallskip\\
\noindent\textit{Regression coefficients.} Let the regression coefficients $\beta_1, \ldots, \beta_p$ be i.i.d.~copies of a random variable $\Pi$ that satisfies $\E \Pi^2 < \infty$. Of particular interest to this paper is an $\epsilon$-sparse prior $\Pi$ in the sense that $\PP(\Pi \neq 0) = \epsilon$ for a constant $0 < \epsilon < 1$. Thus, the realized $\beta_1, \ldots, \beta_p$ are in the linear sparsity regime since the sparsity is approximately equal to $\epsilon p$. 
\smallskip\\
\noindent\textit{Noise.} The noise term $\bm z$  consists of i.i.d.~elements drawn from $\N(0, \sigma^2)$, where the noise level $\sigma \ge 0$ is fixed.
\smallskip\\
For completeness, $\bX, \bbeta$, and $\bz$ are jointly independent. These assumptions are used in the literature on AMP theory and its applications (see, for example, \citet{ampdense,BM12,maleki2013asymptotic, bayati2013estimating,bu2019algorithmic}) and, more recently, have been commonly made in the high-dimensional regression literature \citep{weinstein2017power,mousavi2018consistent,weng2018overcoming, sur2019likelihood}. On top of that, we adopt some adjustments made by \cite{su2017false} that slightly simplify the assumptions on $\bbeta$ and $\bz$. Regarding the assumption on the noise, it is worth noting that we do not exclude the case $\sigma = 0$, which corresponds to noiseless observations. For some of the results in this section, the price-of-competition phenomenon manifests itself most clearly in the noiseless setting.

\subsection{Most heterogeneous effect sizes}
\label{sec:most-favorable-prior}


Our first main theorem considers regression coefficients that are drawn from the following prior distribution:
\begin{definition}
For $M > 0$ and an integer $m > 0$, we call
\begin{equation}\label{eq:eps_m_M_prior}
  \Pi^\Delta =
    \begin{cases}
    0 & \text{w.p. $1-\epsilon$}\\
    M & \text{w.p. $\frac{\epsilon}{m}$}\\
    M^2 & \text{w.p. $\frac{\epsilon}{m}$}\\
    \cdots & \cdots\\
    M^m & \text{w.p. $\frac{\epsilon}{m}$}
    \end{cases}
\end{equation}
the $(\epsilon, m, M)$-heterogeneous prior.
\end{definition}

For notational convenience, we suppress the dependence of $\Pi^\Delta$ on $\epsilon, m, M$. This prior is $\epsilon$-sparse in the sense of the working hypotheses. As is clear, larger values of $m, M$ would render the prior more heterogeneous. Indeed, this paper is primarily concerned with the case where both $M, m \goto \infty$. This corresponds to the regime where the signal-to-noise ratio tends to infinity and, in addition, the true effect sizes are increasingly different. To be complete, the $(\epsilon, m, M)$-heterogeneous prior is only a specific example that attains increasing heterogeneity. See Remark \ref{remark:most_favorable_priors} for more examples.



Following \eqref{eq:fdp_tpp}, $\FDP_{\lambda}(\Pi)$ and $\TPP_{\lambda}(\Pi)$ denote the (random) false discovery proportion and true positive proportion, respectively, of the Lasso estimate at $\lambda$ when the regression coefficients in \eqref{eq:basic_model} are i.i.d.~draws from a prior $\Pi$. For ease of reading, we say a pair $(\TPP, \FDP)$ \textit{outperforms} another pair $(\TPP', \FDP')$ if $\TPP > \TPP'$ and $\FDP < \FDP'$. As noted earlier, all theoretical results in this paper are obtained under the working hypotheses. For conciseness, the statements of our theorems shall not mention this fact anymore. 

\begin{theorem}\label{thm:lower}
Let $C> c > 0$ be fixed. For any $\epsilon$-sparse prior $\Pi$, if both $m$ and $M$ are sufficiently large in the $(\epsilon, m, M)$-heterogeneous prior $\Pi^\Delta$, then the following conclusions are true:
\begin{enumerate}
 
\item[(a)] The event
\[
\bigcup_{c < \lambda,\lambda' < C} \left\{ \left( \TPP_{\lambda'}(\Pi), \FDP_{\lambda'}(\Pi) \right) \text{ outperforms } \left( \TPP_{\lambda}(\Pi^\Delta), \FDP_{\lambda}(\Pi^\Delta) \right) \right\}
\]
happens with probability tending to \textbf{zero} as $n, p \goto \infty$.

\item[(b)] For any constant $\nu > 0$,  no matter how we choose $\widehat{\lambda'}(\by, \bX) \ge c$ adaptively as long as it always satisfies $\TPP_{\widehat{\lambda'}}(\Pi) > \nu$, with probability approaching \textbf{one} there exists $\widehat{\lambda} > 0$ such that
\[
\left( \TPP_{\widehat\lambda}(\Pi^\Delta), \FDP_{\widehat\lambda}(\Pi^{\Delta}) \right) \text{ outperforms } \left( \TPP_{\widehat{\lambda'}}(\Pi), \FDP_{\widehat{\lambda'}}(\Pi) \right).
\]

\end{enumerate}
\end{theorem}

\begin{remark}\label{remark:most_favorable_priors}
The priors for which the theorem holds can be extended in the following way. Consider a sequence of priors $\Pi^\Delta$ satisfying $\Pi^\Delta = 0$ with probability $1 - \epsilon$ and $\Pi^\Delta = M_i \ne 0$ with probability $\epsilon\gamma_i$ for $i = 1, \ldots, m$ such that $\gamma_1 + \cdots + \gamma_m = 1, \max_i \gamma_i \goto 0$, and $\min_{1 \le i \le m} \left| M_i/M_{i-1} \right| \goto \infty$ (set $M_0 = 1$). Alternatively, the nonzero component of the prior can be drawn from a continuous random variable with CDF 
\[
\frac{\log_M x}{m}
\]
for $1 \le x \le M^m$. While the theorem statement is restricted to $(\epsilon, m, M)$-heterogeneous priors for brevity, its proof considers the general case.


\end{remark}

This theorem demonstrates the optimality of heterogeneous and strong effects in terms of the trade-off between the TPP and FDP. Importantly, this optimality is \textit{uniform} in the sense that it holds along the entire Lasso path, no matter how strong the true effects coming from $\Pi$ are. To be sure, the event as a union in (a) is taken over any $(\TPP, \FDP)$ pair from the prior $\Pi$ and any pair from the prior $\Pi^\Delta$. Although each conclusion alone is not a consequence of the other, as we will see from the proof in \Cref{sec:proofs}, the two conclusions are built on top of the fact that the pairs $(\TPP_{\lambda}, \FDP_{\lambda})$ with varying $\lambda$ converge uniformly to a deterministic smooth curve for both $\Pi$ and $\Pi^\Delta$. This fact allows us to obtain the following byproduct:

\begin{proposition}\label{prop:lower}
Under the assumptions of \Cref{thm:lower}, for any sufficiently small  constant $\nu > 0$, the following statement holds with probability tending to one: for any $\lambda, \lambda' > c$ such that $\left| \TPP_{\lambda}(\Pi^\Delta) - \TPP_{\lambda'}(\Pi) \right| < \nu$ and $\TPP_{\lambda'}(\Pi) > 0.001$, we have
\[
\FDP_{\lambda}(\Pi^{\Delta}) < \FDP_{\lambda'}(\Pi).
\]
\end{proposition}

This result makes it self-evident why the prior $\Pi^\Delta$ is a most favorable for the entire Lasso path, though literally, we should interpret this favorability in the limit $m, M \goto \infty$. More precisely, this result implies that given a required power level, the smallest possible FDP is achieved when the effects are increasingly heterogeneous and strong. Of note, the number $0.001$ above can be replaced by any small positive constant, and it does not impede the interpretability of the theorem since we are generally not interested in a model that includes only a tiny fraction of true variables.

An interesting yet unaddressed question is to find an expression of the asymptotic minimum of FDP given $\TPP_{\lambda}(\Pi^\Delta) = u$ in the limit $m, M \goto \infty$. Call this function $q^\Delta(u; \delta, \epsilon)$. From our results, one can easily see that $q^\Delta$ is nothing but the lower \textit{envelope} of instance-specific TPP--FDP trade-off curves over all $\epsilon$-sparse priors. To see this, first note that one can prove that as $n,p \goto \infty$, the pairs $(\TPP_{\lambda}(\Pi), \FDP_{\lambda}(\Pi))$ over all $\lambda$ converge to a smooth curve, which is denoted by $q^{\Pi}(u)$ (see \Cref{sec:proofs}). Recognizing that $\Pi^\Delta$ is also $\epsilon$-sparse and assuming $\lim_{m, M \goto \infty} q^{\Pi^\Delta}$ exists, we must have
\begin{equation}\label{eq:def_q}
q^\Delta(u): = \lim_{m, M \goto \infty} q^{\Pi^\Delta}(u)  \ge \inf_{\Pi: \epsilon\text{-sparse}} q^\Pi(u).
\end{equation}
On the other hand, it follows from \Cref{thm:lower} and in particular Proposition~\ref{prop:lower} that
\[
q^\Pi(u) \ge \lim_{m, M \goto \infty} \left( q^{\Pi^\Delta}(u) + o(1) \right) = q^\Delta(u)
\]
for any $\epsilon$-sparse prior $\Pi$. This display, together with \eqref{eq:def_q}, gives
\begin{equation}\label{eq:envelope_sim}
q^\Delta(u) = \inf_{\Pi: \epsilon\text{-sparse}} q^\Pi(u).
\end{equation}
Interestingly, the right-hand side of \eqref{eq:envelope_sim} has been tackled in \cite{su2017false}, leading to a precise expression. To describe this expression, let $t^\Delta(u)$ be
the largest positive root of the equation in $t$,
\begin{equation}\label{eq:t_lower}
\frac{2(1 - \epsilon)\left[ (1+t^2)\Phi(-t) - t\phi(t) \right] + \epsilon(1 + t^2) - \delta}{\epsilon\left[ (1+t^2)(1-2\Phi(-t)) + 2t\phi(t) \right]} = \frac{1 - u}{1 - 2\Phi(-t)},
\end{equation}
where $\Phi(\cdot)$ and $\phi(\cdot)$ denote the CDF and PDF of the standard normal distribution, respectively. Theorem 2.1 in \cite{su2017false} shows that
\begin{equation}\label{eq:q_lower}
\inf_{\Pi: \epsilon\text{-sparse}} q^\Pi(u) = \frac{2(1-\epsilon)\Phi(-t^\Delta(u))}{2(1-\epsilon)\Phi(-t^\Delta(u)) + \epsilon u}.
\end{equation}
Taken together, \eqref{eq:envelope_sim} and \eqref{eq:q_lower} yield
\begin{equation}\label{eq:q_delta_def}
q^\Delta(u) = \frac{2(1-\epsilon)\Phi(-t^\Delta(u))}{2(1-\epsilon)\Phi(-t^\Delta(u)) + \epsilon u}.
\end{equation}
\begin{remark}
If $u = 0$, treat $\infty$ as a root of the equation and set $0/0=0$ in \eqref{eq:q_lower}. As such, $q^\Delta$ satisfies $q^\Delta(0) = 0$. If $\delta < 1$ and $\epsilon$ is larger than a threshold determined by $\delta$, the function $q^\Delta$ is defined only for $u$ between 0 and a certain number strictly smaller than 1. This is where the celebrated Donoho--Tanner phase transition occurs~\citep{donoho2009observed} (also see \Cref{sec:upper_additional_proof}). Throughout this paper, however, we focus on the regime that is below the Donoho--Tanner phase transition---that is, the case where $\delta \ge 1$, or $\delta < 1$ and $\epsilon$ is small so that the range of $u$ is the unit interval $[0, 1]$. In contrast, above the phase transition, the mapping from the TPP to FDP might not be unique (see Figure 3.1 in \cite{wang2020bridge} and \cite{wang2020complete}).
\end{remark}

In summary, we have the following corollary, which addresses the aforementioned question. 
\begin{corollary}\label{cor:lower}
Under the assumptions of \Cref{thm:lower}, we have
\[
\lim_{m, M \goto \infty} \lim_{n,p \goto \infty}\sup_{\lambda > c} \left| \FDP_{\lambda}(\Pi^\Delta) - q^\Delta\left( \TPP_{\lambda}(\Pi^\Delta) \right) \right| = 0,
\]
where $\lim_{n,p \goto \infty}$ is taken in probability. Moreover, for any $\epsilon$-sparse prior $\Pi$, we have
\[
\FDP_{\lambda}(\Pi) \ge q^\Delta \left( \TPP_{\lambda}(\Pi) \right) - 0.001
\]
for all $\lambda > c$ with probability tending to one.

\end{corollary}

\begin{remark}
As $\lambda \goto \infty$, both $\TPP_{\lambda}(\Pi^\Delta)$ and $\FDP_{\lambda}(\Pi^\Delta)$ tend to 0. Hence, there is no need to impose an upper bound on $\lambda$ when taking the supremum $\sup_{\lambda > c}$. The second conclusion of Corollary~\ref{cor:lower} follows from Proposition~\ref{prop:lower} in conjunction with the continuity of $q^\Delta$. As earlier, $0.001$ can be replaced by any positive constant.
\end{remark}

The second conclusion of Corollary~\ref{cor:lower} is part of Theorem 2.1 in~\cite{su2017false} and demonstrates that true variables and irrelevant variables are always interspersed along the Lasso path. In particular, this is true when the regularization parameter $\lambda$ tends to 0. In this case, indeed, the Lasso would select a significant fraction of false variables with vanishing but nonzero estimated coefficients. This fact necessitates a form of calibration of the Lasso estimates for variable selection~\cite{wang2020bridge,wang2020complete}.

The significance of \Cref{thm:lower} and Corollary~\ref{cor:lower}, however, extends beyond earlier results. Precisely, \cite{su2017false} derived the expression \eqref{eq:q_lower} by constructing a different signal prior $\Pi$ for each power level $u$. Indeed, the nonzero component of the prior constructed in \cite{su2017false} has two different magnitudes with weights depending on $u$, as opposed to an increasing number of different magnitudes as in the $(\epsilon, m, M)$-heterogeneous prior. The increasing level of heterogeneity allows us to give a one-shot construction of most heterogeneous priors at all power levels.

\subsection{Least heterogeneous effect sizes}
\label{sec:least-favor-prior}

We now turn to the opposite question: which effect sizes lead to the worst trade-off between the TPP and FDP along the Lasso path? Inspired by the interpretation of effect size heterogeneity, it is natural to consider the following signal prior as a candidate:
\begin{definition}
For $M > 0$, we call
\begin{equation}\label{eq:least}
  \Pi^\nabla =
    \begin{cases}
    0 & \text{w.p. $1-\epsilon$}\\
    M & \text{w.p. $\epsilon$}
    \end{cases}
\end{equation}
the $(\epsilon, M)$-homogeneous prior.
\end{definition}
This prior would render all true effect sizes equal, thereby being least heterogeneous or most homogeneous among all $\epsilon$-sparse priors. The following theorem confirms our intuition that this homogeneous prior is least favorable for the Lasso as the resulting effect sizes give the least optimal trade-off between false positives and power.

\begin{theorem}\label{thm:upper} 
Let $C > c > 0$ be fixed. In the noiseless setting---that is, $z=0$---for any $\epsilon$-sparse prior $\Pi$ that is non-constant conditional on $\Pi \ne 0$, the following conclusions are true for the $(\epsilon, M)$-homogeneous prior $\Pi^\nabla$:
\begin{enumerate}
\item[(a)] The event
\[
\bigcup_{c < \lambda, \lambda' < C} \left\{ \left( \TPP_{\lambda}(\Pi^\nabla), \FDP_{\lambda}(\Pi^{\nabla}) \right) \text{ outperforms } \left( \TPP_{\lambda'}(\Pi), \FDP_{\lambda'}(\Pi) \right) \right\}
\]
happens with probability tending to \textbf{zero} as $n, p \goto \infty$.

\item[(b)] For any constant $\nu > 0$, no matter how we choose $\widehat{\lambda'}(\by, \bX) \ge c$ adaptively as long as it always satisfies $\TPP_{\widehat{\lambda'}}(\Pi) >  \nu$, with probability tending to \textbf{one} there exists $\widehat\lambda > 0$ such that
\[
\left( \TPP_{\widehat{\lambda'}}(\Pi), \FDP_{\widehat{\lambda'}}(\Pi) \right) \text{ outperforms } \left( \TPP_{\widehat\lambda}(\Pi^\nabla), \FDP_{\widehat\lambda}(\Pi^{\nabla}) \right).
\]

\end{enumerate}
\end{theorem}

This theorem is similar to, but in the opposite sense to, \Cref{thm:lower}. One distinction between the two theorems is that \Cref{thm:upper} assumes the noiseless setting, as opposed to the noisy setting considered in \Cref{thm:lower}. The noiseless setting is equivalent to an infinite value of the signal-to-noise ratio, which allows us to better isolate the impact of effect size heterogeneity from that of the noise term. That said, this theorem remains true in the presence of noise by setting a sufficiently large magnitude $M$ for the true effect sizes.

Just as Proposition~\ref{prop:lower} does, the following result follows from the proof of 
\Cref{thm:upper} presented in \Cref{sec:proofs}.

\begin{proposition}\label{prop:upper}
Under the assumptions of \Cref{thm:upper}, for any sufficiently small constant $\nu > 0$, the following statement holds with probability tending to one: if $\lambda, \lambda' > c$ satisfy $\TPP_{\lambda'}(\Pi) > 0.001$ and $\left| \TPP_{\lambda}(\Pi^\nabla)  - \TPP_{\lambda'}(\Pi) \right| < \nu$, then
\[
\FDP_{\lambda}(\Pi^{\nabla}) > \FDP_{\lambda'}(\Pi).
\]
\end{proposition}

As is clear, this result demonstrates that the prior $\Pi^\nabla$ is least favorable for the entire Lasso path in the noiseless case. Roughly speaking, this proposition shows that if the two Lasso problems agree on the value of the TPP along their paths, then the $(\epsilon, M)$-homogeneous prior $\Pi^\nabla$ must yield a higher FDP. As with Proposition~\ref{prop:lower}, $0.001$ can be replaced by any positive constant. On a related note, the prior \eqref{eq:least} is known to be least favorable for certain estimation problems both in the noiseless and noisy settings \cite{donoho2011noise} (see also Lemma 4.4.1 and Corollary 4.4.3 in \cite{maleki2010approximate}). An interesting direction for future research is to study the relationship between estimation and variable selection with regard to the least favorability of the prior distribution.

The sharp distinction between Propositions~\ref{prop:lower} and \ref{prop:upper} must be attributed to the priors $\Pi^\Delta$ and $\Pi^\nabla$. The cause is, loosely speaking, due to the ``competition'' among variables with about the same effect sizes in entering the Lasso model. However, we find it easier to elucidate the underlying cause when studying the rank of the first false variable and thus defer the detailed discussion to \Cref{sec:heuristics}.

We now proceed to specify the curve on which $(\TPP_{\lambda}(\Pi^{\nabla}), \FDP_{\lambda}(\Pi^{\nabla}))$ lies in the limit. For a fixed $\alpha$, let $\varsigma = \varsigma(\alpha)$ denote the largest root of the equation
\begin{equation}\nonumber
\begin{aligned}
\delta &= 2(1-\epsilon)[(1+\alpha^2)\Phi(-\alpha)-\alpha\phi(\alpha)] - \epsilon(2\alpha+\varsigma)\phi(\varsigma) + \epsilon\varsigma \phi(2\alpha+\varsigma) \\
&+ \epsilon(1+\alpha^2)[\Phi(\varsigma)+\Phi(-2\alpha-\varsigma)] + \epsilon(\varsigma+\alpha)^2[\Phi(-\varsigma) + \Phi(-2\alpha-\varsigma)],
\end{aligned}
\end{equation}
and let  $t^\nabla = t^\nabla(u;\delta, \epsilon)$ be the largest root of the following equation in $\alpha$:
\begin{equation}\nonumber
\Phi(\varsigma(\alpha)) + \Phi(-2\alpha-\varsigma(\alpha)) = u.
\end{equation}
With all of these in place, define
\begin{equation}\label{eq:q_upper}
q^\nabla(u;\delta,\epsilon) = \frac{2(1-\epsilon)\Phi(-t^\nabla(u))}{2(1-\epsilon)\Phi(-t^\nabla(u)) + \epsilon u}.
\end{equation}
The derivation of the expression is given in Lemma \ref{lem:upper} in \Cref{sec:upper_additional_proof}. The following result shows that this function describes the limiting trade-off between the TPP and FDP in the case of minimal effect size heterogeneity:

\begin{corollary}\label{cor:upper}
Under the assumptions of \Cref{thm:upper}, we have
\[
\lim_{n, p \goto \infty} \sup_{\lambda > c} \left| \FDP_{\lambda}(\Pi^\nabla) - q^\nabla\left( \TPP_{\lambda}(\Pi^\nabla) \right) \right| = 0.
\]
Moreover, for any $\epsilon$-sparse prior $\Pi$, we have
\[
\FDP_{\lambda}(\Pi) \le q^\nabla \left( \TPP_{\lambda}(\Pi) \right) + 0.001
\]
for all $\lambda > c$ with probability tending to one.

\end{corollary}

As with Theorem~\ref{thm:upper}, Corollary~\ref{cor:upper} holds for any $M > 0$ because of the noiseless setting.

Taken together, Corollaries~\ref{cor:lower} and \ref{cor:upper} give the following result:
\begin{theorem}\label{thm:slides}
Let $c > 0$ be any small constant. In the noiseless setting, for any $\epsilon$-sparse prior $\Pi$, we have
\[
q^\Delta \left( \TPP_{\lambda}(\Pi) \right) - 0.001 \le \FDP_{\lambda}(\Pi) \le q^\nabla \left( \TPP_{\lambda}(\Pi) \right) + 0.001
\]
for all $\lambda > c$ with probability tending to one.
\end{theorem}

\begin{figure}[!htb]
\centering
\begin{minipage}{0.48\textwidth}
\centering
\includegraphics[width = \linewidth]{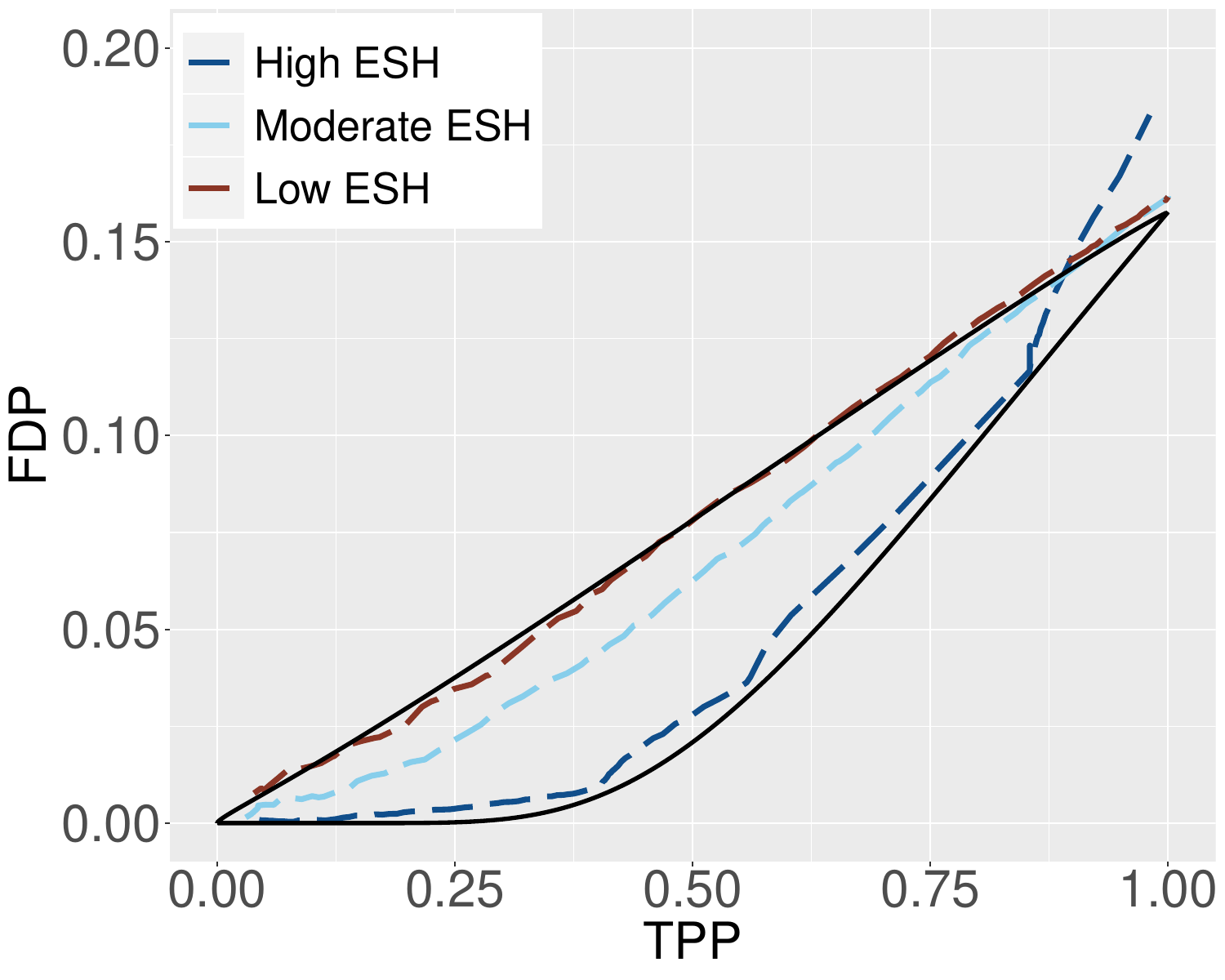}
\end{minipage}
\begin{minipage}{0.48\textwidth}
\centering
\includegraphics[width = \linewidth]{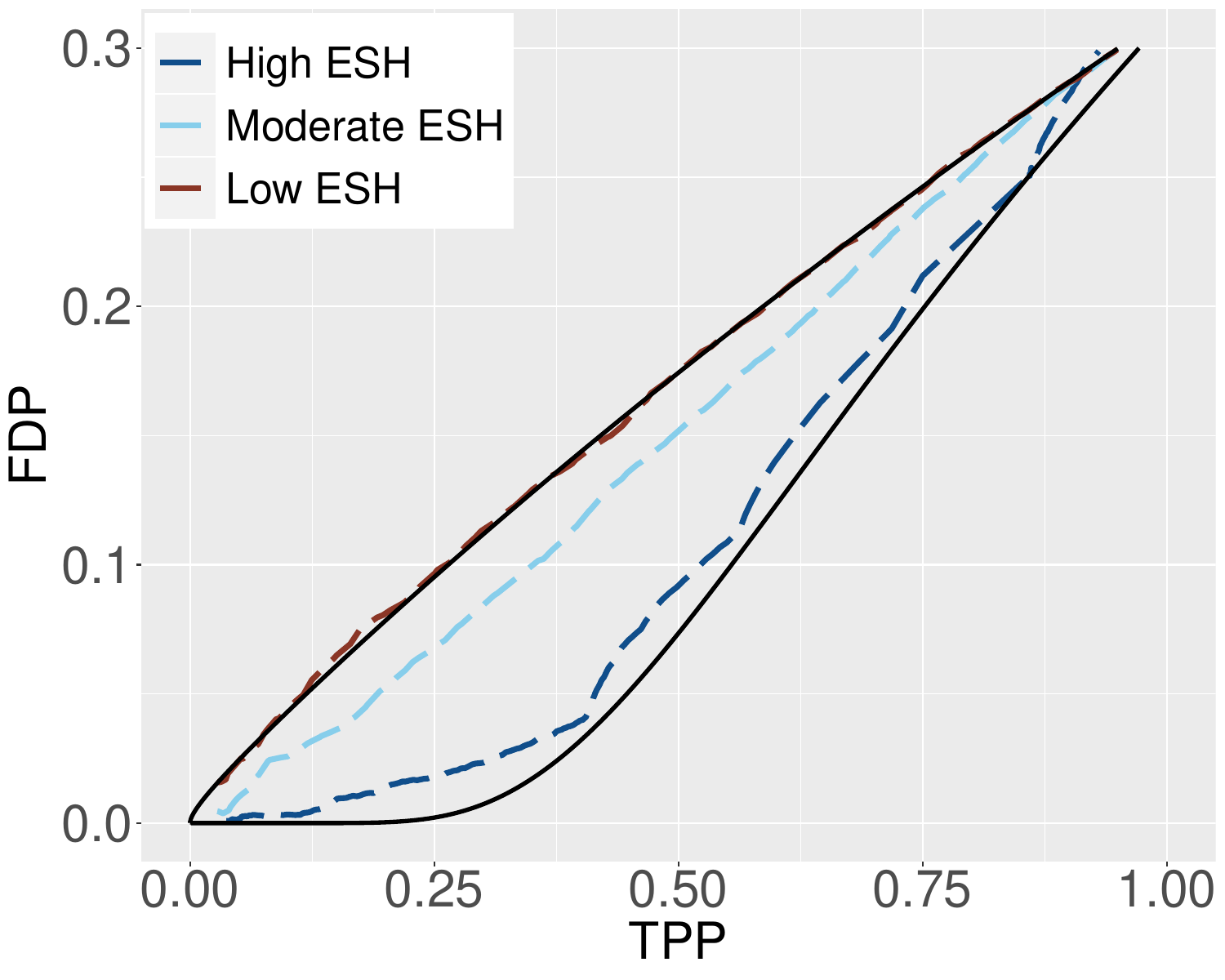}
\end{minipage}
\caption{Illustration of the interpretation of the Lasso Crescent via \Cref{thm:slides}. The design of size $n \times p$ has i.i.d.~$\N(0,1/n)$ entries and the noise level is set to 0. Specifically, we use $n = p = 1000$, and sparsity $k = 200$ in the left panel, and $n = 800, p = 1200$, $k = 200$ in the right panel. The ``high effect size heterogeneity (ESH)'' setting: the 200 coefficients take 4 different values; The ``moderate ESH'' setting: the first 100 coefficients are set to 100 and the second 100 coefficients are set to 50; The ``low ESH'' setting: the 200 coefficients are set to 100. The dashed lines are averaged over 200 independent runs of the Lasso path. The two boundaries $q^\Delta$ and $q^\nabla$ are in solid black lines.}
\label{fig:multi_upper_lower}
\end{figure}

The two curves $q^\Delta$ and $q^\nabla$ enclose a crescent-shaped region, which we call the Lasso Crescent. This theorem shows any $(\TPP, \FDP)$ pairs along the entire Lasso path would essentially lie in the corresponding Lasso Crescent that is specified by the shape $n/p$ of the design and the sparsity ratio $k/p$ of the effect sizes, and this region is tight. \Cref{fig:multi_upper_lower} presents two instances of the Lasso Crescent, with simulations showing good agreement between the predicted and observed behaviors.\footnote{\texttt{R} and \texttt{Matlab} code to calculate $q^\Delta$ and $q^\nabla$ is available at \url{https://github.com/huawang-wharton/effectsizeheterogeneity}.}

\section{The First False Selection}
\label{sec:heuristics}

In this section, we examine the impact of effect size heterogeneity on model selection by the Lasso from a different perspective: when is the first false variable selected? Intuitively, the later the first false variable occurs, the better the method performs. Using a mix of new and old results, this section will show that the first false variable occurs much earlier when effect size heterogeneity is minimal than when it is maximal.

Denote the rank of the first falsely selected variable by
\[
T := \#\{j: \widehat\beta_j(\lambda^\ast - 0) \ne 0 \}  = \#\{j: \widehat\beta_j(\lambda^\ast) \ne 0 \} + 1.
\]
Above, $\lambda^\ast$ is the first time along the Lasso path that a false variable is about to be selected:
\[
\lambda^\ast = \sup\{\lambda: \text{there exists } 1 \le i \le p \text{ such that } \widehat\beta_i(\lambda) \ne 0, \beta_i = 0 \},
\]
and $\lambda^\ast - 0$ informally represents a value that is infinitesimally smaller than $\lambda^\ast$. In words, $T$ is equal to one plus the number of true variables before the first false variable. 

The problem of the rank of the first false selection has been considered by \citet{su2018first} in the case where all nonzero regression coefficients are equal. This corresponds to minimal effect size heterogeneity. While we continue employing the working hypotheses as earlier, in this section the regression coefficients $\bbeta$ are assumed to be deterministic.

\begin{proposition}\citep[Theorem 2]{su2018first}\label{prop:logT}
Under the working hypotheses, let $\beta_j = M$ for $1 \le j \le k$ and $\beta_j = 0$ for $k+1 \le j \le p$, where $k/p \goto \epsilon$ and $M \goto \infty$ as $n, p \goto \infty$. Then, the rank $T$ of the first false variable selected by the Lasso satisfies
\[
\log T = (1 + o_{\P}(1)) \sqrt{\frac{2\delta\log p}{\epsilon}},
\]
where $o_{\P}(1)$ tends to 0 in probability.
\end{proposition}
This result also applies to forward stepwise regression and least angle regression \citep{efron2004least}. Note that this proposition considers the regime where the signal-to-noise ratio $M/\sigma \goto \infty$ as $\sigma$ is fixed. If $M/\sigma$ is bounded, one has $\log T \le (1 + o_{\P}(1)) \sqrt{\frac{2\delta\log p}{\epsilon}}$~\citep[Theorem 1]{su2018first}. Indeed, the original theorem predicts that 
\[
\log T\le (1+o_{\P}(1)) \left( \sqrt{2n(\log p)/k} - n/(2k) + \log(n/(2p\log p)) \right),
\] 
which is reduced to the upper bound above since $n/p \goto \delta$ and $k/p \goto \epsilon$ under our working hypotheses.

Turning to most heterogeneous effect sizes, we have the result below.
\begin{proposition}\label{prop:wainwright}
Under the working hypotheses, let $\beta_j = M^{k+1-j}$ for $1 \le j \le k$ and $\beta_j = 0$ for $k+1 \le j \le p$, where $k/p \goto \epsilon$. If $M$ is sufficiently large, then there exists $\lambda$ depending on $n$ such that
\[
\#\left\{ j: \widehat\beta_j(\lambda) \ne 0, \beta_j \ne 0 \right\} = (1 + o_{\P}(1))\frac{n}{2\log p}, \text{ and } \#\left\{ j: \widehat\beta_j(\lambda) \ne 0, \beta_j = 0 \right\} = 0
\]
as $n, p \goto \infty$.
\end{proposition}

The proof of this proposition is given in the appendix. Regarding how large $M$ should be, precisely, this proposition holds if $M$ satisfies $M \ge n^a$ as $n \goto \infty$ for any constant $a > \frac{1}{2}$. It is also worth mentioning that the proof is adapted from the proof of Theorem 1 in \citet{wainwright2009information}. The effect sizes in Proposition~\ref{prop:wainwright} are essentially the same as an $(\epsilon, m, M)$-heterogeneous prior \eqref{eq:eps_m_M_prior} with $m \goto \infty$.

Proposition~\ref{prop:wainwright} asserts that all $(1 + o_{\P}(1))\frac{n}{2\log p}$ selected variables are true at some point along the Lasso path. If the Lasso does not kick out any selected variables before that point,\footnote{It is well-known that along the Lasso path, a selected variable can be dropped out as $\lambda$ decreases~\citep{efron2004least}. However, we did not observe this phenomenon before the first $(1 + o_{\P}(1))\frac{n}{2\log p}$ variables are selected in all of our simulations.} this result implies that $T \ge (1 + o_{\P}(1))\frac{n}{2\log p}$. Recognizing the fact that
\[
\e^{(1 + o_{\P}(1)) \sqrt{\frac{2\delta\log p}{\epsilon}}} \ll  (1 + o_{\P}(1)) \frac{n}{2\log p},
\]
most heterogeneous effect sizes are more \textit{favorable} than least heterogeneous effect sizes for the Lasso not only in terms of the TPP--FDP trade-off as shown in the previous section, but also in terms of the rank of the first false variable.

Unlike \Cref{thm:lower} and \Cref{thm:upper}, the two propositions here are silent on whether their bounds can be extended to general $\epsilon p$-sparse effect sizes. The following theorem gives a partial \textit{affirmative} answer to this question, which broadly applies to all regression coefficients with sparsity no more than $\epsilon p$, as opposed to the exact sparsity level $\epsilon p$. 

\begin{theorem}\label{prop:rank_opt}
Under the working hypotheses, for arbitrary regression coefficients $\bbeta$ with sparsity satisfying $k \le \epsilon p$, the rank $T$ of the first false variable selected by the Lasso satisfies
\[
T \le (1 + o_{\P}(1)) \frac{n}{2\log p}
\]
as $n, p \goto \infty$.

\end{theorem}

Together with Proposition~\ref{prop:wainwright}, this theorem indicates that maximal effect size heterogeneity is most favorable for the Lasso in terms of the rank of the first false variable. Importantly, the sharp bound $(1 + o_{\P}(1)) \frac{n}{2\log p}$ is the maximum number of true variables before a false selection for essentially all sparsity levels, no matter how strong and how heterogeneous the effect sizes are. This novel result is a contribution of independent interest to high-dimensional statistics. The proof is given in \Cref{sec:matching-bound} and does not involve any elements from AMP theory.

In regard to Proposition~\ref{prop:logT}, however, it is tempting to ask whether minimal effect size heterogeneity is least favorable from the same standpoint; that is, whether or not
\[
\log T \ge (1 + o_{\P}(1)) \sqrt{\frac{2\delta\log p}{\epsilon}}
\]
for any $\epsilon p$-sparse regression coefficients in the noiseless case. We leave this question for future work.

In passing, we briefly explain how and why effect size heterogeneity has a significant impact on model selection by the Lasso, shedding light on the price-of-competition phenomenon. To ease the elaboration, we assume the noiseless setting ($\bz = \bm{0}$) and denote by $S$ the set of all true variables. Consider the Lasso solution $\widehat\bbeta(\lambda)$ at some $\lambda$ where no false selection occurs (that is, the support $\widehat S$ of $\widehat\bbeta$ is a subset of $S$). Our explanation relies crucially on the fact that a variable $\bX_j$ ($j \notin \widehat S$) is likely to be the next selected variable if its inner product with the residual, $\bX_j^\top (\by - \bX \widehat\bbeta)$, is the largest in magnitude. Note that (denote by $\bX_Q$ the matrix that is formed by the columns corresponding to $Q$ for a subset $Q$ of $\{1, \ldots, p\}$)
\[
\begin{aligned}
\bX_j^\top (\by - \bX \widehat\bbeta) = \bX_j^\top (\by - \bX_{\widehat S} \widehat\bbeta_{\widehat S}) = \bX_j^\top \bX_{S\setminus \widehat S} \bbeta_{S\setminus \widehat S} + \bX_j^\top \bX_{\widehat S} (\bbeta_{\widehat S} - \widehat\bbeta_{\widehat S}).
\end{aligned}
\]

Now, we argue that the largest $\bX_j^\top \bX_{S\setminus \widehat S} \bbeta_{S\setminus \widehat S}$ in absolute value in the case of high effect size heterogeneity is likely to be from a true variable $\bX_j$ ($j \in S\setminus \widehat S$) and, conversely, it is likely to be from an irrelevant variable $\bX_j$ ($j \notin S$) if effect size heterogeneity is low. Informally, regarding $\widehat S$ as deterministic, then $\bX_j^\top \bX_{S\setminus \widehat S} \bbeta_{S\setminus \widehat S}$ is
approximately normally distributed with variance $\|\bX_{S\setminus \widehat S} \bbeta_{S\setminus \widehat S}\|^2/n \approx \|\bbeta_{S\setminus \widehat S}\|^2/n$ and mean
\[
\begin{cases}
0 &\quad \text{ if } j \notin  S\\
\beta_j &\quad \text{ if } j \in S\setminus \widehat S.
\end{cases}
\]
In the setting of Proposition~\ref{prop:wainwright} where true effect sizes are widely different from each other, the standard deviation $\|\bbeta_{S\setminus \widehat S}\|/\sqrt{n}$ is much smaller than $\sup_{j \in S\setminus \widehat S} \beta_j$. Consequently, the unselected variable with the largest effect size $\sup_{j \in S\setminus \widehat S} \beta_j$ tends to stand out, with essentially no ``competition'' among all unselected variables, thereby being the next selected variable. In the setting of Proposition~\ref{prop:logT}, however, the standard deviation $\|\bbeta_{S\setminus \widehat S}\|/\sqrt{n}$ is comparable to the largest unselected effect sizes, which are in fact of the same size. Another way to put this is that the overall effect is evenly distributed across true variables, and the resulted competition renders any variable dwarfed by the noise. Accordingly, some noise variable $\bX_j$ is very likely to have a larger inner product $\bX_j^\top (\by - \bX \widehat\bbeta)$ in
magnitude than any unselected true variable does. As such, a false selection is likely to occur very early when effect size heterogeneity is low.

\subsection{Proof of \Cref{prop:rank_opt}}
\label{sec:matching-bound}

Let $\nu > 0$ be any small constant and denote by $\A_{\nu}$ the event that the rank of the first false variable
\[
T \ge (1 + \nu) \frac{n}{2\log p}.
\]
The proof follows if one can show $\P(\A_{\nu}) \goto 0$ for all $\nu > 0$ as $n, p \goto \infty$. Recall that $S$ denotes the support $\supp(\bbeta)$. If the sparsity $|S| = k < (1 + \nu )\frac{n}{2\log p} -  1 =  (1 + \nu + o(1) )\frac{n}{2\log p}$, the event $\A_{\nu}$ is an empty set because $T$ is always no greater than $|S| + 1 < (1 + \nu )\frac{n}{2\log p}$, leading to $\P(\A_{\nu}) = 0$.

Now, we turn to the more challenging case where $k \ge (1 + \nu )\frac{n}{2\log p} -  1$ and the remainder of the proof aims to show $\P(\A_{\nu}) \goto 0$. Consider the solution $\widetilde{\bbeta}(\lambda)$ to the restricted Lasso problem
\begin{equation}\label{eq:res_Lasso}
\widetilde{\bbeta}(\lambda) := \argmin_{\bm{b} \in \R^k} \frac{1}{2}\|\bm{y} - \bX_S \bb \|^2+\lambda \|\bm{b}\|_1.
\end{equation}
Let
\[
\overline\lambda = \sup\left\{ \lambda: \|\widetilde{\bbeta}(\lambda)\|_0 \ge (1 + \nu )\frac{n}{2\log p} - 1 \right\}
\]
be the first time that the restricted Lasso selects $(1 + \nu )\frac{n}{2\log p} - 1$ variables and denote by $\widehat S$ the support of $\widetilde{\bbeta}(\overline\lambda - 0)$ (here $\overline\lambda - 0$ is infinitesimally smaller than $\overline\lambda$). In particular, this set must satisfy
\begin{equation}\label{eq:s_size}
(1 + \nu )\frac{n}{2\log p} - 1 \le |\widehat S| \le (1 + \nu )\frac{n}{2\log p}.
\end{equation}
On the event $\A_{\nu}$, the support of the full Lasso solution is a subset of $S$. Therefore, $\widetilde\bbeta(\overline\lambda)$ defined in \eqref{eq:res_Lasso} is also the solution to the full Lasso problem at $\overline\lambda$:
\[
\widetilde{\bbeta}(\overline\lambda) = \argmin_{\bm{b} \in \R^p} \frac{1}{2}\|\bm{y} - \bX \bb \|^2 + \overline\lambda \|\bm{b}\|_1.
\]
Note that $\widetilde{\bbeta}(\overline\lambda)$ may be $k$-dimensional as in \eqref{eq:res_Lasso} or $p$-dimensional by setting the remaining $p-k$ entries to zero, depending on the context. Writing $\widetilde{\bbeta}$ as a shorthand for $\widetilde{\bbeta}(\overline\lambda)$, as a consequence, we have $\left| \bX_j^\top (\by - \bX_S \widetilde{\bbeta}) \right| \le \overline\lambda$
for all $j \notin S$ on the event $\A_{\nu}$, thereby certifying
\begin{equation}\nonumber
\P(\A_{\nu}) \le \P\left(\left| \bX_j^\top (\by - \bX_S \widetilde{\bbeta}) \right| \le \overline\lambda \text{ for all } j \notin S \right).
\end{equation}

To prove $\P(\A_{\nu}) \goto 0$, therefore, it suffices to show that 
\begin{equation}\label{eq:one_suffice}
\max_{j \notin S}\left| \bX_j^\top (\by - \bX_S \widetilde{\bbeta}) \right| > \overline\lambda
\end{equation}
with probability tending to one. Making use of the independence between $\bX_j$ and $\by - \bX_S \widetilde{\bbeta}$, $\bX_j^\top (\by - \bX_S \widetilde{\bbeta})$'s are $p - k$ i.i.d.~normal random variables with mean 0 and variance $\|\by - \bX_S \widetilde{\bbeta}\|^2/n$, conditional on $\by - \bX_S \widetilde{\bbeta}$. This gives
\begin{equation}\label{eq:bound_proj}
\begin{aligned}
\max_{j \notin S}\left| \bX_j^\top (\by - \bX_S \widetilde{\bbeta}) \right| &= (1 + o_{\P}(1)) \frac{\|\by - \bX_S \widetilde{\bbeta}\|}{\sqrt n} \sqrt{2\log (p-k)}\\
& \ge  (1 + o_{\P}(1)) \frac{\|\bX_{\widehat S}(\bX_{\widehat S}^\top \bX_{\widehat S})^{-1} \bX_{\widehat S}^\top (\by - \bX_S \widetilde{\bbeta})\|}{\sqrt n} \sqrt{2\log (p-k)},
\end{aligned}
\end{equation}
where the inequality follows since $\bX_{\widehat S}(\bX_{\widehat S}^\top \bX_{\widehat S})^{-1} \bX_{\widehat S}^\top$ is a projection. For the moment, take the inequality
\begin{equation}\label{eq:to_proof}
\|\bX_{\widehat S}(\bX_{\widehat S}^\top \bX_{\widehat S})^{-1} \bX_{\widehat S}^\top (\by - \bX_S \widetilde{\bbeta})\| \ge (1 + c) \overline\lambda\sqrt{\frac{n}{2\log p}}
\end{equation}
as given for some constant $c > 0$ possibly depending on $\nu$, with probability tending to one. Combining \eqref{eq:bound_proj} and \eqref{eq:to_proof} yields
\[
\begin{aligned}
\max_{j \notin S}\left| \bX_j^\top (\by - \bX_S \widetilde{\bbeta}) \right| &\ge (1 + o_{\P}(1))\sqrt{2 \log(p-k)} \times (1 + c)\overline\lambda \frac1{\sqrt{2\log p}}\\
&= (1 + c + o_{\P}(1))\overline\lambda \sqrt{\frac{\log(p-k)}{\log p}} \\
&\ge (1 + c + o_{\P}(1))\overline\lambda \sqrt{\frac{\log(p - \epsilon p)}{\log p}} \\
&= (1 + c + o_{\P}(1) )\overline\lambda \\
\end{aligned}
\]
with probability tending to one, which ensures \eqref{eq:one_suffice}. 


We proceed to complete the proof of this theorem by verifying \eqref{eq:to_proof}. The Karush--Kuhn--Tucker condition for the Lasso asserts that
\[
\bX_{\widehat S}^\top (\by - \bX_S \widetilde{\bbeta}) = \overline\lambda\sgn(\widetilde{\bbeta}_{\widehat S}) \in \overline\lambda \{1, -1\}^{|\widehat S|},
\]
from which we get
\[
\left\| \bX_{\widehat S}^\top (\by - \bX_S \widetilde{\bbeta}) \right\| = \overline\lambda \sqrt{|\widehat S|}.
\]
A classical result in random matrix theory (see Lemma~\ref{lm:rip} in the appendix) shows that the singular values of $\bX_{\widehat S}(\bX_{\widehat S}^\top \bX_{\widehat S})^{-1}$ are all bounded below by $\frac{1}{\sqrt{1 + \theta}}$ with probability $1 - 1/p^2$, where
\begin{equation}\label{eq:eta_s}
\begin{aligned}
\theta = C\sqrt{\frac{(1 + \nu )\frac{n}{2\log p} \log(p/((1 + \nu )\frac{n}{2\log p}))}{n}} \asymp \sqrt{\frac{\log \log p}{\log p}}
\end{aligned}
\end{equation}
for an absolute constant $C$. This allows us to get
\begin{equation}\label{eq:effect_norm}
\begin{aligned}
\|\bX_{\widehat S}(\bX_{\widehat S}^\top \bX_{\widehat S})^{-1} \bX_{\widehat S}^\top (\by - \bX_S \widetilde{\bbeta})\| &\ge \frac1{\sqrt{1 + \theta}}\|\bX_{\widehat S}^\top (\by - \bX_S \widetilde{\bbeta})\|\\
&= \sqrt{\frac{{|\widehat S|}}{{1 + \theta}}} \cdot \overline\lambda
\end{aligned}
\end{equation}
with probability tending to one. Recognizing that $\theta < \nu/2$ for sufficiently large $p$ and plugging \eqref{eq:eta_s} and \eqref{eq:s_size} into \eqref{eq:effect_norm}, we obtain
\[
\begin{aligned}
 \|\bX_{\widehat S}(\bX_{\widehat S}^\top \bX_{\widehat S})^{-1} \bX_{\widehat S}^\top (\by - \bX_S \widetilde{\bbeta})\| &\ge \sqrt{\frac{(1 + \nu )\frac{n}{2\log p} - 1}{1 + \nu/2}}  \cdot \overline\lambda\\
& = (1 + c) \overline\lambda\sqrt{\frac{n}{2\log p}}
\end{aligned}
\]
with probability approaching one, where $c = \sqrt{\frac{1 + \nu}{1 + \nu/2}} - 1 > 0$. This proves \eqref{eq:to_proof}, thereby completing the proof of \Cref{prop:rank_opt}.


\section{Proofs for the Lasso Crescent}
\label{sec:proofs}

To prove Theorems~\ref{thm:lower} and \ref{thm:upper}, we start by introducing AMP theory at a \textit{minimal} level. In the case of the Lasso, loosely speaking, tools from AMP theory enable the characterization of the asymptotic joint distribution of the Lasso estimate $\widehat\bbeta(\lambda)$ and the true regression coefficients $\bbeta$ under the working hypotheses~\citep{DMM09,ampdense,BM12}. The distribution is determined by several parameters that can be solved from two equations (see \eqref{basic} below). It is important to note, however, that this body of literature only allows for the analysis of the Lasso at a \textit{fixed} value of $\lambda$. As such, these tools are not directly applicable to the full Lasso path that this paper deals with.

To overcome this difficulty, we leverage a recent extension on AMP theory that allows us to work on the Lasso problem uniformly over its penalty parameter~\citep{su2017false}. Under the working hypotheses, let $\tau >0$ and $\alpha > \alpha_0$ be the unique solution to
\begin{equation}\label{basic}
\begin{aligned}\
	&\tau^2=\sigma^2+\frac{1}{\delta}\E(\eta_{\alpha\tau}(\Pi+\tau W)-\Pi)^2\\
	&\lambda=\left(1-\frac{1}{\delta}\P(|\Pi+\tau W| > \alpha\tau)\right)\alpha\tau,
\end{aligned}
\end{equation}
where $\eta_c(x) : = \sgn(x) \cdot \max\{|x| - c, 0\}$ is the soft-thresholding function, $W$ is a standard normal random variable that is independent of $\Pi$, and $\alpha_0 = 0$ if $\delta > 1$ and otherwise is the unique root of $(1+t^2)\Phi(-t) - t\phi(t) = \frac{\delta}{2}$ in $t \ge 0$. Let $\Pi^\star$ be distributed the same as $\Pi$ conditional on $\Pi \ne 0$, and define the two deterministic functions
\begin{equation}\label{eq:fdp_tpp_infty}
\begin{aligned}
	\tpp_{\lambda}(\Pi)~=&~  \PP(|\Pi^{\star}+\tau W| > \alpha\tau)\\
	\fdp_{\lambda}(\Pi)~=&~ \frac{ 2(1-\epsilon) \Phi(-\alpha)}{ 2(1-\epsilon) \Phi(-\alpha)+ \epsilon \PP(|\Pi^{\star}+\tau W| > \alpha\tau)}.
\end{aligned}
\end{equation}
Above, $\Pi^\star$ remains independent of $W$. For convenience, we use $\stackrel{\PP}{\longrightarrow}$ to denote convergence in probability. With the notation in place, now we state a lemma that our proofs rely on.

\begin{lemma}\citep[Lemma A.2]{su2017false_supp}\label{lem:fdp_tpp_fix_lambda}
Fix $0 <  c  <  C $. Under the working hypotheses, we have
\begin{equation}\label{eq:fdp}
\sup_{ c  <\lambda <  C } \left| \TPP_\lambda(\Pi) - \tpp_\lambda(\Pi) \right|  \stackrel{\PP}{\longrightarrow} 0 \quad \text{and} \quad \sup_{ c  <\lambda <  C } \left| \FDP_\lambda(\Pi) - \fdp_\lambda(\Pi) \right|  \stackrel{\PP}{\longrightarrow} 0.
\end{equation}
\end{lemma}

Lemma~\ref{lem:fdp_tpp_fix_lambda} offers all the elements the present paper needs from AMP theory. In addition to the use of this lemma, notably, our proofs of Theorems~\ref{thm:lower} and \ref{thm:upper} involve several technical novelties that we shall highlight in Sections~\ref{sec:q_upper} and \ref{sec:lower_curve}. Relating to the literature, the convergence of $\TPP_\lambda(\Pi)$ and $\FDP_\lambda(\Pi)$ for a single $\lambda$ has been established earlier in \citet{BM12,supp}. 

We use $q^\Pi(\cdot)$ to represent the $\lambda$-parameterized curve $(\tpp_\lambda,\fdp_\lambda)$ in the sense that
\[
\fdp_\lambda(\Pi) = q^\Pi(\tpp_\lambda(\Pi)).
\]
Formally, Lemma~\ref{lem:fundamental_tppfdp} in \Cref{sec:fundamental_prop_proof} demonstrates that the instance-specific trade-off curve $q^\Pi$ is continuously differentiable and strictly increasing. Relating to \Cref{sec:lasso-crescent}, Corollary~\ref{cor:lower} implies that, taking the $(\epsilon, m, M)$-heterogeneous prior $\Pi^\Delta$, $q^{\Pi^\Delta}$ converges to $q^\Delta$ as $m, M \goto \infty$. Likewise, from Corollary~\ref{cor:upper} it follows that $q^{\Pi^\nabla}(\cdot)$ is identical to $q^\nabla(\cdot)$ in the noiseless setting.

\subsection{The upper boundary}
\label{sec:q_upper}

Our first aim is to prove \Cref{thm:upper} along with Proposition~\ref{prop:upper}. The proof is built on top of the following important lemma, which considers a non-constant $\Pi^\star$.

\begin{lemma}\label{tradeoff}

Let $\Pi$ be any $\epsilon$-sparse prior that is non-constant conditional on $\Pi \ne 0$. In the noiseless setting $\sigma=0$, we have
\[
q^\Pi(u) < q^\nabla(u)
\]
for all $0 < u < 1$.

\end{lemma}

Taking this lemma as given for the moment, we prove part (a) of \Cref{thm:upper}. 

\begin{proof}[Proof of \Cref{thm:upper}(a)]\label{pf:thm_upper}

We start by pointing out the following fact: there exists a constant $\upsilon > 0$ such that for all $c < \lambda, \lambda' < C$, the two inequalities
\begin{equation}\label{eq:not_hold}
\tpp_{\lambda}(\Pi^\nabla) > \tpp_{\lambda'}(\Pi) - \upsilon \text{ and } \fdp_{\lambda}(\Pi^{\nabla}) < \fdp_{\lambda'}(\Pi) + \upsilon
\end{equation}
cannot hold simultaneously.

Assuming this fact for the moment, it is a stone's throw away to prove \Cref{thm:upper}. Lemma~\ref{lem:fdp_tpp_fix_lambda} ensures that,  with probability tending to one as $n, p \goto \infty$, the four terms 
\[
\left| \TPP_\lambda(\Pi^\nabla) -  \tpp_\lambda(\Pi^\nabla) \right|, \left| \FDP_\lambda(\Pi^\nabla) -  \fdp_\lambda(\Pi^\nabla) \right|, \left| \TPP_{\lambda'}(\Pi) -  \tpp_{\lambda'}(\Pi) \right|,
\]
and $\left| \FDP_{\lambda'}(\Pi) -  \fdp_{\lambda'}(\Pi) \right|$ are all smaller than $\upsilon/2$ for all $c < \lambda, \lambda' < C$. In this event, 
\[
\TPP_{\lambda}(\Pi^\nabla) > \TPP_{\lambda'}(\Pi)
\] 
implies 
\[
\tpp_{\lambda}(\Pi^\nabla) > \tpp_{\lambda'}(\Pi) - \upsilon
\]
and, likewise, $\FDP_{\lambda}(\Pi^{\nabla}) < \FDP_{\lambda'}(\Pi)$ implies $\fdp_{\lambda}(\Pi^{\nabla}) < \fdp_{\lambda'}(\Pi) + \upsilon$. Recognizing that the two inequalities in \eqref{eq:not_hold} cannot both hold, therefore, in this event the following inequalities
\[
\TPP_{\lambda}(\Pi^\nabla) > \TPP_{\lambda'}(\Pi) \text{ and } \FDP_{\lambda}(\Pi^{\nabla}) < \FDP_{\lambda'}(\Pi)
\]
cannot hold simultaneously for all $c < \lambda, \lambda' < C$. In words, $\left( \TPP_{\lambda}(\Pi^\nabla), \FDP_{\lambda}(\Pi^{\nabla}) \right)$ does not outperform  $\left( \TPP_{\lambda'}(\Pi), \FDP_{\lambda'}(\Pi) \right)$, and this applies to all $c < \lambda, \lambda' < C$ with probability tending to one.

We conclude the proof by verifying \eqref{eq:not_hold}. To this end, first find $0 < u_1 < u_2 < 1$ such that the asymptotic powers $\tpp_{\lambda}(\Pi^\nabla), \tpp_{\lambda'}(\Pi)$ are always between $u_1$ and $u_2$ for $c < \lambda, \lambda' < C$. Next, set
\begin{equation}\label{eq:4.8}
\upsilon' := \inf_{u_1 \le u_\le u_2} \left( q^\nabla(u) - q^\Pi(u) \right).
\end{equation}
From Lemma~\ref{tradeoff}, we must have $\upsilon' > 0$. Since $q^\nabla$ is a continuous function on the closed interval $[0, 1]$, its uniform continuity yields
\begin{equation}\label{eq:4.9}
\left| q^\nabla(u) - q^\nabla(u') \right| < \frac{\upsilon'}{2}
\end{equation}
as long as $|u - u'| \le \upsilon''$ for some $\upsilon'' > 0$. 

As the final step, we show that \eqref{eq:not_hold} cannot hold simultaneously by taking $\upsilon = \min\{\upsilon'/2, \upsilon'' \}$. To see this, suppose we already have $\tpp_{\lambda}(\Pi^\nabla) > \tpp_{\lambda'}(\Pi) - \upsilon$, from which we get
\[
\begin{aligned}
\fdp_{\lambda}(\Pi^{\nabla}) &= q^\nabla(\tpp_{\lambda}(\Pi^{\nabla})) \\
&\ge q^\nabla\left( \tpp_{\lambda}(\Pi^{\nabla}) + \upsilon \right) - \frac{\upsilon'}{2}\\
& > q^\nabla\left( \tpp_{\lambda'}(\Pi)\right) - \frac{\upsilon'}{2}.
\end{aligned}
\]
Above, the first inequality follows from \eqref{eq:4.9}. We proceed by leveraging \eqref{eq:4.8} and obtain
\[
\begin{aligned}
\fdp_{\lambda}(\Pi^{\nabla}) & > q^\nabla\left( \tpp_{\lambda'}(\Pi)\right) - \frac{\upsilon'}{2} \\
& \ge q^\Pi\left( \tpp_{\lambda'}(\Pi)\right) + \upsilon' - \frac{\upsilon'}{2}\\
& = q^\Pi\left( \tpp_{\lambda'}(\Pi)\right) + \frac{\upsilon'}{2}.
\end{aligned}
\]
Finally, note that
\[
q^\Pi\left( \tpp_{\lambda'}(\Pi)\right) + \frac{\upsilon'}{2} \ge  q^\Pi\left( \tpp_{\lambda'}(\Pi)\right) + \upsilon = \fdp_{\lambda'}(\Pi) + \upsilon.
\]
Taken together, these calculations reveal that $\tpp_{\lambda}(\Pi^\nabla) > \tpp_{\lambda'}(\Pi) - \upsilon$ implies that $\fdp_{\lambda}(\Pi^{\nabla}) > \fdp_{\lambda'}(\Pi) + \upsilon$. As such, the two inequalities in \eqref{eq:not_hold} cannot hold at the same time. This completes the proof.

\end{proof}

The same reasoning in the proof above can be used to prove part (b) of \Cref{thm:upper} and Proposition~\ref{prop:upper}. More precisely, the first step is to establish the desired result for the deterministic functions $\tpp_{\lambda}$ and $\fdp_{\lambda}$ using Lemma~\ref{tradeoff}, followed by the second step that shows the uniform convergence using Lemma~\ref{lem:fdp_tpp_fix_lambda}. In particular, part (b) of \Cref{thm:upper} relies on the strictly increasing property of $q^\nabla$. Moreover, note that a lower bound on $\TPP_{\lambda'}(\Pi)$ can be translated into an upper bound on $\lambda'$~\citep[Lemma D.1]{su2017false_supp}.

Before turning to the proof of Lemma~\ref{tradeoff}, we propose the following preparatory lemma.
\begin{lemma}\citep[Lemma C.1]{su2017false}\label{concave}
For any fixed $\alpha > 0$, define a function $y = f(x)$ in the parametric form
\begin{equation}
\begin{aligned}
&x(t)=\P(|t+W>\alpha|)\\
&y(t)=\E (\eta_\alpha(t+W)-t)^2
\notag
\end{aligned}
\end{equation}
for $t\geq 0$, where $W$ is a standard normal random variable. Then $f$ is strictly concave.
\end{lemma}

\begin{proof}[Proof of Lemma \ref{tradeoff}]

We parameterize the curve $(\tpp_\lambda, \fdp_\lambda)$ using $\alpha > \alpha_0$. Explicitly, treating $\alpha$ as the free parameter instead of $\lambda$, we can solve $\tau$ from the AMP equation \eqref{basic}. Define 
	$$
	\begin{aligned}
	\fd_\alpha(\Pi) &= 2(1-\epsilon)\Phi(-\alpha) 
	\\\td_\alpha(\Pi) & =\epsilon \P(|\Pi^{\star}+\tau W|>\alpha\tau).
	\end{aligned}
	$$
This allows us to express the asymptotic power and FDP as functions of $\alpha$:	
	\begin{align*}
	\tpp_\alpha(\Pi) & =\frac{\td_\alpha(\Pi)}{\epsilon}\\
	\fdp_\alpha(\Pi) &= \frac{\fd_\alpha(\Pi)}{\fd_\alpha(\Pi)+\td_\alpha(\Pi)}.
	\end{align*}	

To prove Lemma~\ref{tradeoff}, for each $\alpha > \alpha_0$, it is sufficient to find a certain value of $M$ such that
	\begin{align}\label{eq:fd_td_compare}
	\fd_\alpha(\Pi) = \fd_\alpha(\Pi^\nabla) ~\text{ and }~
	\td_\alpha(\Pi) > \td_\alpha(\Pi^\nabla),
	\end{align}
where $\Pi^\nabla$ is the $(\epsilon, M)$-homogeneous prior \eqref{eq:least}. To see this fact, suppose on the contrary that 
\begin{equation}\label{eq:to_contrad0}
q^{\Pi}(u) \ge q^\nabla(u)
\end{equation}
for some $0 < u < 1$. Let $\alpha$ satisfy $\tpp_\alpha(\Pi) = u$. From \eqref{eq:fd_td_compare} we obtain
\begin{equation}\label{eq:to_contrad1}
u = \tpp_\alpha(\Pi)=\frac{\td_\alpha(\Pi)}{\epsilon}>\frac{\td_\alpha(\Pi^\nabla)}{\epsilon} = \tpp_\alpha(\Pi^\nabla) := u^\nabla
\end{equation}
and
\[
\fdp_\alpha(\Pi) = \frac{\fd_\alpha(\Pi)}{\fd_\alpha(\Pi)+\td_\alpha(\Pi)}<\frac{\fd_\alpha(\Pi^\nabla)}{\fd_\alpha(\Pi^\nabla)+\td_\alpha(\Pi^\nabla)} = \fdp_\alpha(\Pi^\nabla),
\]
which gives
\[
q^\Pi(u) = \fdp_\alpha(\Pi) < \fdp_\alpha(\Pi^\nabla) = q^\nabla(u^\nabla).
\]
This inequality combined with \eqref{eq:to_contrad0} yields
\[
q^\nabla(u) < q^\nabla(u^\nabla),
\]
which, together with the fact that $q^\nabla$ is an increasing function, leads to $u < u^\nabla$. This is a contradiction to \eqref{eq:to_contrad1}. Therefore, \eqref{eq:to_contrad0} cannot hold for any $0 < u < 1$.

The remainder of the proof aims to establish \eqref{eq:fd_td_compare} by constructing a certain prior $\Pi^\nabla$. Explicitly, it suffices to show
\begin{equation}\label{eq:fd_td_comparex}
\td_\alpha(\Pi) > \td_\alpha(\Pi^\nabla)
\end{equation}
because the equality in \eqref{eq:fd_td_compare} holds regardless of the choice of $\Pi^\nabla$. To construct $\Pi^\nabla$, we first write $\td_\alpha(\Pi)$ as
\begin{equation}\label{eq:ori_eq}
\td_\alpha(\Pi) = \epsilon \int \P(|t+W|>\alpha) \mathrm{d}\pi(t),
\end{equation}
where $\mathrm{d}\pi(t)$ denotes the measure of $\Pi^{\star}/\tau$. Since $\P(|t+W|>\alpha)$ is a strictly increasing function of $t$, there must exist $t' > 0$ such that
\begin{equation}\label{eq:tnabla}
\td_\alpha(\Pi) = \epsilon \P(|t' + W|>\alpha).
\end{equation}
Following \eqref{eq:least}, we let $\Pi^\nabla = t' \tau$ with probability $\epsilon$ and $\Pi^\nabla = 0$ otherwise.

Now, let $\tau^\nabla$ denote the solution to \eqref{basic} given $\alpha$ and $\Pi^{\nabla}$. That is (note that $\sigma =0$), 
\begin{equation*}
	(1-\epsilon)\E \eta_{\alpha}(W)^2 + \epsilon \E\left(\eta_{\alpha}\left(\frac{t' \tau}{\tau^\nabla}+W\right)-\frac{t' \tau}{\tau^\nabla}\right)^2 =\delta.
\end{equation*}
Our next step is to show
\[
\tau^\nabla > \tau. 
\]
To this end, we invoke Lemma~\ref{concave} and the strict concavity of $f$ gives
\begin{equation}\label{eq:concave_big}
\begin{aligned}
\E\left(\eta_{\alpha}\left(\frac{t' \tau}{\tau}+W\right) - \frac{t' \tau}{\tau}\right)^2 &\equiv f\left( \P(|t' + W|>\alpha)\right)\\
& = f\left(\int \P(|t+W|>\alpha) \mathrm{d}\pi(t)\right)\\
& > \int f\left(\P(|t+W|>\alpha)\right) \mathrm{d}\pi(t)\\
& = \int \E\left(\eta_{\alpha}\left(t + W\right) - t \right)^2 \mathrm{d}\pi(t)\\	&=\E\left(\eta_{\alpha}\left(\frac{\Pi^{\star}}{\tau}+W\right)-\frac{\Pi^{\star}}{\tau}\right)^2,
\end{aligned}
\end{equation}
where the second equality follows from the definition of $t'$ in \eqref{eq:ori_eq} and \eqref{eq:tnabla}, and the inequality is strict because $\Pi^\star$ is not constant. Together with the AMP equation for $\Pi$
\begin{equation}\nonumber
	(1-\epsilon)\E \eta_{\alpha}(W)^2 + \epsilon \E\left(\eta_{\alpha}\left(\frac{\Pi^{\star}}{\tau}+W\right)-\frac{\Pi^{\star}}{\tau}\right)^2 =\delta,
\end{equation}
\eqref{eq:concave_big} implies
\begin{equation*}
	(1-\epsilon)\E \eta_{\alpha}(W)^2 + \epsilon \E\left(\eta_{\alpha}\left(\frac{t' \tau}{\tau}+W\right)-\frac{t' \tau}{\tau}\right)^2 >\delta
\end{equation*}
or, equivalently,
\begin{equation}\label{eq:delta_lar}
(1-\epsilon)\E \eta_{\alpha}(W)^2 + \E \left[ \left(\eta_{\alpha}\left(\frac{\Pi^\nabla}{\tau}+W\right)-\frac{\Pi^\nabla}{\tau}\right)^2; \Pi^\nabla \ne 0 \right]  > \delta.
\end{equation}
By definition, however, $\tau^\nabla$ must satisfy
\begin{equation}\label{eq:hope_last}
(1-\epsilon)\E \eta_{\alpha}(W)^2 + \E \left[ \left(\eta_{\alpha}\left(\frac{\Pi^\nabla}{\tau^\nabla}+W\right)-\frac{\Pi^\nabla}{\tau^\nabla}\right)^2; \Pi^\nabla \ne 0 \right]  = \delta.
\end{equation}
A comparison between \eqref{eq:delta_lar} and \eqref{eq:hope_last} immediately gives $\tau^\nabla > \tau$.

Having shown $\tau^\nabla > \tau$, we complete the proof by noting
\begin{align*}
\td_\alpha(\Pi)  &=\epsilon\P(|t' + W|>\alpha) \\
&> \epsilon \P\left(\left|\frac{t' \tau}{\tau^\nabla}+ W\right|>\alpha\right) \\
&= \epsilon \P\left(\left| \Pi^\nabla + \tau^\nabla W \right| > \alpha\tau^\nabla \Big| \Pi^\nabla \ne 0 \right) \\
&=\td_\alpha(\Pi^\nabla).
\end{align*}
This verifies \eqref{eq:fd_td_comparex}.

\end{proof}


\subsection{The lower boundary}\label{sec:lower_curve}

Now we turn to the proof of \Cref{thm:lower}. As with the architecture of the proof of \Cref{thm:upper}, our strategy is to first prove the theorem for the deterministic functions $\tpp_{\lambda}$ and $\fdp_{\lambda}$, and then apply Lemma~\ref{lem:fdp_tpp_fix_lambda} to carry over the results to the random functions $\TPP_{\lambda}$ and $\FDP_{\lambda}$. Having said this, it is important to note that the proof presents a novel element to the literature. Below, we shall highlight the novel part of the proof of \Cref{thm:lower} and leave the rest to the appendix.

As shown in \cite{su2017false}, the trade-off curve $q^\Pi$ of any $\epsilon$-sparse prior $\Pi$ obeys
\[
q^\Pi(u) > q^\Delta(u)
\]
for $0 < u < 1$ in both the noiseless and noisy settings, where the curve $q^\Delta$ is defined in \eqref{eq:q_delta_def}. If the $(\TPP, \FDP)$ pairs from the $(\epsilon, m, M)$-heterogeneous prior $\Pi^\Delta$ form the curve $q^\Delta$ asymptotically as $n, p \goto \infty$, the proof of \Cref{thm:lower} would follow immediate, just as \Cref{thm:upper}. For any values of $m$ and $M$, however, the $\lambda$-parameterized curve $(\tpp_\lambda(\Pi^\Delta), \fdp_\lambda(\Pi^\Delta))$ does not agree with $q^\Delta$. This is in contrast to the proof of \Cref{thm:upper}, where the $(\TPP, \FDP)$ pairs from the $(\epsilon, M)$-homogeneous prior \eqref{eq:least} converge to the curve $q^\nabla$ for any value of $M \ne 0$, thanks to the assumed noiseless setting. 

To tackle this challenge, our strategy is to \textit{uniformly} approximate $q^\Delta$ using a more general prior for effect sizes that takes the form
\begin{equation}\label{eq:most_favorable_prior}
\Pi^\Delta(\bM,\bm{\gamma})=
\begin{cases}
0 & \text{w.p. $1-\epsilon$}\\
M_1 & \text{w.p. $ \epsilon\gamma_1$}\\
M_2 & \text{w.p. $ \epsilon\gamma_2$}\\
M_3 & \text{w.p. $ \epsilon\gamma_3$}\\
\cdots & \cdots\\
M_m & \text{w.p. $ \epsilon\gamma_m$},
\end{cases}
\end{equation} 
where $0 < M_1 < M_2 < \cdots < M_m$ and $\gamma_1 + \cdots + \gamma_m = 1$ with $\gamma_i > 0$. Fixing $\bm\gamma = (\gamma_1, \ldots, \gamma_m)$ while letting $M_1 \goto \infty$ and $M_{i+1}/M_{i} \to \infty$ for all $i$, we have the following lemma:

\begin{lemma}\label{lem:n_points_convergence}
The curve $q^{\Pi^\Delta(\bM,\bm{\gamma})}$ converges to a function that agrees with $q^\Delta$ at $m - 1$ points on $(0, 1)$.
\end{lemma}

\begin{figure}[!htb]
\centering
\includegraphics[width = 0.85\linewidth]{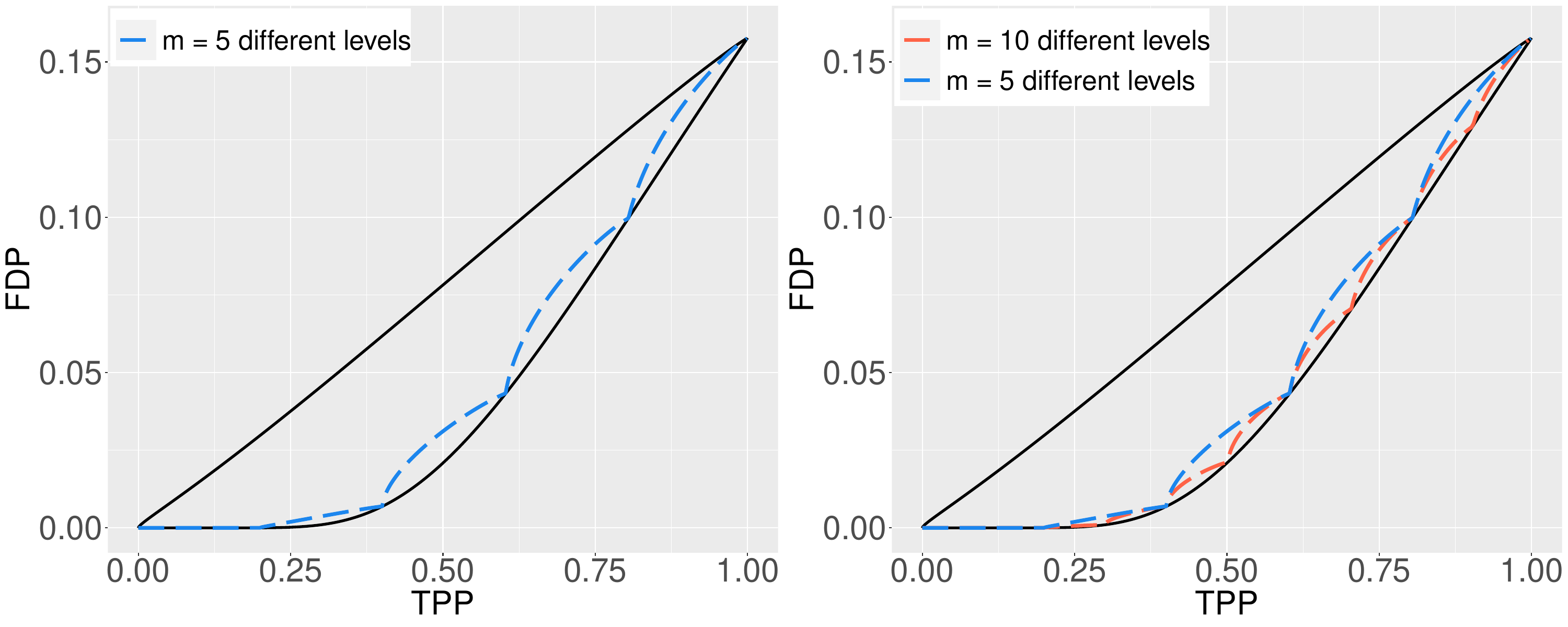}
\caption{Illustration of Lemma~\ref{lem:n_points_convergence}, showing the convergence to the lower curve $q^\Delta$. Left: $m = 5$ different levels in the prior \eqref{eq:most_favorable_prior} with $\gamma_1 = \cdots =\gamma_5 = 0.2$, and the associated trade-off curve touches the lower boundary at 4 points; Right: the case $m = 10$ and $\gamma_1 =\cdots=\gamma_{10} = 0.1$ is added as a comparison to the left case.}
\label{multiple points}
\end{figure}

For convenience, denote by $q^{\Delta(\bm\gamma)}$ the limiting curve of $q^{\Pi^\Delta(\bM,\bm{\gamma})}$ as $M_1 \goto \infty$ and $M_{i+1}/M_{i} \to \infty$. \Cref{multiple points} provides an illustration of this limiting curve. To see why $q^{\Delta(\bm\gamma)}$ is close to $q^\Delta$, note that Lemma~\ref{lem:n_points_convergence} ensures that there exist $0 < u_1 < u_2 < \cdots < u_{m-1} < 1$ such that 
\[
q^{\Delta(\bm\gamma)}(u_i) = q^{\Delta}(u_i)
\]
for $i = 1, \ldots, m-1$. In fact, the two functions also agree at $u_0 := 0$ and $u_m := 1$. Recognizing that both functions are increasing, for any $u_i \le u \le u_{i+1}$ we get
\[
0 \le q^{\Delta(\bm\gamma)}(u) - q^{\Delta}(u) \le q^{\Delta(\bm\gamma)}(u_{i+1}) - q^{\Delta}(u_i) = q^{\Delta}(u_{i+1}) - q^{\Delta}(u_i).
\]
Making use of the uniform continuity of $q^\Delta$, the desired conclusion follows if we show that the gap $u_{i+1} - u_i$ is small for all $i = 0, \ldots, m-1$. The proof of Lemma~\ref{lem:n_points_convergence}, indeed, reveals that this is true if $\max \gamma_i$ is sufficiently small. See the proof of this lemma and the remaining details in \Cref{sec:additional_proofs}.

In passing, we remark that \eqref{eq:most_favorable_prior} in the special case $m = 2$ has been considered in \cite{su2017false}. Explicitly, the lower boundary $q^\Delta$ is formed as the lower envelope of the instance-specific trade-off curves induced by the $\epsilon$-sparse priors. See the discussion following \eqref{eq:envelope_sim} in \Cref{sec:lasso-crescent}.



\section{Illustrations}
\label{sec:illustrations}

In this section, we present simulation studies to illustrate the impact of effect size heterogeneity beyond the working hypotheses, with a focus on how the impact depends on the design matrix and the noise level.

\subsection{Design matrix}

We perform four simulation studies to examine the impact of effect size heterogeneity on the Lasso method under various synthetic design matrices. Overall, the simulation results show that effect size heterogeneity remains an influential factor in determining the performance of the Lasso far beyond Gaussian designs.

\begin{figure}[htbp]
	\centering
	\includegraphics[width = 0.8\linewidth]{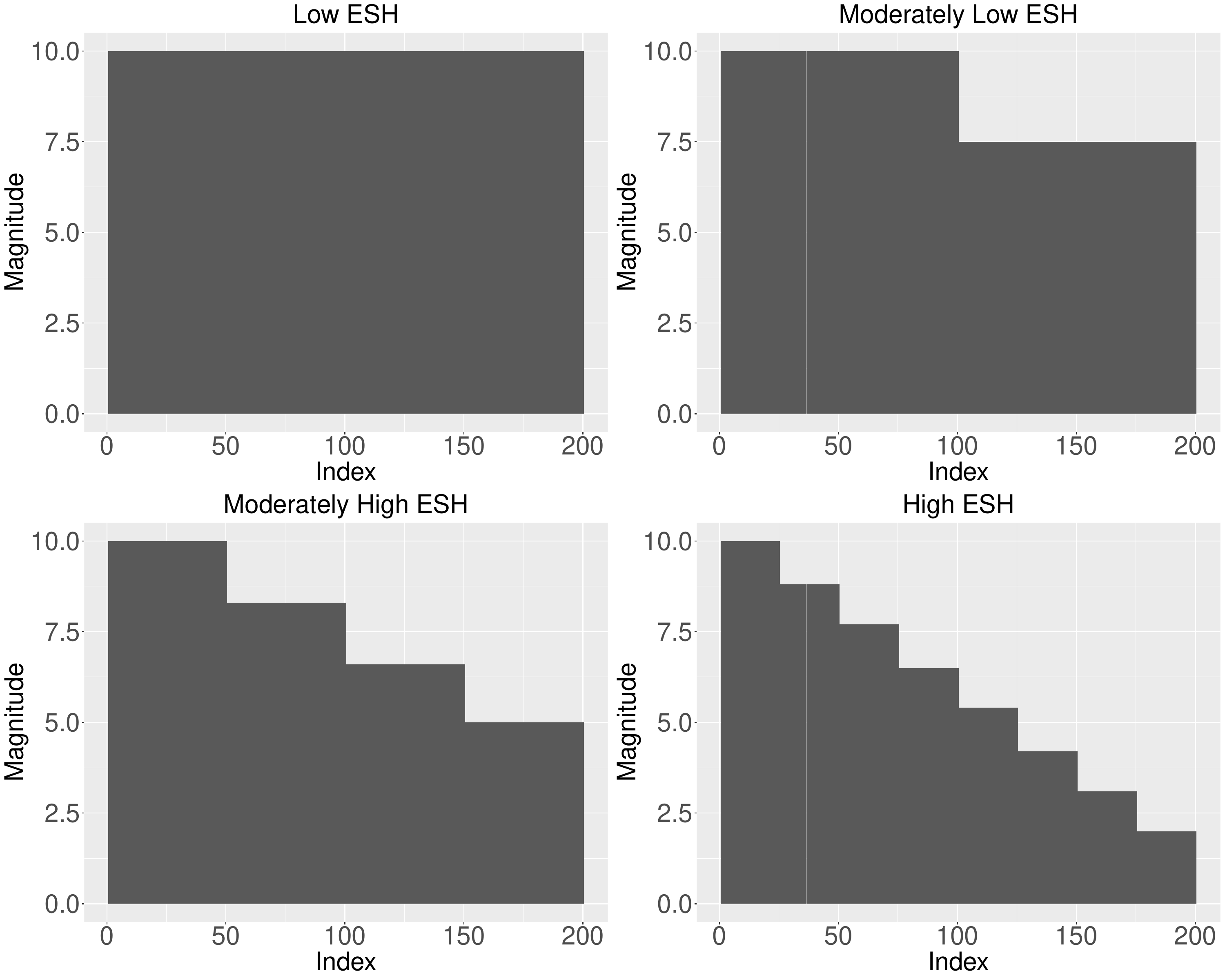}
	\caption{Four sets of effect sizes ranked in increasing order of their effect size heterogeneity. The corresponding regression coefficients in $\R^{1000}$ with sparsity 200 are used in the experiments of Figures~\ref{fig:GaussianBernoulli}, \ref{fig:varied_X_Gen}, and \ref{fig:noisy}.}
	\label{fig:config}
\end{figure}

\begin{figure}[htbp]
				\centering
			\begin{minipage}{0.42\linewidth}
				\centering
				\includegraphics[width = \linewidth]{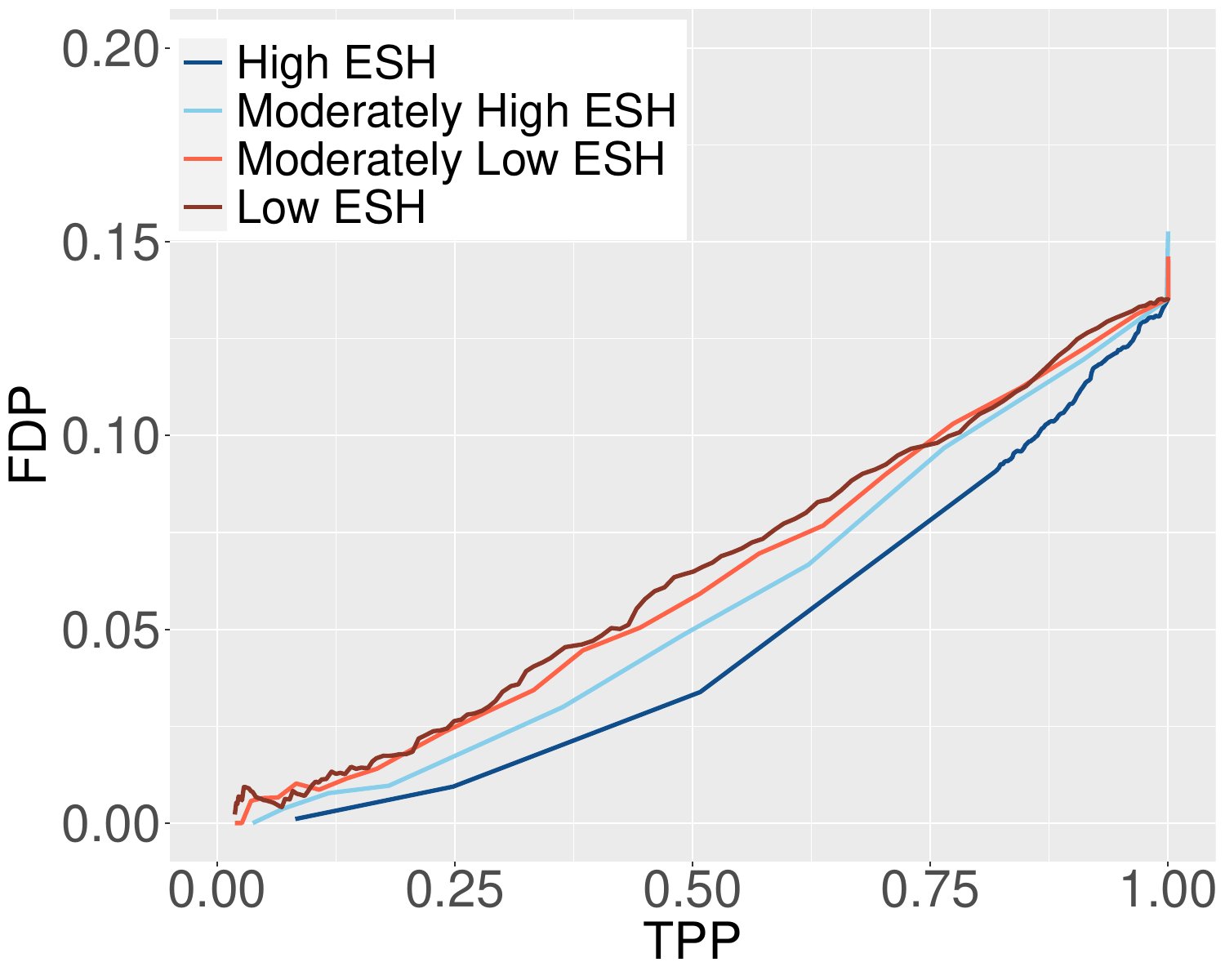}
				\label{fig:non-sparse estim}
			\end{minipage}
			\begin{minipage}{0.42\linewidth}
				\centering
				\includegraphics[width = \linewidth]{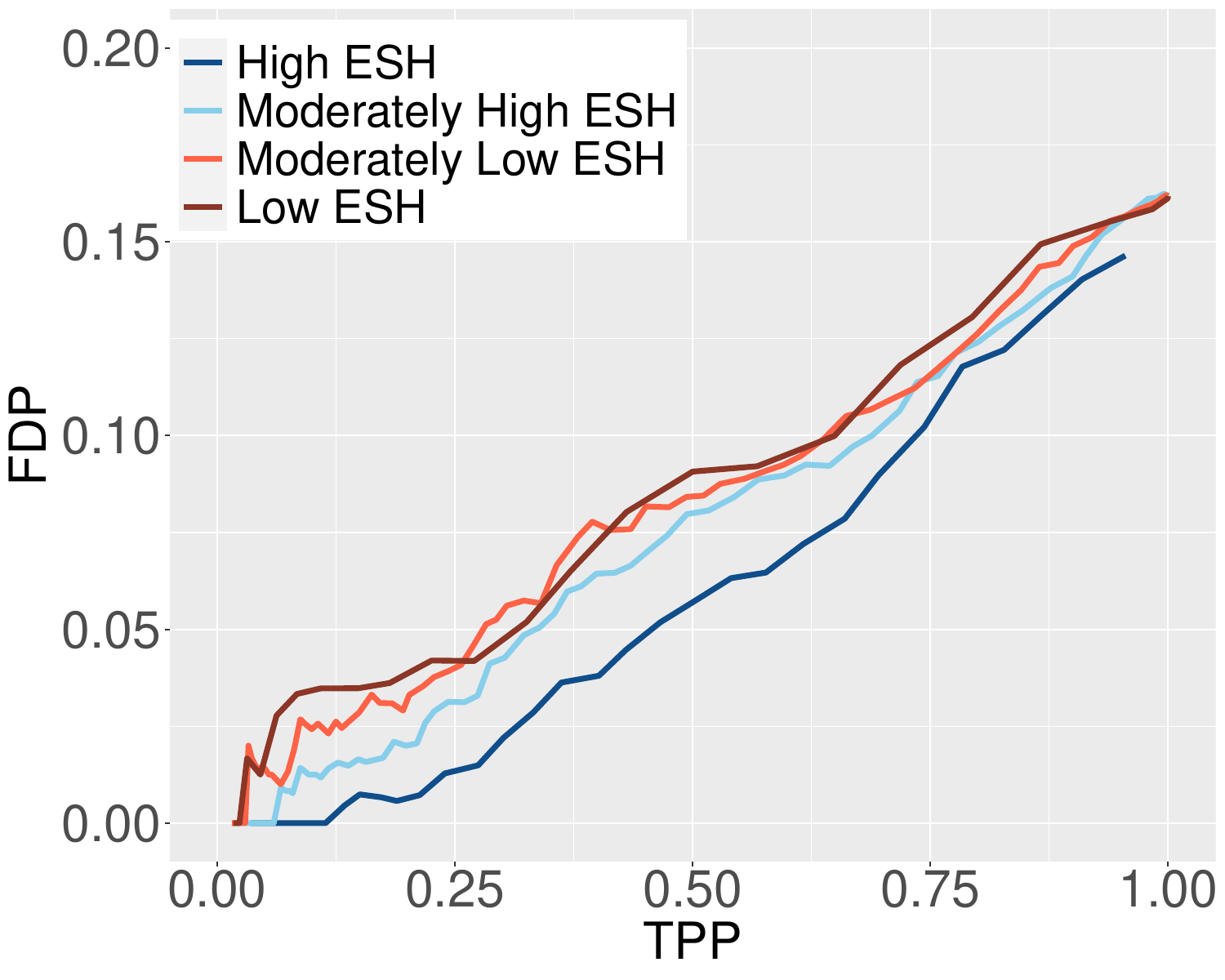}
				\label{non-sparse real}
			\end{minipage}
			\vspace{0in}
			\caption{The TPP--FDP trade-off along the Lasso path under a correlated Gaussian design and a Bernoulli design (Study 1). We set $n= p= 1000$, $k = 200$ and $\sigma = 0$ in both simulations. Left: Gaussian design matrix, each row having covariance $\bm\Sigma$ taking the form $\Sigma_{ij} = 0.5^{|i - j|}$. Right: design matrix with i.i.d.~Bernoulli entries taking the value $1/\sqrt{1000}$ or $-1/\sqrt{1000}$ with equal probability. The four sets of regression coefficients are described in \Cref{fig:config}. The mean FDP is obtained by averaging over $200$ replicates.}
			\label{fig:GaussianBernoulli}
		\end{figure}

\textit{Study 1.} We consider a design matrix of size $1000 \times 1000$ that has each row independently drawn from $\N(\bm 0, \bm \Sigma)$, where $\Sigma_{ij} = 0.5^{|i - j|}/1000$ and another design matrix of the same size that has independent Bernoulli entries, which take the value $1/\sqrt{1000}$ with probability half and otherwise $-1/\sqrt{1000}$. The sparsity is fixed to $k = 200$ while we consider four scenarios of the 200 true effects corresponding to low, moderately low, moderately high, and high effect size heterogeneity (see \Cref{fig:config}). The results on the TPP--FDP trade-off are presented in \Cref{fig:GaussianBernoulli}



\begin{figure}[htbp]
\centering
	\includegraphics[width = 0.44\linewidth]{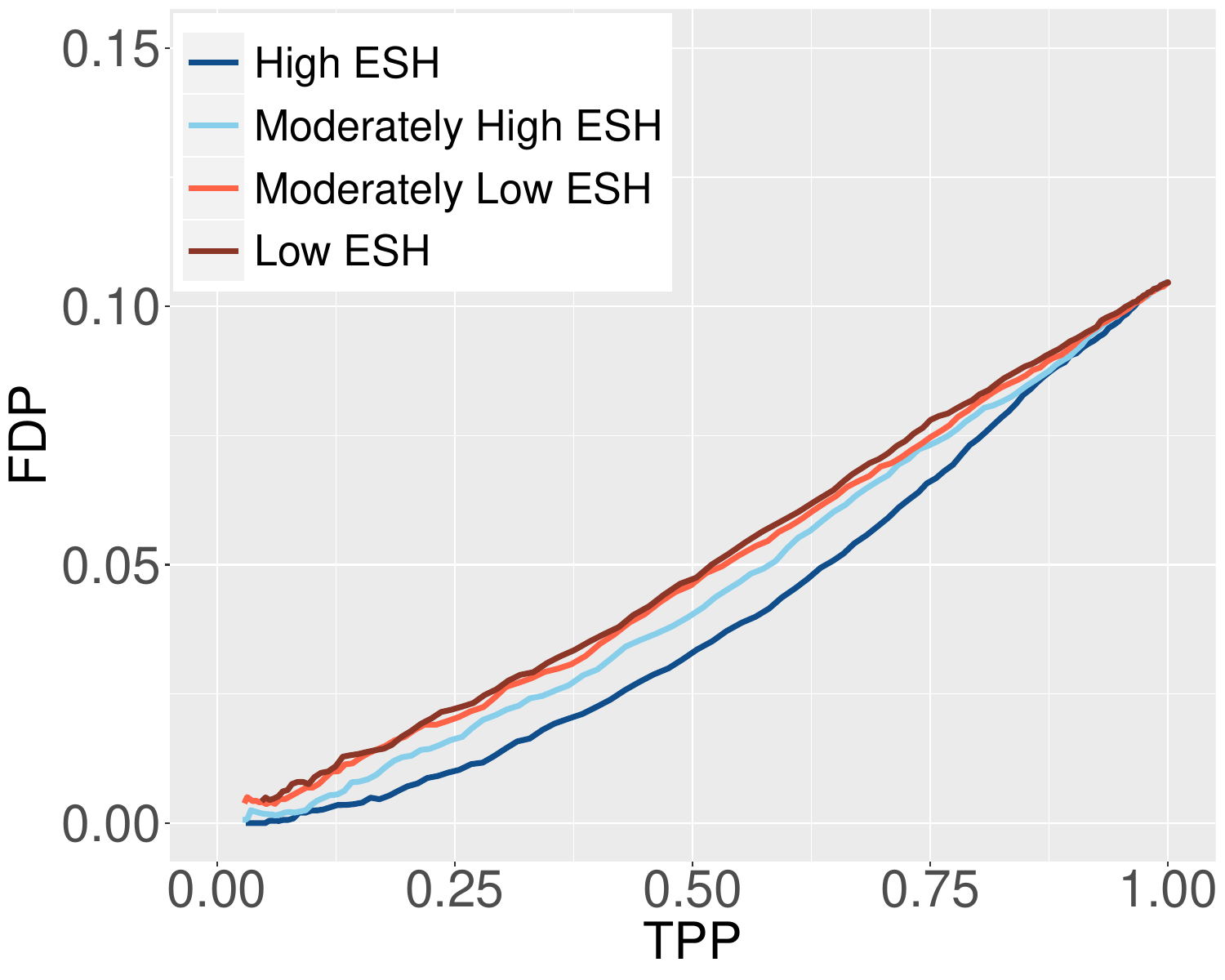}
\caption{The TPP--FDP trade-off for the genotype dataset (Study 2). The four curves correspond to the four sets of effect sizes described in Figure \ref{fig:config}. The noise term is set to be $\bm 0$. The results are obtained by averaging over $200$ replicates.}
	\label{fig:varied_X_Gen}
\end{figure}

\textit{Study 2.} In this study, we use a dataset of size $1000 \times 892$ that is simulated from the admixture of the African-American and European populations, based on the HapMap genotype data \citep{Hap} (see more details in \cite{BFST,SortedL1}). The variables can only take $0, 1,$ or $2$ according to the genotype of a genetic marker. To improve the conditioning of the design matrix, we add i.i.d.~$\N(0, 1/1000)$ perturbations to all the entries. Each column is further standardized to have mean 0 and unit norm. We use the effect sizes described in \Cref{fig:config} to generate a synthetic response $\by$ following the linear model~\eqref{eq:basic_model}, with noise $\bz = \bm{0}$. The results are plotted in Figure~\ref{fig:varied_X_Gen}.

\begin{figure}[!htp]
	\centering
	\begin{minipage}{0.42\linewidth}
		\centering
		\includegraphics[width = \linewidth]{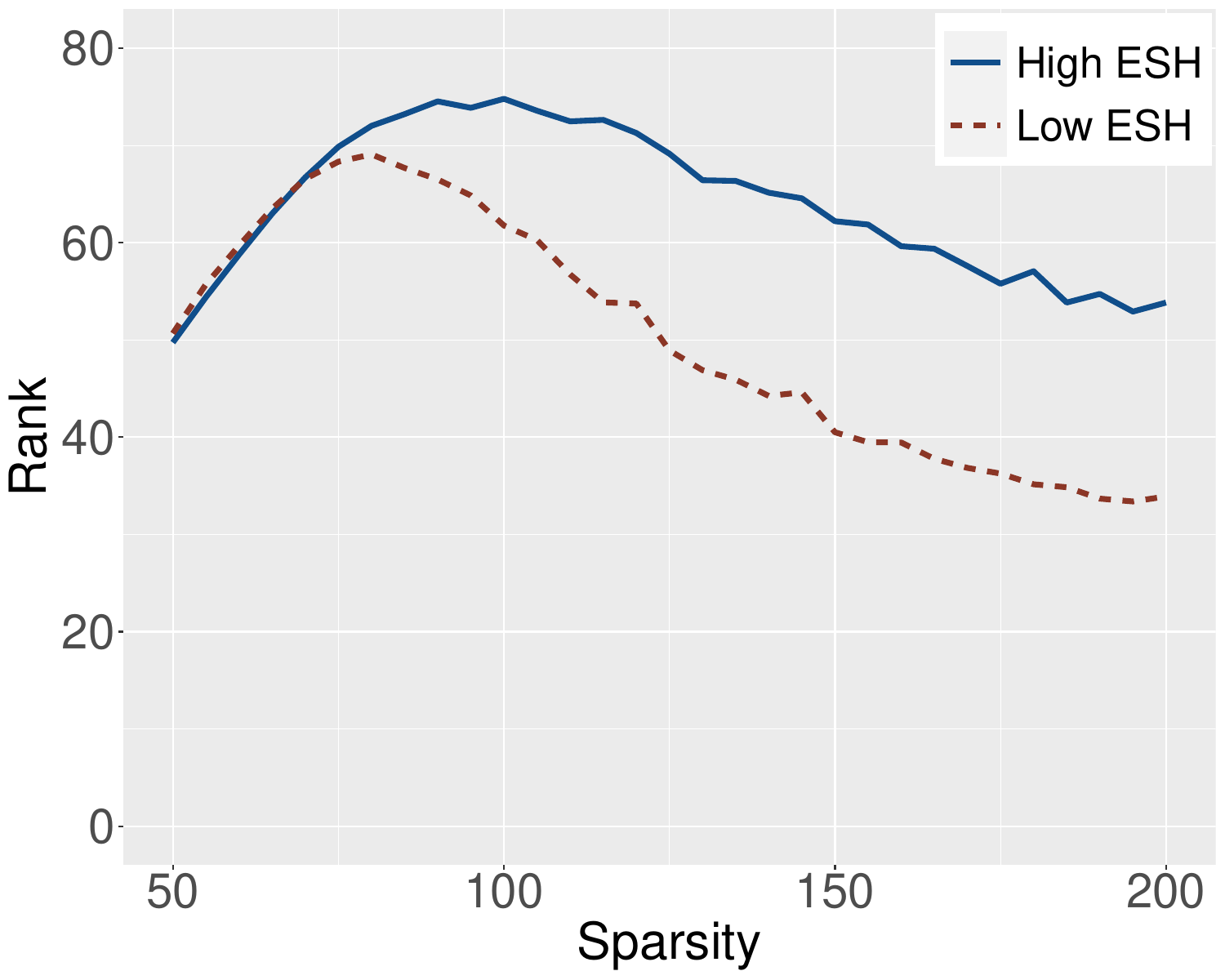}
		\label{fig:gaussian}
	\end{minipage}
	\begin{minipage}{0.42\linewidth}
		\centering
		\includegraphics[width = 1\linewidth]{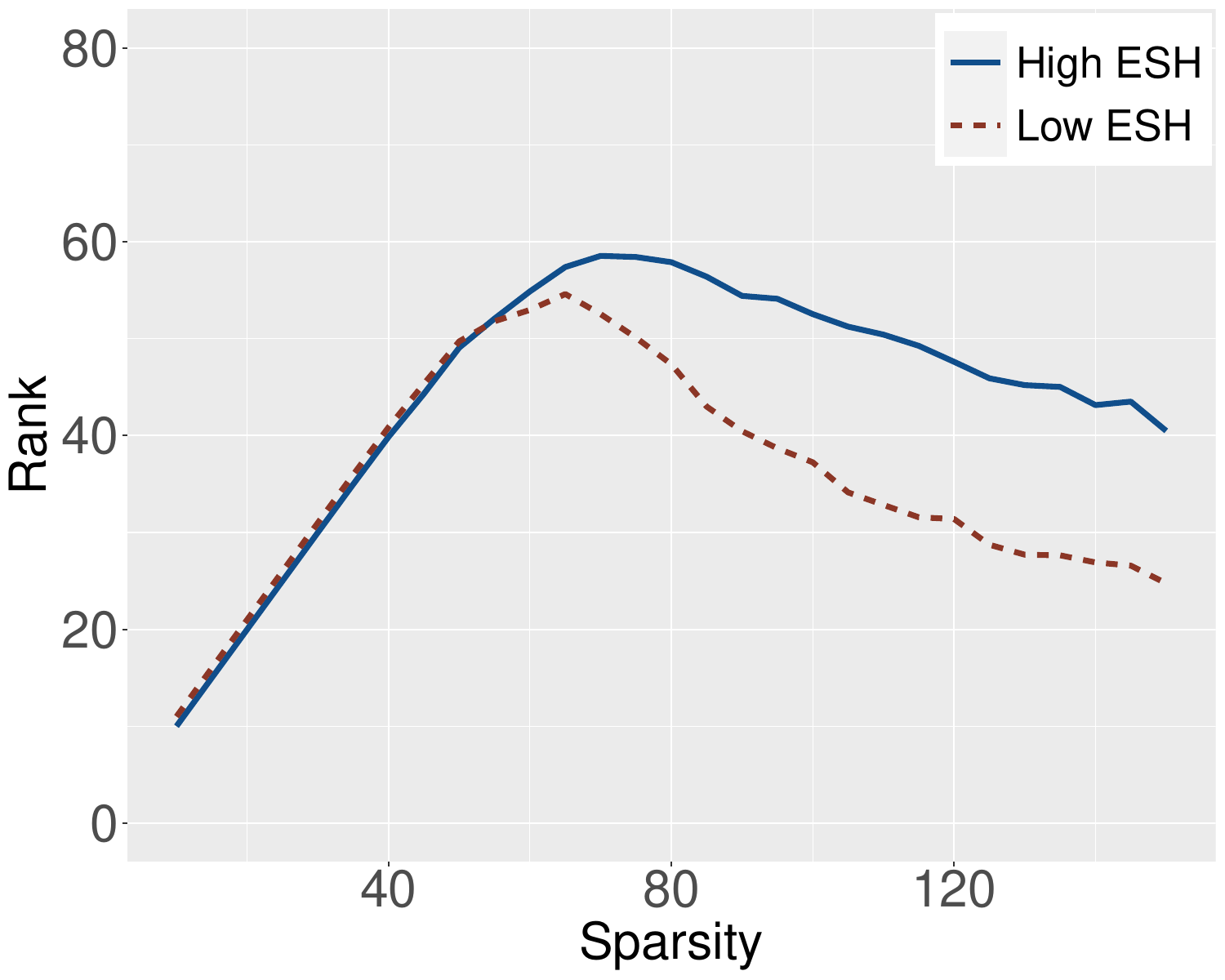}
		\label{fig:bernoulli}
	\end{minipage}
	\caption{The rank of the first spurious variable with varying sparsity (Study 3). Left: design matrix of size $1000 \times 1000$ consists of i.i.d.~$\mathcal{N}(0, \frac{1}{1000})$ entries. Right: design matrix of size $800 \times 1200$, with i.i.d.~Bernoulli entries that take the value \(1 / \sqrt{500}\) with probability \(1 / 2\) and value \(-1 / \sqrt{500}\) otherwise. Each curve is averaged over $200$ independent replicates.}
	\label{fig:rank1}
\end{figure}

\textit{Study 3.} Working under Gaussian and Bernoulli designs, we now empirically examine the rank of the first false variable. This study considers a varying sparsity level $k$ and sets the effect sizes to $\beta_j = 100$ for $j = 1, \ldots, k$ (low effect size heterogeneity) or $\beta_j = j$ for $j = 1, \ldots, k$ (high effect size heterogeneity). Each noise component \(z_{i}\) follows \(\N(0,1)\) independently. Figure~\ref{fig:rank1} shows the results under an independent Gaussian random design and an independent Bernoulli design.

\begin{figure}[htb]
	\centering
	\begin{minipage}{0.42\linewidth}
		\centering
		\includegraphics[width = 1\linewidth]{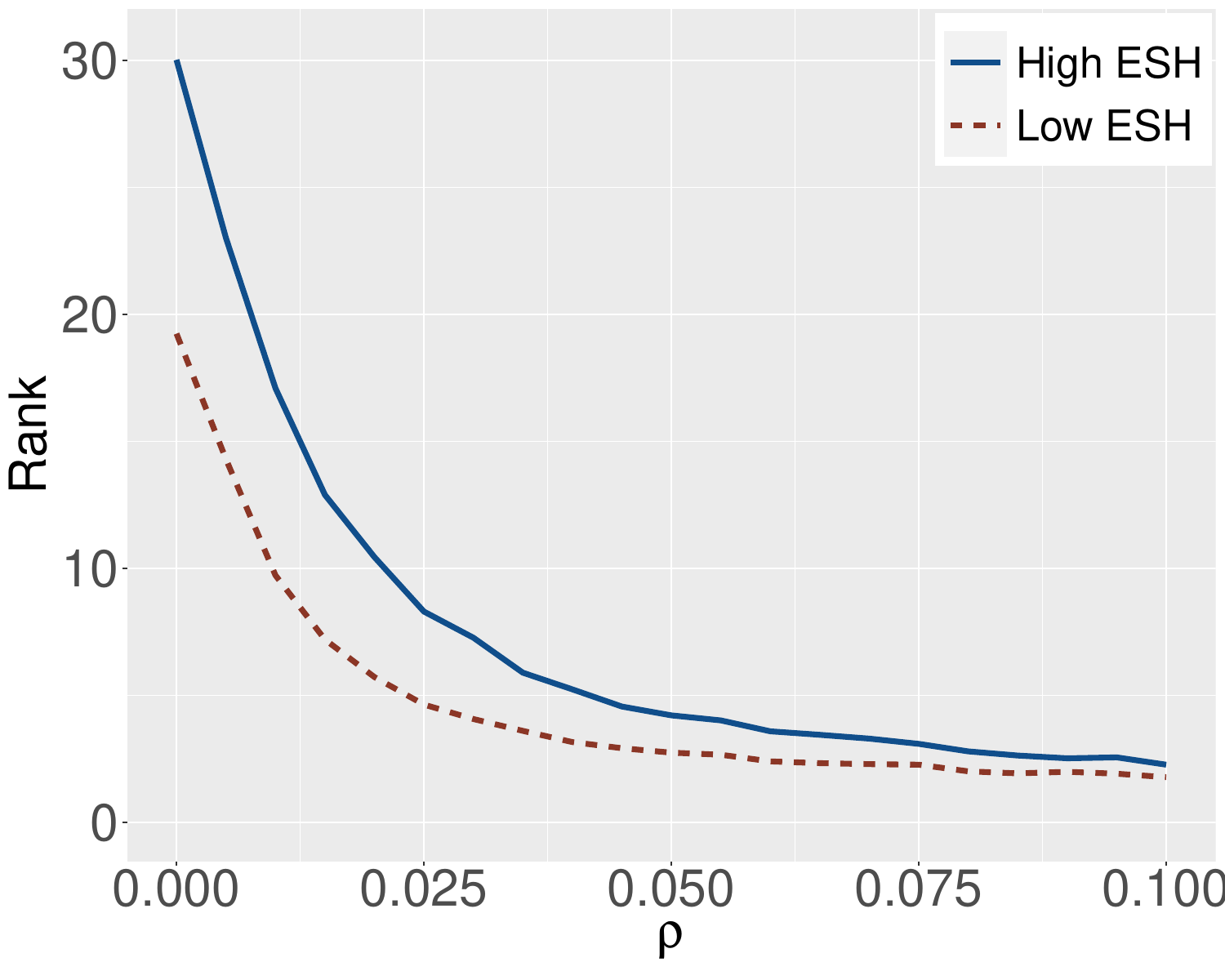}
		\label{fig:autoregressive}
	\end{minipage}
	\begin{minipage}{0.42\linewidth}
		\centering
		\includegraphics[width = 1\linewidth]{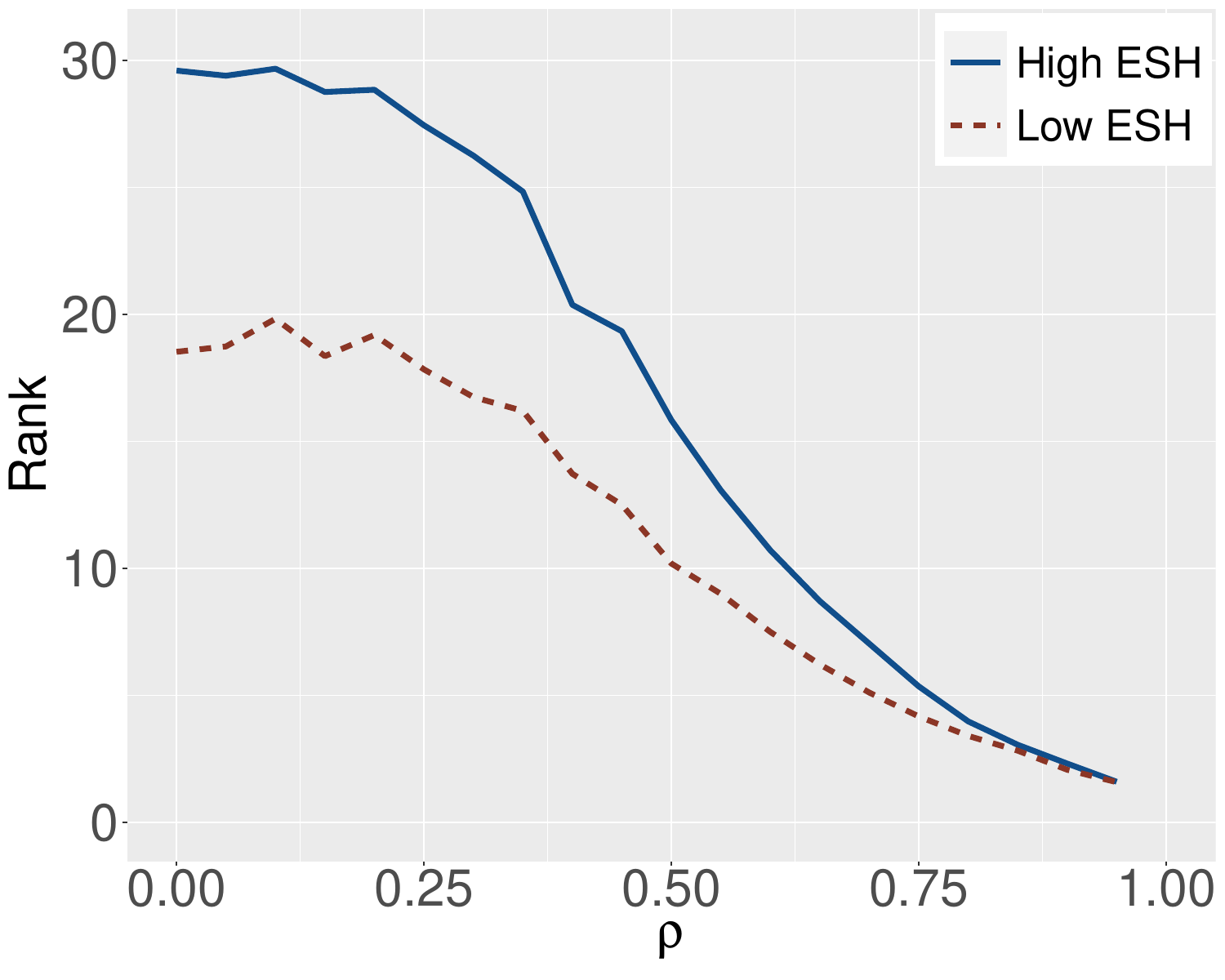}
		\label{fig:toepliz}
	\end{minipage}	
	\caption{The rank of the first spurious variable (Study 4). Left: Gaussian design with an equi-correlation covariance matrix, with the non-diagonal correlation $\rho$ varying from $0$ to $0.1$. Right: Gaussian design with covariance $\bm\Sigma$ taking the form $\Sigma_{i j}=\rho^{|i-j|} /1000$. Each curve is averaged over $200$ independent replicates.}
	\label{fig:rank2}
\end{figure}

\textit{Study 4.} This scenario uses \(500 \times 1000\) design matrices that have each row drawn independently from $\N(0, \bm{\Sigma})$. In the left panel of \Cref{fig:rank2}, the \(1000 \times 1000\) covariance matrix \(\bm{\Sigma}\) is set to \(\Sigma_{i j}=\rho / 1000\) if \(i \neq j\) and $\Sigma_{jj}=1 /1000$. In the right panel, the covariance satisfies $\Sigma_{i j}=\rho^{|i-j|} /1000$ for all $i, j$, with $\rho$ varying from 0 to $0.95$. The effect sizes are set to \(\beta_{j}= 100 \sqrt{ 2 \log p } \) for \(j \leqslant k\) (low effect size heterogeneity) or, in the low effect size heterogeneity case, the true effect sizes are set to a decreasing sequence from $ 100 \sqrt{ 2 \log p} $ to $0$. The noise \(\bz\) consists of independent standard normal variables. As is clear, both \Cref{fig:rank1} and \Cref{fig:rank2} show that the rank of the first false variable is larger when effect size heterogeneity is high, aligning with our analysis in \Cref{sec:heuristics}.

\subsection{Noise level}

\begin{figure}[htb]
	\centering
	\begin{minipage}{0.42\linewidth}
		\centering
		\includegraphics[width = 1\linewidth]{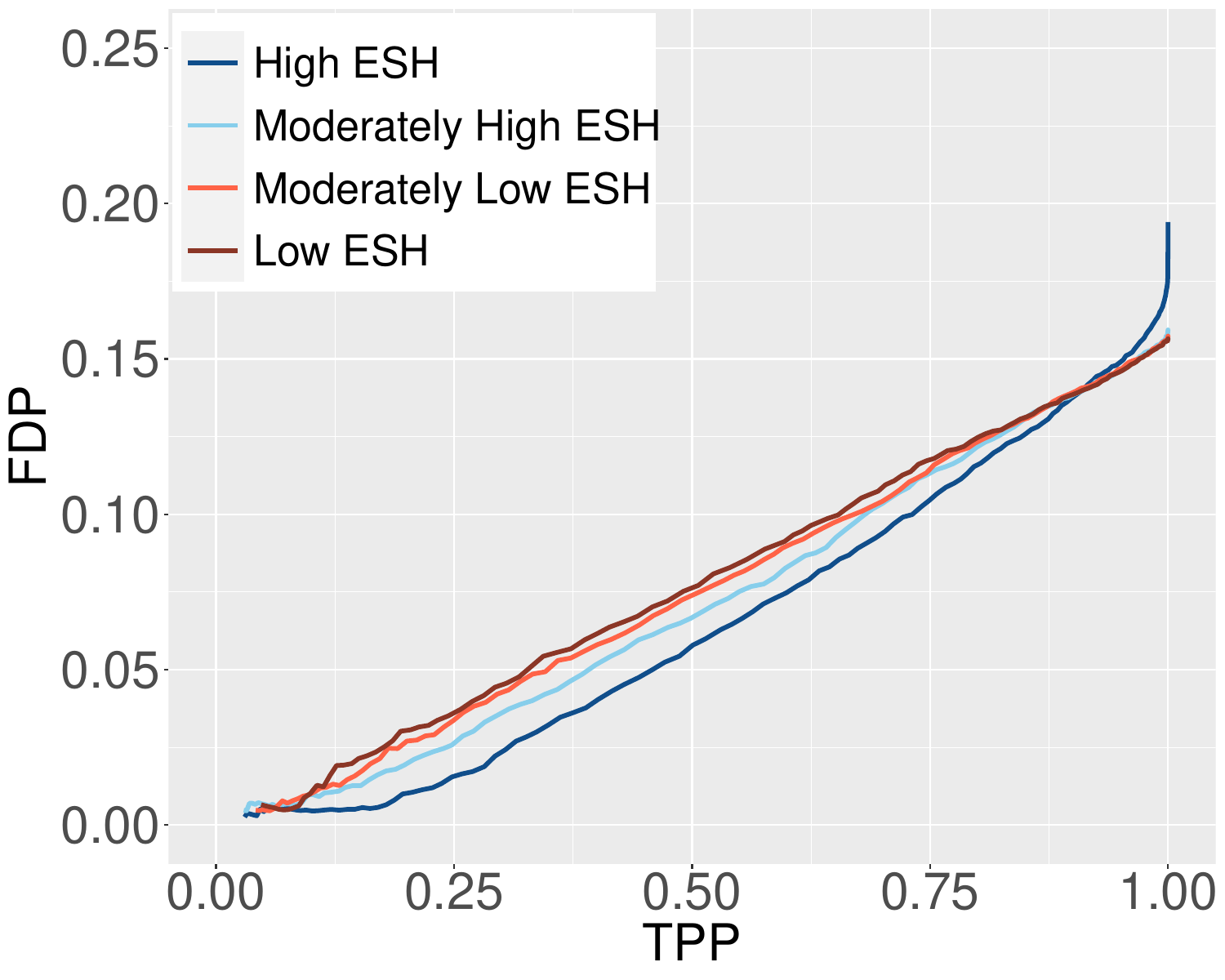}
	\end{minipage}
	\begin{minipage}{0.42\linewidth}
		\centering
		\includegraphics[width = 1\linewidth]{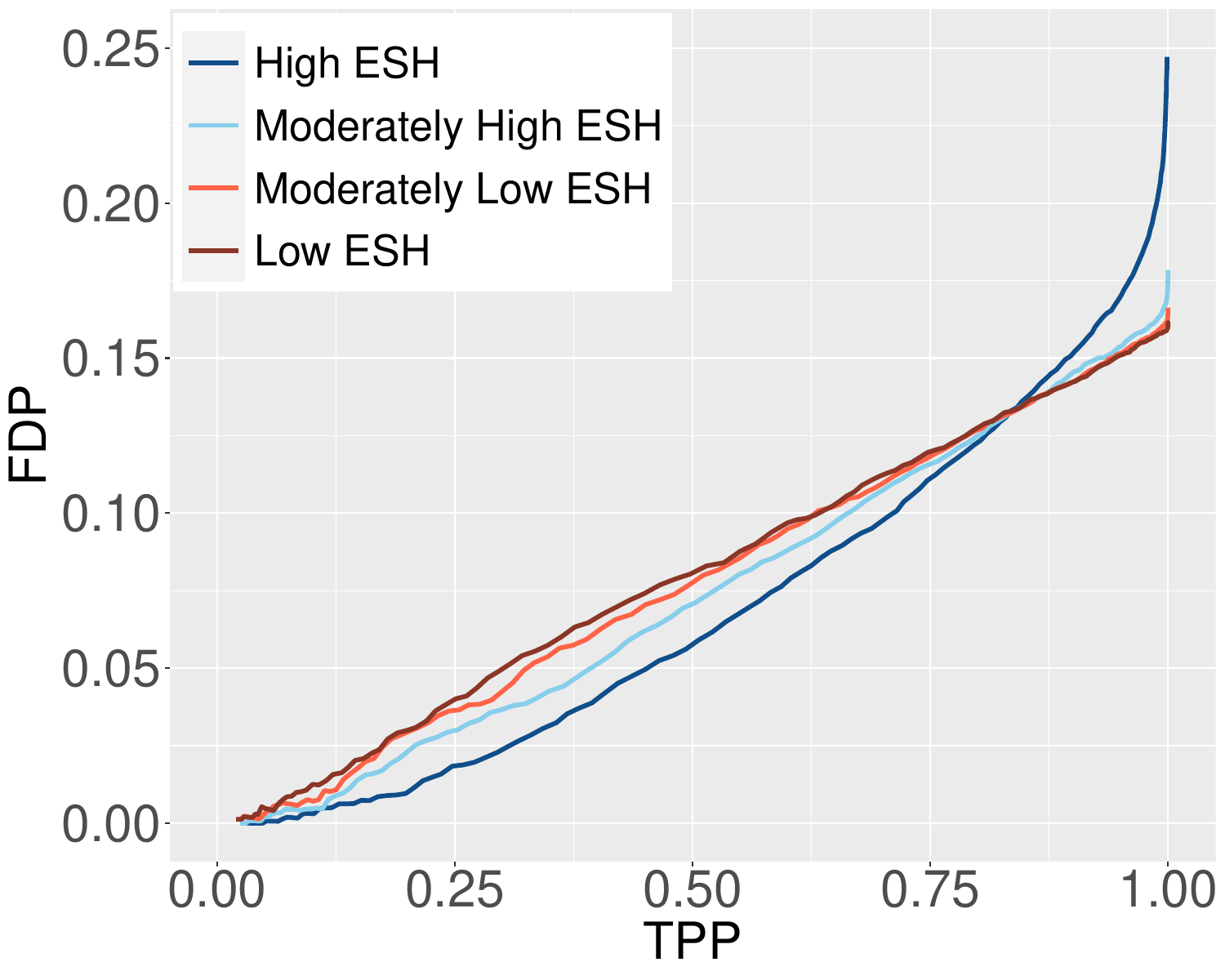}
	\end{minipage}	
\centering
\begin{minipage}{0.42\linewidth}
	\centering
	\includegraphics[width = 1\linewidth]{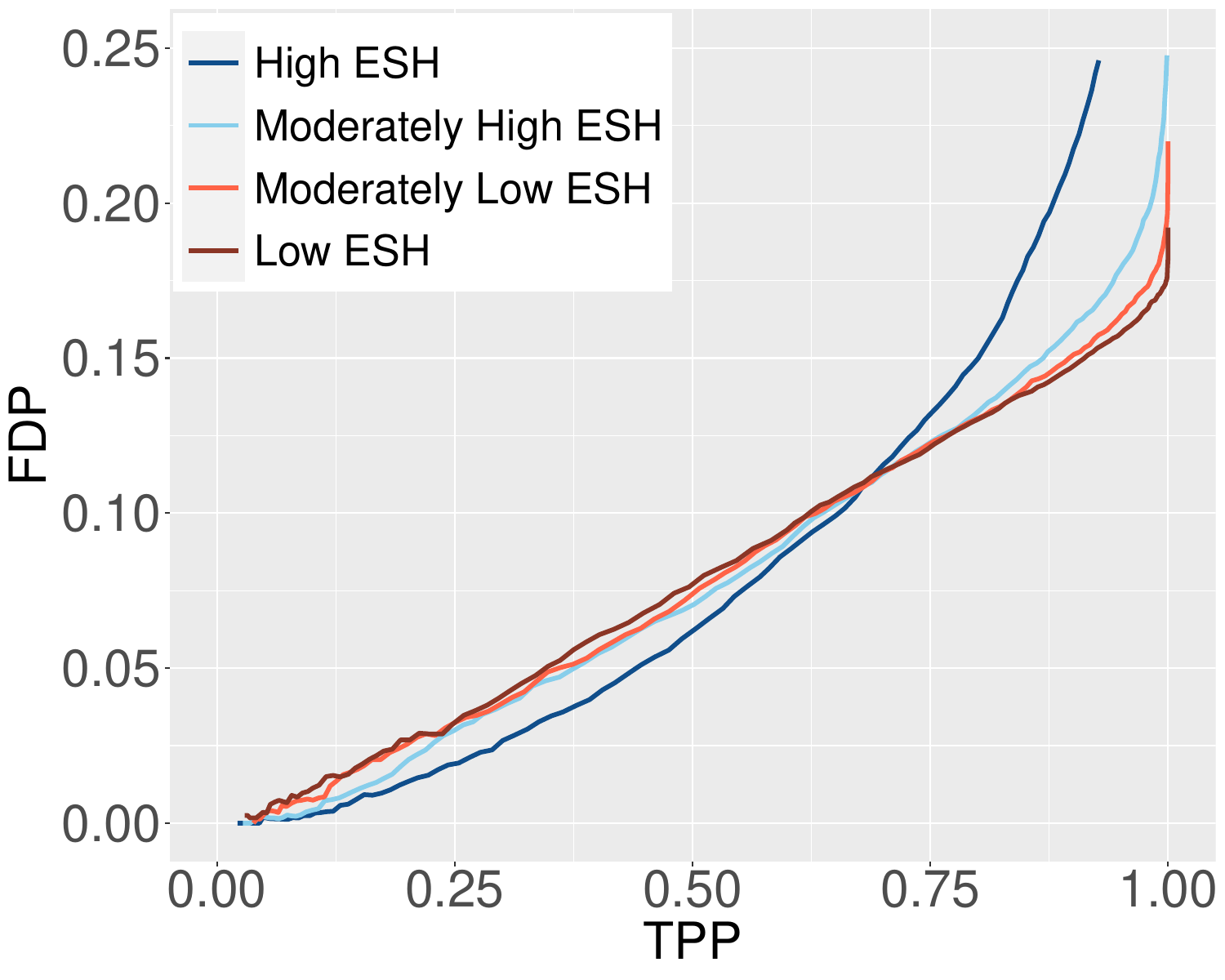}
\end{minipage}
\begin{minipage}{0.42\linewidth}
	\centering
	\includegraphics[width = 1\linewidth]{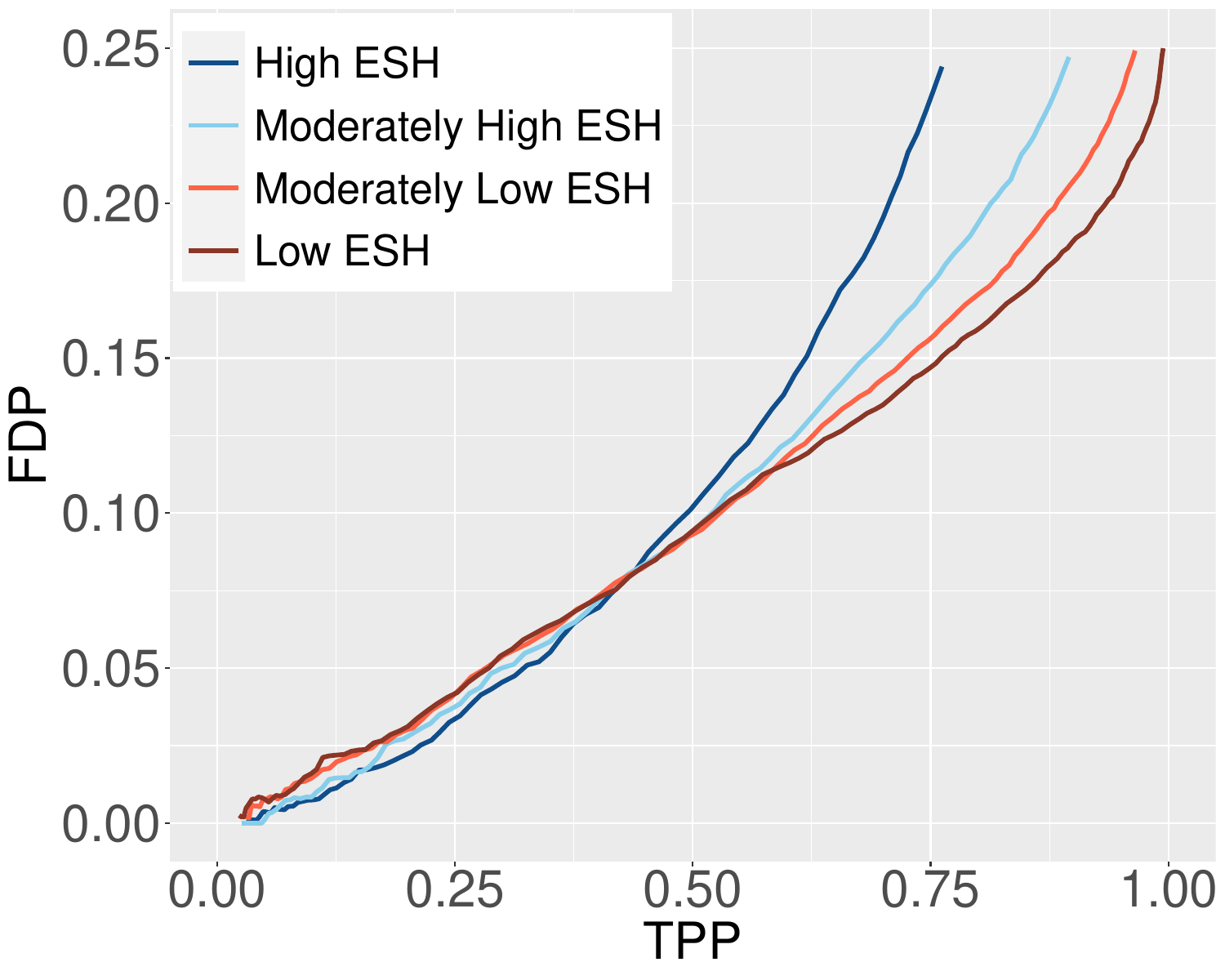}
\end{minipage}	
	\caption{The TPP--FDP trade-off plot with varying noise levels. The design matrix is specified by $n = p = 1000$, with i.i.d.~Gaussian entries. The regression coefficients are from \Cref{fig:config}, and the noise vector has i.i.d.~$\N(0, \sigma^2)$ entries, where $\sigma$ is set to $0.1, 0.2, 0.5$ and 1 in the top-left, top-right, bottom-left, and bottom-right panels, respectively. The mean FDP is obtained by averaging over $100$ replicates.}
	\label{fig:noisy}
\end{figure}

While \Cref{thm:upper} concerning the regime of low effect size heterogeneity only applies to the noiseless case, we make an attempt to show the impact of effect size heterogeneity in the noisy setting via simulations. Under an independent Gaussian random design of size $1000 \times 1000$, we set the nonzero regression coefficients to the four sets of effect sizes as depicted in \Cref{fig:config}. The noise term \(\bz\) consists of independent $\N(0, \sigma^2)$ entries with $\sigma = 0.1, 0.2, 0.5, 1.0$. The results are displayed in \Cref{fig:noisy}.

As with our previous simulation results, higher effect size heterogeneity tends to give rise to a better trade-off between the TPP and FDP from the beginning of the Lasso path. Interestingly, we observe a crossing point in each of the four panels of \Cref{fig:noisy} where higher heterogeneity undergoes a transition from giving a better trade-off down to a worse trade-off. In particular, the crossing point occurs earlier as the noise level $\sigma$ goes up. While it requires further research to understand this transition in a concrete manner, our observation is that the unselected effect sizes in the late stage of the Lasso path tend to be relatively small compared to the noise level, especially the effect sizes depicted in the bottom-right panel of \Cref{fig:noisy}, which have relatively high effect size heterogeneity. Intuitively, this crossing point is where signal-to-noise ratio becomes the dominant factor in place of effect size heterogeneity.

\section{Discussion}
\label{sec:discussion}

In this paper, we have proposed a concept termed effect size heterogeneity for measuring how diverse the nonzero regression coefficients are. Working under Gaussian random designs, we demonstrate that effect size heterogeneity has a significant impact on model selection consistency of the Lasso when the sparsity is linear in the ambient dimension. In short, we prove that the Lasso attains the optimal trade-off between true and false positives uniformly along its path when the effect sizes are strong and heterogeneous, and attains the worst trade-off when the effects are about the same size in the noiseless case. We also identify similar dependence of the rank of the first noise variable on effect size heterogeneity. While the two theoretical results are proved under certain assumptions, our simulations show that effect size heterogeneity has a significant impact on the Lasso estimate in a much wider range of settings.

Moving forward, this paper opens up several directions for future research. First, it is important to develop methods that incorporate the level of effect size heterogeneity for solving high-dimensional regression problems. In particular, one would be tempted to improve on the Lasso when the level of effect size heterogeneity is low. Interestingly, the SLOPE method has inadvertently addressed this question as its sorted $\ell_1$ penalty generally increases as the heterogeneity gets higher~\citep{slope,slopeminimax,bu2021characterizing}. Another related method developed from a Bayesian angle is the spike-and-slab Lasso procedure~\citep{rovckova2018spike}, which enables the adaptation to a mixture of large and small effects. Nevertheless, it is highly desirable to have methods that leverage effect size heterogeneity more directly. Moreover, a pressing question is to give a quantitative and formal definition of effect size heterogeneity. From a practical standpoint, regression coefficients are seldom exactly zero and thus it might be
more appropriate to consider the Type S error, which occurs when a nonzero effect is selected but with the incorrect sign~\citep{gelman2000type,barber2019knockoff}. This reality should prompt one to investigate how effect size heterogeneity interacts with the trade-off between the resulted directional FDP and power. Another question of practical importance is to examine carefully how the impact of effect size heterogeneity depends on the noise level. As an aside, given that Proposition~\ref{prop:logT} remains true for forward stepwise regression and least angle regression, we conjecture that Proposition~\ref{prop:wainwright} and Theorem~\ref{prop:rank_opt} also hold for the two model selection procedures. More broadly, it is of interest to investigate whether effect size heterogeneity retains its impact on other $\ell_1$ regularized methods such as the two-stage Lasso~\cite{wang2020bridge} and the Dantzig selector.



\subsection*{Acknowledgments}
We would like to thank Ma{\l}gorzata Bogdan, Emmanuel Cand\`es, Edward George, Pragya Sur, and Nancy Zhang for stimulating discussions. We are grateful to the referees and associate editor for their constructive comments that helped improve the presentation of this work. This work was supported in part by NSF through CAREER DMS-1847415, CCF-1763314, and CCF-1934876, the Wharton Dean's Research Fund, and a Facebook Faculty Research Award.

{\small
\bibliographystyle{abbrvnat}
\bibliography{ref}
}

\clearpage
\appendix
\section{Technical Proofs}
\label{sec:technical-proofs}

\subsection{Technical proofs for section~\ref{sec:heuristics}}
\label{sec:proofs-section-}

\subsubsection{Proof of proposition \ref{prop:wainwright}}
\begin{proof}[Proof of Proposition \ref{prop:wainwright}]
	We use the ``primal-dual witness'' argument in the Lasso literature (for example, see Theorem 2 in \cite{wainwright2009information}).
	As a reminder, here we consider the standard form of Lasso as in (\ref{eq:Lasso}).
	\begin{equation*}
	\bth = \argmin_{\bm{b} \in \R^p} \frac{1}{2}\|\by- \bX\bm{b}\|^2 + \lambda\|\bm{b}\|_1
	\end{equation*}
	with the model specified by (\ref{eq:basic_model}),
	\begin{equation*}
	\by = \bX\bbeta + \bz.
	\end{equation*}
	We define a pair $(\bth, \zh)\in \R^p\times \R^p$ to be {\it primal-dual optimal} if $\bth$ is a minimizer of (\ref{eq:Lasso}), and $\zh\in \partial \|\bth\|_1$, satisfying the zero-subgradient condition 
	\begin{equation*}
 	\X^T(\X\bth - \by) + \lambda \zh = 0.
	\end{equation*}
	For the convenience of analysis, we denote $\lambda_n = \frac{\lambda}{n}.$ Thus the condition above is equivalent to 
	\begin{equation}\label{eq:kkt}
	\frac{1}{n} \X^T(\X\bth - \by) + \lambda_n \zh = 0.
	\end{equation}
	By the sufficiency of KKT condition, we know that if there exists some $\zh$ such that the pair  $(\bth, \zh)\in \R^p\times \R^p$ satisfies (\ref{eq:kkt}), then $\bth$ is the solution to the Lasso. So $\zh$ can be seen as a ``dual witness'' showing $\bth$ is indeed a solution. We are therefore going to construct a ``dual witness'' vector $\zh$ to prove a certain $\bth$ is the solution to the Lasso.
	
	To concretely give our construction of $\bth$ and $\zh$, we fix an arbitrary small $\xi>0$, and then let $s = [(1-\xi)\frac{n-1}{1-\xi + 2\log p}]$.\footnote{This is only for technical convenience. One can easily verify that this condition is equivalent to $s = (1-o_\PP(1))\frac{2n}{\log p}$.} Denote $S_0 \equiv \{1, 2, ..., s\}$, $S_1 \equiv \{s+1, s+2, .., k\},$ and $S = S_0\cup S_1$, thus we have $S^C = \{k+1, ..., p\}$.
	Let $M(n) = n^a$ for some $a>\frac{1}{2}$, and let $\lambda_n =n^b$ for some $b$ that satisfies $(k-s)a -1<b<(k-s+1)a -\frac{3}{2}$. We omit the dependence of $M$ on $n$ in the following proof. For clarity, for any subset of $T\subset\{1,2,.., p\}$, we always use the notation $\bw_T$ to denote the restricted vector $(w_i)_{i\in T}$ of a vector $\bw$, and the notation $\bX_T$ to denote the restricted column matrix $(x_{i,j})_{j\in T}$ of a matrix $\bX$. We consider the following procedure to construct the pair $(\bth, \zh)$,
	\begin{enumerate}
		\item Let $\widehat{\bbeta}_{S_0^C} = 0$;
		\item Solve $(\bth_{S_0}, \zh_{S_0})\in\R^s\times\R^s$ from the following oracle sub-problem
		\begin{equation}\label{eq:oracle_lasso_problem}
		\bth_{S_0} \in \argmin_{\bm{b} \in \R^s} \left\{  \frac{1}{2}\|\by-\Xs\bm{b}\|^2_2+ \lambda \|\bm{b}\|_1 \right\},
		\end{equation}
		and choose $\zh_{S_0}\in \partial \|\bth_{S_0}\|_1$ such that 
		\begin{equation*}
		\frac{1}{n} \Xs^T(\Xs\bth_{S_0} - \by) + \lambda_n \zh_{S_0} = 0;
		\end{equation*}
		\item Given $\bth_{S_0}, \zh_{S_0}$, and $\widehat{\bbeta}_{S_0^C} = 0$, compute $\zh_{S_0^C}\in \R^{p-s}$ by equation (\ref{eq:kkt}), and check whether the {\it strict dual feasibility} condition $\|\zh_{S_0^C}\|_\infty < 1$ holds.
	\end{enumerate}
	The primal-dual witness construction guarantees that if a pair $(\bth,\zh)$ satisfies all the three conditions above, then $\bth$ is the unique solution of the Lasso \cite{wainwright2009information}. Once we prove our construction satisfies the conditions above, the second claim of Proposition \ref{prop:wainwright} is an easy corollary as we explicitly require $\bth_j = 0$ for all $j\in S_0^C$, and this gives $$\#\left\{ j: \widehat\beta_j(\lambda) \ne 0, \beta_j = 0 \right\} = 0.$$ And from this construction, it is also not hard to prove the first claim of the Proposition \ref{prop:wainwright}.
	With this protocol in mind, we proceed to prove that we can construct such a pair of $(\bth,\zh)$. Now, we solve $\bth_{S_0}, \zh_{S_0}$ from the subproblem in condition 2.
	Then, we set $\widehat{\bbeta}_{S_0^C} = 0$ as in condition 1, and solve $\zh_{S^C}\in \R^{p-s}$ from (\ref{eq:kkt}). To prove Proposition \ref{prop:wainwright}, we only need to prove that with this construction, the strict dual feasibility condition holds with high probability as $n, p\to \infty$.
	
	To prove this, we first simplify condition (\ref{eq:kkt}) by substituting $\widehat{\bbeta}_{S_0^C} = 0$, and write it in block matrix form as follows,	
	\begin{equation*}
	\frac{1}{n}
	\begin{bmatrix}
	\Xs^T\Xs & \Xs^T\Xss & \Xs^T\XSc\\ 
	\Xss^T\Xs & \Xss^T\Xss & \Xss^T\XSc\\ 
	\XSc^T\Xs & \XSc^T\Xss & \XSc^T\XSc
	\end{bmatrix}
	\cdot 
	\begin{bmatrix}
	\bth_{S_0} - \bs\\ 
	-\bss\\ 
	0
	\end{bmatrix}
	- \frac{1}{n} 
	\begin{bmatrix}
	\Xs^T\bz \\ 
	\Xss^T\bz\\ 
	\XSc^T\bz
	\end{bmatrix} +
	\lambda_n \begin{bmatrix}
	\zh_{S_0} \\ 
	\zh_{S_1}\\ 
	\zh_{S^C}
	\end{bmatrix} = 
	\begin{bmatrix}
	0 \\ 
	0\\ 
	0
	\end{bmatrix},
	\end{equation*}
	or equivalently, 
	\begin{eqnarray}
	&&\frac{1}{n}\Xs^T\Xs(\bth_{S_0}-\bs) - \frac{1}{n}\Xs^T\Xss\bss - \frac{1}{n}\Xs^T\bz + \lambda_n\zh_{S_0} = 0, \label{eq:A1}\\
	&&\frac{1}{n}\Xss^T\Xs(\bth_{S_0}-\bs) - \frac{1}{n}\Xss^T\Xss\bss - \frac{1}{n}\Xss^T\bz + \lambda_n\zh_{S_1} = 0, \label{eq:A2}\\
	&&\frac{1}{n}\XSc^T\Xs(\bth_{S_0}-\bs) - \frac{1}{n}\XSc^T\Xss\bss - \frac{1}{n}\XSc^T\bz + \lambda_n\zh_{S^C} = 0. \label{eq:A3}
	\end{eqnarray}
	By (\ref{eq:A1}), we have
	\begin{equation}\label{eq:beta_difference}
	\bth_{S_0}-\bs = \left(\Xs^T\Xs\right)^{-1} \left[\Xs^T\Xss\bss + \Xs^T\bz\right] - \lambda_n n\left(\Xs^T\Xs\right)^{-1}\zh_{S_0}.
	\end{equation}
	By substituting (\ref{eq:beta_difference}) into (\ref{eq:A2}) and (\ref{eq:A3}), we can solve $\widehat{w}_j$ for any $j\in S_0^C$ as
	\begin{align}
	\widehat{w}_j = & ~-\frac{1}{\lambda_n n}\Xj^T\Xs(\bth_{S_0}-\bs) + \frac{1}{\lambda_n n}\left[\Xj^T\Xss\bss + \Xj^T\bz\right] \notag\\
	=&~\Xj^T\Xs\left(\Xs^T\Xs\right)^{-1}\zh_{S_0} \notag\\&~- \frac{1}{\lambda_n n}\Xj^T\Xs\left(\Xs^T\Xs\right)^{-1}\left[\Xs^T\Xss\bss + \Xs^T\bz\right] \notag\\&~+ \frac{1}{\lambda_n n}\left[\Xj^T\Xss\bss + \Xj^T\bz\right] \notag\\
	=&~\underbrace{\Xj^T\Xs\left(\Xs^T\Xs\right)^{-1}\zh_{S_0}}_{v_j} + \underbrace{\Xj^T\Ps\left[\frac{\bz}{\lambda_n n}+ \Xss\frac{\bss}{\lambda_n n}\right]}_{u_j}, \label{eq:vjuj}
	\end{align}
	where $\Ps = \bI - \Xs\left(\Xs^T\Xs\right)^{-1}\Xs^T$. As mentioned previously, our goal is to show the strict dual feasibility condition $\max_{j\in S_0^C} |\widehat{w}_j| < 1$ holds with high probability. We will prove it by analyzing the two terms $u_j$ and $v_j$ separately. Specifically, we prove that $v_j<1-\frac{\xi}{16}$ and $u_j\to0$ with high probability.
	
	Denote $M_n = \frac{1}{n}\zh_{S_0}^T\left(\Xs^T\Xs\right)^{-1}\zh_{S_0}$.
	Conditioning on the event $E = \{\zh_{S_0} = sign(\bbeta_{S_0})\}$ and its compliment gives us
	\begin{align*}
	\PP\left(\max_{j\in S_0^C} |v_j|\ge c\right) &\le \PP\left(\max_{j\in S_0^C} |v_j|\ge c~\biggr|E\right)\P(E) + \PP(\max_{j\in S_0^C} |v_j|\ge c,\text{ and } E^C)\\
	&\le \PP\left(\max_{j\in S_0^C} |v_j|\ge c~\biggr|E\right) + \PP(E^C).
	\end{align*}
	It can be seen  through Lemma \ref{lem:sign_consistence} that the second term of the last display tends to zero; For the first term, we let $T(\vartheta)$ denote the event $\{|M_n-\EE M_n|\ge \vartheta\EE M_n\}$. Similar as before, conditioning on the event $T(\vartheta)$ and its complement gives for any $c \in (0,1)$,
	\begin{align*}
	\PP\left(\max_{j\in S_0^C} |v_j|\ge c\biggr|E\right)&\le \PP\left(\max_{j\in S_0^C} |v_j|\ge c~\biggr|T(\vartheta)^C\cap E\right) + \PP(T(\vartheta)\cap E).
	\end{align*}
	By Lemma \ref{lem:vj} and Lemma \ref{lem:Mn}, the second term in the last display goes to $0$ as $n\to \infty$ faster than $\frac{2}{\vartheta^2(n-s-3)}$. And for the first term, we tackle it by considering $\max_{j\in S^C_0} v_j$ and $\min_{j\in S^C_0} v_j$ separately. Denote the event $T = T(\vartheta)^C\cap E$ for convenience, we have
	\begin{align*}
	\PP\left(\max_{j\in S_0^C} v_j \ge c \biggr|T\right) \le \PP\left(\max_{j\in S_0^C} \widetilde{v}_j \ge c\right),
	\end{align*}
	where $\widetilde{v}_j$ are i.i.d. from $\Normal(0, (1+\vartheta)\EE [M_n|E]) = \Normal(0, (1+\vartheta)\frac{s}{n-s-1})$.
	This inequality follows from Lemma \ref{lem:normalincrease}, which states that the probability of the event $\{\max_{i\in S_0^C} v_i\ge c\}$ increases as the mean and variance of $v_i$ increase for Gaussian variables $v_i$. Given the event $T$, the maximum variance of $v_j$ is simply $(1+\vartheta)\EE [M_n|E]$, and thus we have the inequality above. Set $c = \varpi + \EE \max_{j\in S_0^C} \widetilde{v}_j$. From Lemma \ref{lem:maxbound}, we have 
	\begin{align*}
	\PP\left(\max_{j\in S_0^C} v_j \ge c \biggr|T\right)\le\PP\left(\max_{j\in S_0^C} \widetilde{v}_j > \varpi + \EE \max_{j\in S_0^C} \widetilde{v}_j\right)\le & (p-s) \exp\left(-\frac{\varpi^2}{2(1+\vartheta)\EE [M_n|E]}\right).
	\end{align*}
	A similar argument for $\min_{j\in S_0^C} \widetilde{v}_j$ gives us
	\begin{align*}
	\PP\left(\min_{j\in S_0^C} {v}_j < -c\biggr|T\right)= \PP\left(-\min_{j\in S_0^C} {v}_j > \varpi + \EE \max_{j\in S_0^C} \widetilde{v}_j\biggr|T\right)
	\le & (p-s) \exp\left(-\frac{\varpi^2}{2(1+\vartheta)\EE [M_n|E]}\right).
	\end{align*}
	Combining the two inequalities above yields
		\begin{align}\label{eq:vjc}
	\PP\left(\max_{j\in S_0^C} |{v}_j| > \varpi + \EE \max_{j\in S_0^C} \widetilde{v}_j\biggr|T\right)\le & 2(p-s) \exp\left(-\frac{\varpi^2}{2(1+\vartheta)\EE [M_n|E]}\right) \to 0.
	\end{align}
	By Lemma \ref{lem:Mn} and (\ref{eq:vjc}), we get
	\begin{equation}\label{eq:vj}
	\PP\left(\max_{j\in S_0^C} |v_j|\ge \varpi + \EE \max_{j\in S_0^C} \widetilde{v}_j\right) \le \frac{2}{\vartheta^2(n-s-3)} + 2(p-s) \exp\left(-\frac{\varpi^2(n-s-1)}{2s(1+\vartheta)}\right).
	\end{equation}
	With the relation of $s = [(1-\xi)\frac{n-1}{1-\xi +2\log p}]$, we can set $\varpi = \frac{\xi}{8}, \vartheta = \frac{(1-\frac{\xi}{4})^2}{1-\xi}-1 >0$, and obtain	
	\begin{align}
	\varpi + \sqrt{(1+\vartheta)\frac{s}{n-s-1} 2\log(p-s)} 
	\le& \varpi + \sqrt{(1+\vartheta)\frac{(1-\xi)\frac{n-1}{1-\xi +2\log p} +1}{n-(1-\xi)\frac{n-1}{1-\xi +2\log p}-1} 2\log(p)}\notag\\
	=& \varpi + \sqrt{(1+\vartheta) \frac{(1-\xi)(n-1) + (1-\xi+\log p)}{2\log p (n-1)}2 \log p} \notag\\
	= & \varpi + \sqrt{(1+\vartheta)\left( (1-\xi)+\frac{ (1-\xi+\log p)}{(n-1)}\right)}\notag \\
	\le & \frac{\xi}{8} + (1-\frac{\xi}{4}) + \frac{C}{n-1}\notag\\
	< & 1 - \frac{\xi}{16} ~~~\text{for some large $n$.} \label{eq:lessthan1}
	\end{align}
	Combining this with the well-known fact of the expectation of the maximum of i.i.d.~Gaussian variables that $\E \max_{j\in S_0^C}|\widetilde{v}_j|\le\sqrt{(1+\vartheta)\frac{s}{n-s-1} 2\log(p-s)} $ and (\ref{eq:vj}), we know
	\begin{align}\label{eq:A10}
	\P(\max_{j\in S_0^C}|v_j|\ge1-\frac{\xi}{16})\le \frac{2}{\vartheta^2(n-s-3)} + 2(p-s) \exp\left(-\frac{\varpi^2(n-s-1)}{2s(1+\vartheta)}\right) + \P(E^C).
	\end{align}
	This bound is good enough for our purpose. We now proceed to obtain a similar bound of $u_j$.
	
	By the Cauchy--Schwarz inequality, we have
	\begin{equation}\label{eq:ujbound}
	|u_j| \le \left\|\Xj^{T}\Ps\right\|\cdot \left\|\frac{\bz}{\lambda_n n}+ \Xss\frac{\bss}{\lambda_n n}\right\|.
	\end{equation}
	Therefore, we can bound $|u_j|$ if we can control the two norms in \eqref{eq:ujbound} separately. 
	Because all eigenvalues of $\Ps$ are less than $1$, we have
	\begin{equation*}
	\left\|\Xj^{T}\Ps\right\|^2 \le \|\Xj\|^2 = \sum_{i=1}^n \mathbf{X}_{ij}^2 \stackrel{\mathcal{D}}{=} \frac{1}{n}\sum_{i=1}^n W_i^2,
	\end{equation*}
	where $W_i \stackrel{i.i.d.}{\sim} \Normal(0,1)$. The summation $\sum_{i=1}^n W_i^2$ is thus a $\chi^2$-distribution with degree of freedom $n$. By Lemma \ref{lem:chisquare}, for any $t_1>0$, we have
	\begin{equation*}
	\P\left(\left|\frac{\|\Xj\|^2}{n}-1\right|\ge t_1\right)\le 2 \exp(-nt_1^2/8).
	\end{equation*}
	Combining this with $\left\|\Xj^{T}\Ps\right\| \le \|\Xj\|$ gives us
	\begin{equation}\label{eq:xj}
	\|\Xj^{T}\Ps\|_2\le \sqrt{1+t_1} ,~~~\text{w.p.~} \ge ~1- 2 \exp(-nt_1^2/8).
	\end{equation}
	Finally, let $\tuj$ denote
	\begin{equation*}
	\tuj = \left\|\frac{\bz}{\lambda_n n}+ \Xss\frac{\bss}{\lambda_n n}\right\|.
	\end{equation*}
	It is easy to realize that 
	\begin{equation*}
	e_i\cdot ({\bz}+ \Xss{\bss})\sim \Normal(0, \sigma'^2), ~~\text{for all } 1\le i \le s,
	\end{equation*}
	where $e_i\in \R^n$ is the $i$-th standard unit vector and $\sigma'^2 =\sigma^2 + \frac{M^{2(k-s+1)}-1}{n(M^2-1)}$. By easy calculation, we obtain
	\begin{equation*}
	\E(\|\tuj\|^2) = n \cdot \frac{\sigma'^{2} }{\lambda_n^2n^2}.
	\end{equation*}
	Applying Lemma \ref{lem:chisquare} with $Z_i = e_i\cdot ({\bz}+ \Xss{\bss})$ again, we know for any $t_2\ge0$
	\begin{equation*}
	\P\left(\left|\frac{\|\tuj\|^2}{ \E(\|\tuj\|^2)}-1\right|\ge t_2\right)\le 2 \exp(-nt_2^2/8),
	\end{equation*}
	which is equivalent to
	\begin{equation}\label{eq:ujtilde}
	\|\tuj\| \le \sqrt{1+t_2}\cdot \frac{\sigma' }{\lambda_n \sqrt{n}}, ~~\text{w.p.}~\ge ~1- 2 \exp(-nt_2^2/8).
	\end{equation}
	Using (\ref{eq:ujbound}) by combining the two bounds (\ref{eq:xj}) and (\ref{eq:ujtilde}), we get
	\begin{equation}\label{eq:uj}
	\P(\max_{j\in S_0^C}|u_j|\ge \sqrt{1+t_1}\sqrt{1+t_2}\cdot \frac{\sigma' }{\lambda_n \sqrt{n}}) \le 2(p-s) (\exp(-nt_1^2/8)+ 2 \exp(-nt_1^2/8)).
	\end{equation}
	Now, we can set $M = n^a$ for some $a>\frac{1}{2}$, $\lambda_n =n^b$ for some $b$ that satisfies $(k-s)a -1<b<(k-s+1)a -\frac{3}{2}$, and $t_1=t_2=1$ to obtain
	\begin{align}
	\sqrt{1+t_1}\sqrt{1+t_2}\cdot \frac{\sigma' }{\lambda_n \sqrt{n}}\notag
	=& 2 \frac{\sqrt{\sigma^2 + \frac{M^{2(k-s+1)}-1}{n(M^2-1)}}}{\lambda_n\sqrt{n}}\notag \\
	\lesssim & 2 n ^{(k-s)a -1 - b} \to 0.\label{eq:to0}
	\end{align}
	Particularly, when $n$ is large enough, $\sqrt{1+t_1}\sqrt{1+t_2}\cdot \frac{\sigma' }{\lambda_n \sqrt{n}}$ is less than $\frac{\xi}{32}$, which in turn gives
	\begin{equation}\label{eq:A16}
	\P(\max_{j\in S_0^C}|u_j|\ge \frac{\xi}{32}) \lesssim 2(p-s) (\exp(-nt_1^2/8)+ 2 \exp(-nt_1^2/8)).
	\end{equation}
	And thus by a union bound and then (\ref{eq:vjuj}, \ref{eq:A10}, \ref{eq:A16}), we have
	\begin{align}
	&~\P\left(\max_{j\in S_0^C}|\widehat{w}_j|\ge1-\frac{\xi}{32}\right)\notag\\
	\le&~\P\left(\max_{j\in S_0^C}|v_j|\ge 1-\frac{\xi}{16}\right) + \P\left(\max_{j\in S_0^C}|u_j|\ge \frac{\xi}{32}\right)+\P(E^C) \notag\\
	\lesssim&~~\frac{2}{\vartheta^2(n-s-3)} + 2(p-s)\left( \exp(-\frac{\varpi^2(n-s-1)}{2s(1+\vartheta)})+ \exp(-nt_2^2/8) +  \exp(-nt_1^2/8) \right).\label{eq:ziineq}
	\end{align}
	We observe that as $n\to\infty$, 
	\begin{equation}
	\frac{2}{\vartheta^2(n-s-3)} + 2(p-s)\left( \exp(-\frac{\varpi^2(n-s-1)}{2s(1+\vartheta)})+ \exp(-nt_2^2/8) +  \exp(-nt_1^2/8) \right) \to 0, \label{eq:prob0}
	\end{equation}
	which simply implies that as $n\to\infty$,
	\begin{equation*}
	\P(\max_{j\in S_0^C}|\widehat{w}_j|\ge 1) \to 0.
	\end{equation*}
	Thus, we have proven under our construction, strict dual feasibility holds. And as we pointed out at the beginning of the proof, the second part of the proposition \ref{prop:wainwright} holds, since we set $\bth_j= 0$ for any $j\in S_0^C$. Therefore, we obtain
	\begin{equation*}
	\#\left\{ j: \widehat\beta_j(\lambda) \ne 0, \beta_j = 0 \right\} = 0.
	\end{equation*}
	
	Now, we proceed to prove the first part of the proposition, that is,
	\begin{equation*}
	\#\left\{ j: \widehat\beta_j(\lambda) \ne 0, \beta_j \ne 0 \right\} = s =(1 - o_{\P}(1))\frac{n}{2\log p}.
	\end{equation*}
	Observe that the second equality is due to our assumption on $s$. So we only need to prove the first equality, that is, for all $j \in S_0$, $\widehat\beta_j$ are non-zero, and thus the total number of non-zero $\bth$'s is exactly $s$.
	To show this, observe that if $\beta_j -\widehat{\beta}_j <\beta_j$, it is clear that $\widehat{\beta}_j >0$. Therefore, it suffices to show the following inequality
	\begin{align*}
	\max_{j \in S_0}(\beta_j - \widehat{\beta}_j) < \min_{j \in S_0} \beta_{j} = M^{k-s+1}\equiv \rho
	\end{align*}
	holds with probability tending to 1. We denote
	\begin{equation}\label{eq:yi}
	Y_i = -e_i^T \cdot \left(\Xs^T\Xs\right)^{-1}\Xs^T \left[\Xss\bss + \bz\right] + e_i \cdot \lambda_n n \left(\frac{\Xs^T\Xs}{n}\right)^{-1}\zh_{S_0},
	\end{equation}
	where $e_i\in \R^n$ is the $i$-th standard unit vector. 
	By (\ref{eq:beta_difference}), we know
	\begin{align*}
	\max_{j \in S_0}(\beta_j - \widehat{\beta}_j) = \max_{1\le i\le n} Y_i.
	\end{align*}
	So it is equivalent to show that $\max_{1\le i\le n} Y_i\ge \rho$ holds with probability tending to zero. By Lemma \ref{lem:yi}, we know for $E_i = \E(Y_i\big|\Xs)$, 
	and $V_i = \VAR(Y_i\big|\Xs)$, the event 
	\begin{equation*}
	A \equiv \bigcup_{i=1}^s \big\{|E_i|\ge (1 + \sqrt{n})|\E(E_i)|, \text{  or  }|V_i| \ge 2\E(V_i) \big\}
	\end{equation*}
	has probability
	\begin{equation*}
	\P(A)\le \frac{sK}{n-s} \to 0, \text{  as}~ n\to\infty.
	\end{equation*}
	By conditioning on the event $A$ and its complement, we have
	\begin{align*}
	\P(\max_{i\in S_0} Y_i \ge \rho) &\le \P(\max_{i\in S_0}Y_i \ge \rho \big| A^C) + \P(A) \notag\\
	&\le \P(\max_{i\in S_0} \tyi \ge \rho) + \frac{K}{\frac{n}{s}-1},
	\end{align*}
	where $\tyi \stackrel{i.i.d.}{\sim} \Normal((1 + \sqrt{n})\E(E_i), 2\E(V_i))$ and the second inequality used the fact in Lemma \ref{lem:normalincrease} that the probability of the event $\{\max_{i\in S_0 } Y_i \ge \rho\}$ increases as the mean and variance increase as long as the mean is less than $\rho$, which can be directly verified by
	\begin{equation}
	(1 + \sqrt{n})\E(E_i) =(1 + \sqrt{n})\frac{\lambda_n n^2}{n-s-1}< \rho.\label{eq:meansmall}
	\end{equation}
	Markov's inequality then gives us
	\begin{align}
	\P(\max_{i\in S_0} \tyi \ge \rho) &\le \frac{1}{\rho}\E\left(\max_{i\in S_0} |\tyi|\right)\notag\\
	&\le\frac{1}{\rho} \left( \E(\tyi) + \E\left(\max_{i\in S_0} |\tyi - \E(\tyi)|\right) \right)\notag \\
	& \le \frac{1}{\rho} \left(  \frac{(1+\sqrt{n})\lambda_n n^2}{n-s-1} + 3 \sqrt{\frac{2\sigma'^2 \log s}{n-s-1}}\right),\label{eq:probzero}
	\end{align}
	where the last inequality uses the bound on Gaussian maxima in Lemma \ref{lem:gaussian_maxima}.
	By the relation in (\ref{eq:meansmall}), we can easily verify that the probability in (\ref{eq:probzero}) converges to zero under our conditions of $M = n^a$, $\rho = M^{k-s+1}$ and $\lambda_n =n^b$, where $a$ and $b$ satisfy $a>\frac{1}{2}$ and $(k-s)a -1<b<(k-s+1)a -\frac{3}{2 }$, as $n\to\infty$. 
\end{proof}

\subsubsection{Miscellaneous lemmas for \Cref{sec:heuristics}}
We first state a well-known result in the random matrix theory (see, for example, Theorem 5.2 in \cite{baraniuk2008simple}) that we use in the proof of \Cref{prop:rank_opt}. Then we list all the necessary lemmas for proving Proposition \ref{prop:wainwright}.
\begin{lemma}\label{lm:rip}
	Under the working assumptions, for any deterministic $1 \le m \le p/2$, the matrix spectrum norm $\|\cdot\|_2$ satisfies
	\[
	\max_{|S| \le m} \left\|\bX_S^\top \bX_S - \bI\right\|_2 \le C\sqrt{\frac{m \log(p/m)}{n}}
	\]
	with probability $1 - 1/p^2$, where $C$ is a universal constant and $T$ is any set of column indices.
\end{lemma}

\begin{lemma}\label{lem:vj}
	For $v_j$ defined in (\ref{eq:vjuj}), and any $i, j\in S_0^C$, we have the following facts,
	\begin{enumerate}
		\item $\EE(v_j|\Xs) = 0$;
		\item $\VAR(v_j|\Xs) = \frac{1}{n} \zh_{S_0}^T(\Xs^T\Xs)^{-1}\zh_{S_0}$;
		\item $\Cov(v_j, v_i|\Xs) = 0$, if $i\ne j$.
	\end{enumerate}
\end{lemma}
\begin{proof}[Proof of Lemma \ref{lem:vj}]
	
	Since $v_j = \Xj^T\Xs\left(\Xs^T\Xs\right)^{-1}\zh_{S_0}$, $j\in S_0^C$, and $X_j \indep \Xs$, fact $1$ follows from $X_j$ being a centered Gaussian variable.
	
	For fact 2 and fact 3, we observe
	\begin{align*}
	\Cov(V_j, V_i|\Xs) =& \EE\left(\zh_{S_0}^T\left(\left.\Xs^T\Xs\right)^{-1}\Xs^T\Xj\mathbf{X}_i^T\Xs\left(\Xs^T\Xs\right)^{-1}\zh_{S_0}\right|\Xs\right) \\
	=&\zh_{S_0}^T\left(\Xs^T\Xs\right)^{-1}\Xs^T\EE\left(\left.\Xj\mathbf{X}_i^T\right|\Xs\right)\Xs\left(\Xs^T\Xs\right)^{-1}\zh_{S_0}\\
	=& \begin{cases}
	\frac{1}{n}\zh_{S_0}^T\left(\Xs^T\Xs\right)^{-1}\zh_{S_0} ~~~~&\text{if } i= j,\\
	0 ~~~~& \text{if } i\ne j.
	\end{cases}
	\end{align*}
\end{proof}

\begin{lemma}\label{lem:wishart}
	Consider $\Xs\in \R^{n\times s}$, and suppose each of its column $\mathbf{X}_{j}\stackrel{i.i.d.}{\sim} \Normal(0, \bm{\Sigma})$ where $\bm{\Sigma}\in \R^{s\times s}$ is positive definite. Then $\Xs^T\Xs$ is a Wishart distribution of degree of freedom $n$, and $(\Xs^T\Xs)^{-1}$ is the inverse Wishart distribution, with expectation and variance 
	\begin{enumerate}
		\item $\E(\Xs^T\Xs)^{-1} = \frac{\bm{\Sigma}^{-1}}{n-s-1}$;
		\item $\VAR[(\Xs^T\Xs)^{-1}_{i,j}] = \frac{(n-s+1)(\Sigma^{-1}_{i,j})^2 + (n-s-1)\Sigma^{-1}_{i,i}\Sigma^{-1}_{j,j}}{(n-s)(n-s-1)^2(n-s-3)}$.
	\end{enumerate}
\end{lemma}
\begin{proof}[Proof of Lemma \ref{lem:wishart}]
	See for example  Lemma 7.7.1 of \cite{anderson1962introduction} and the formula for the second moment of the inverse Wishart matrices in \cite{siskind1972second}.
\end{proof}

\begin{lemma}\label{lem:Mn}
	Let $M_n = \frac{1}{n}\zh_{S_0}^T\left(\Xs^T\Xs\right)^{-1}\zh_{S_0}$. Conditioned on the event $E$, that is, $\zh_{S_0} = sign(\bbeta_{S_0})$, we have the following facts:
	\begin{enumerate}
		\item $\EE(M_n|E) = \frac{s}{n-s-1}$;
		\item $\VAR(M_n|E) = \frac{2s^2}{(n-s-1)^2(n-s-3)}$;
		\item $\forall \vartheta >0$, $\PP[|M_n-\EE(M_n)| \ge \vartheta \EE(M_n)|E] \le \frac{2}{\vartheta^2(n-s-3)} $.
	\end{enumerate} 
\end{lemma}
\begin{proof}[Proof of Lemma \ref{lem:Mn}]
	Observe that $\Xs^T\Xs$ follows the Wishart distribution with variance $\frac{1}{n}\bI_{S_0}$, and thus by Lemma \ref{lem:wishart}, the matrix $(\Xs^T\Xs)^{-1}$ is the inverse Wishart distribution with mean 
	\begin{equation}\label{eq:wishart}
	\EE(\Xs^T\Xs)^{-1} = \frac{n}{n-s-1} \bI_{S_0}.
	\end{equation}
	Notice that $\widehat{w}_i = \pm 1$ for all $i\in S_0$, and when conditioned on $E$, it is equal to $sign(\bbeta_{S_0})$, which is independent of $X_{S_0}$. Therefore, we have
	\begin{equation*}
	\EE(M_n|E) = \frac{1}{n} \frac{n}{n-s-1} \zh_{S_0}^\top \bI_{S_0}\zh_{S_0} = \frac{s}{n-s-1}.
	\end{equation*}
	
	To calculate the second moment of the inverse Wishart matrices (\cite{siskind1972second}), we have that for $n - s - 3 > 0 $,
	\begin{align*}
	\EE(M_n^2|E) =& \frac{1}{n^2}\frac{1}{(n-s)(n-s-3)}(n\cdot \zh_{S_0}^\top \bI_{S_0}\zh_{S_0})^2 \frac{n-s}{n-s-1}\notag\\
	=& \frac{s^2}{(n-s-1)(n-s-3)}.
	\end{align*}
	Therefore, combining the two equations above, we obtain
	\begin{align*}
	\Var(M_n|E) =&  \frac{s^2}{(n-s-1)(n-s-3)} - \frac{s^2}{(n-s-1)^2} \notag \\
	=& \frac{2s^2}{(n-s-1)^2(n-s-3)}.
	\end{align*}
	
	For the third statement, Markov's inequality gives us
	\begin{align*}
	\PP(|M_n-\EE(M_n)| \ge \vartheta \EE(M_n | E) \le & \frac{\VAR(M_n | E)}{\vartheta^2 (\EE(M_n|E))^2} \\
	= & \frac{\frac{2s^2}{(n-s-1)^2(n-s-3)}}{\vartheta^2\frac{s^2}{(n-s-1)^2}}\\
	= & \frac{2}{\vartheta^2(n-s-3)}.
	\end{align*}
\end{proof}

\begin{lemma}\label{lem:maxbound}
	Consider i.i.d.~Gaussian random variables $z_j \sim \Normal(0, \sigma^2)$, where $j = 1,..,l$ for some $l\ge 2$. We have for any $\varpi>0$,
	\begin{equation*}
	\PP\left(\max_{1\le j \le l} z_j > \varpi + \EE\left(\max_{1\le j \le l} z_j\right)\right) \le \e^{-\frac{\varpi^2}{2\sigma^2}}.
	\end{equation*}
\end{lemma}

\begin{proof}[Proof of Lemma \ref{lem:maxbound}]
	By the Gaussian tail bound 
	\begin{align*}
	\PP(z_j > \varpi)\le \frac{\sigma}{\sqrt{2\pi}\varpi}\e^{-\frac{\varpi^2}{2\sigma^2}},
	\end{align*}
	and the well-known fact for the expectation of maximum of i.i.d.~Gaussian variables
	\begin{align*}
	\EE\left(\max_{1\le j \le l} z_j\right) \le \sigma\sqrt{2\log l},
	\end{align*}
	we have the following union bound
	\begin{align*}
	\PP\left(\max_{1\le j \le l} z_j > \varpi + \EE[\max_{1\le j \le l} z_j]\right) \le & l\frac{1}{\sqrt{2\pi}(\frac{\varpi}{\sigma}+\sqrt{2\log l)}}\cdot \e^{\frac{(\varpi + \sqrt{2\log(l)} \sigma)^2}{2\sigma^2}}\\
	=& \frac{1}{\sqrt{2\pi}(\frac{\varpi}{\sigma}+\sqrt{2\log l)}}\exp\left(-\frac{\varpi^2}{2\sigma^2}\right) \cdot \exp\left(\sqrt{2\log l}\frac{\varpi}{\sigma}\right)\\
	\le &\exp\left(-\frac{\varpi^2}{2\sigma^2}\right)
	\end{align*}
	holds as long as $l\ge 2$.
\end{proof}

\begin{lemma}\label{lem:chisquare}
	Consider $Z_i\stackrel{i.i.d.}{\sim} \Normal(0, \theta^2)$, and denote $Z = \sum_{i=1}^n Z_i^2$. For $t>0$, we have the inequality
	$$
	\P\left(\left|\frac{Z}{\E(Z)}-1\right|\ge t\right) \le 2\exp(-nt^2/8).
	$$
\end{lemma}

\begin{proof}[Proof of Lemma \ref{lem:chisquare}]
	Let $\widetilde{Z}_i $ be defined as
	$$\widetilde{Z}_i = \frac{Z_i}{\theta}\stackrel{i.i.d.}{\sim} \Normal(0, 1).$$ 
	We have
	$$
	\E(Z) = \sum_{i=1}^n \theta^2 = n \theta^2,
	$$
	and	therefore
	$$
	\frac{Z}{\E(Z)} = \frac{\sum_{i=1}^n \widetilde{Z}_i}{n} \stackrel{\mathcal{D}}{=} \frac{\chi^2}{n}.
	$$
	By easy calculation, we obtain
	\begin{align*}
	\E\left(\e^{\lambda (\widetilde{Z}_i^2-1)}\right) & =\frac{1}{\sqrt{2 \pi}} \int_{-\infty}^{+\infty} \e^{\lambda\left(z^{2}-1\right)} \e^{-z^{2} / 2} d z \\ &=\frac{\e^{-\lambda}}{\sqrt{1-2 \lambda}} \le \e^{2\lambda^2}.
	\end{align*}
	This means $Z$ is sub-exponential with parameter $(2,4)$ (definition of a sub-exponential variable is standard, so we refer the reader to, for example, the Definition 2.2 and Example 2.11 in \cite{wainwright2019high}), and thus $Z$ is a sub-exponential variable with parameter $(2\sqrt{n},4)$. By the Bernstein inequality, we obtain
	$$
	\P\left(\left|\frac{Z}{\E(Z)}-1\right|\ge t\right) \le 2\exp(-nt^2/8).
	$$
\end{proof}

\begin{lemma}\label{lem:yi}
	For $Y_i$ defined in (\ref{eq:yi}), we have the following facts:
	\begin{enumerate}
		\item Denote $E_i = \E(Y_i\big|\Xs)$, then we have $|\E(Y_i)| = |\E(E_i)| = \frac{\lambda_n n^2}{n-s-1}$;
		\item Denote $V_i = \VAR(Y_i\big|\Xs)$, then we have $\E(V_i)  = 1$;
		\item For $n$ sufficiently large, the inequality $\P\big(|E_i|\ge (1 + \sqrt{n})|\E(E_i)|, ~\text{or  }|V_i| \ge 2\E(V_i)\big) \le \frac{K}{n-s}$ holds for some constant $K$ independent of $n$ and $s$.
	\end{enumerate}
\end{lemma}
\begin{proof}[Proof of Lemma \ref{lem:yi}]The idea of the following proof is adapted from Lemma 6 of \cite{wainwright2009information}.
	
	For part (a), since $\Xss \indep \Xs$ and $\bz \indep \Xs$, we get
	\begin{equation*}
	E_i = \E(Y_i\big| \Xs) = -\lambda_nn e_i^T \left({\Xs^T\Xs}\right)^{-1} \zh_{S_0}.
	\end{equation*}
	Thus, we have
	\begin{align*}
	|\E(Y_i)| &= \left|\E\left(-\lambda_n n e_i \left({\Xs^T\Xs}\right)^{-1} \zh_{S_0}\right)\right|\\
	&= \left|-\lambda_n n e_i^T \frac{n}{n-s-1}\bI_{S_0}^{-1}\zh_{S_0}\right|\\
	&= \frac{\lambda_n n^2}{n-s-1},
	\end{align*}
	where the second equality is by (\ref{eq:wishart}) for the mean of the inverse Wishart distribution.
	
	Next, we turn to prove part (b). We observe that each entry of vector $(\Xss\bss + \bz)$ is i.i.d.~distributed as $\Normal(0, \sigma'^2)$, and is independent of $\Xs$, where we denote $\sigma'^2 = \sigma^2 + \frac{M^{2(k-s+1)}-1}{n(M^2-1)}$. So we have 
	\begin{align*}
	\VAR(Y_i\big|\Xs) &= \E\left[ (e_i^T\cdot\left(\Xs^T\Xs\right)^{-1} \Xs^T\left(\Xss\bss + \bz\right))^2\big|\Xs\right] \notag\\
	&=e_i^T\left(\Xs^T\Xs\right)^{-1} \Xs^T\E\left[ (\Xss\bss + \bz)(\Xss\bss + \bz)^T\big|\Xs\right] \Xs \left(\Xs^T\Xs\right)^{-1} e_i \notag\\
	&= e_i^T\left(\Xs^T\Xs\right)^{-1} \Xs^T
	(\sigma'^2 \bI_{S_0})
	\Xs \left(\Xs^T\Xs\right)^{-1} e_i \notag\\
	&=\sigma'^2e_i^T\left(\Xs^T\Xs\right)^{-1} e_i.
	\end{align*}
	Thus by (\ref{eq:wishart}) again, we obtain
	\begin{align*}
	\E(V_i) &= \E(\sigma'^2e_i^T\left(\Xs^T\Xs\right)^{-1} e_i)\notag\\
	& = \sigma'^2e_i^T \frac{n}{n-s-1}\bI_{S_0}^{-1} e_i\notag\\
	&= \frac{n\sigma'^2}{n-s-1}.
	\end{align*}
	
	To prove part (c), we use the formula for the second moment of the inverse Wishart distribution in the part (2) of Lemma \ref{lem:wishart}. With $E_i = \E(Y_i\big|\Xs)$, we get
	\begin{align*}
	\E(E_i^2) = \E(\E(Y_i\big|\Xs))^2 = &\frac{\lambda_{n}^{2} n^{2}}{(n-s)(n-s-3)} \left[\left(e_{i}^{T}\left(\frac{1}{n} \bI_{S_0}\right)^{-1} \zh_{S_0}\right)^{2} \right. \notag\\   & \left.+\frac{1}{n-s-1}\left(\zh_{S_0}^{T}\left(\frac{1}{n} \bI_{S_0}\right)^{-1} \zh_{S_0}\right)\left(e_{i}^{T}\left(\frac{1}{n} \bI_{S_0}\right)^{-1} e_{i}\right)\right]\notag \\
	=&\frac{\lambda_{n}^{2} n^{2}}{(n-s)(n-s-3)}\left[ n^2 +\frac{1}{n-s-1} \cdot ns \cdot n \right] \notag\\
	=&\frac{\lambda_{n}^{2} n^{4}(ns+n-s^2-2s+1)}{(n-s)(n-s-3)(n-s-1)}.
	\end{align*}
	Thus, we have
	\begin{align}
	\VAR(E_i) 	=&\frac{\lambda_{n}^{2} n^{4}(ns+n-s^2-2s+1)}{(n-s)(n-s-3)(n-s-1)} - \frac{\lambda_{n}^{2} n^{4}}{(n-s-1)^2}\notag\\=& \frac{\lambda_{n}^{2} n^{4}(n-1)}{(n-s)(n-s-3)(n-s-1)}.
	\end{align}
	
	By Chebyshev's inequality, we can get
	\begin{align}
	\P(|E_i| \ge (1 + \sqrt{n}) \E(E_i)) &\le \P(|E_i-\E(E_i)|\ge \sqrt{n}\E(E_i)) \notag\\
	&\le \frac{\VAR(E_i)}{n(\E(E_i))^2} \notag\\
	& = \frac{ns+n-s^2-2s+1}{4n(n-s)(n-s-3)} \notag\\
	& \le \frac{K_1}{n-s},\label{eq:K1}
	\end{align}
	for some constant $K_1$ when $n$ is large enough.
	
	Similarly, by Lemma \ref{lem:wishart} (2) again for $i = j$, and $\Sigma = \frac{1}{n}\bI$, we have
	\begin{align*}
	\VAR(V_i^2) = & \sigma'^4 \VAR\left[(e_i^T\left(\Xs^T\Xs\right)^{-1} e_i)^2\right] \\
	=& \sigma'^4 \frac{(n-s+1+n-s-1)n^2}{(n-s)(n-s-1)^2(n-s-3)}
	\\=&\frac{2\sigma'^{4} n^{2}}{(n-s)(n-s-1)(n-s-3)},
	\end{align*}
	and thus
	\begin{align}
	\P(V_i \ge 2E(V_i)) = \P(V_i - \E(V_i)\ge E(V_i))\le&\frac{\VAR(V_i)}{(\E(V_i))^2}\notag\\
	=& \frac{\frac{2\sigma'^{4} n^{2}}{(n-s)(n-s-1)(n-s-3)}}{\left(\frac{n\sigma'^2}{n-s-1}\right)^2}\notag\\
	\le& \frac{K_2}{n-s},\label{eq:K2}
	\end{align}
	for some constant $K_2$ for large $n$.
	Therefore combining (\ref{eq:K1},\ref{eq:K2}) with union bound, the statement in part 3 holds with $K = K_1+K_2$.
\end{proof}

\begin{lemma}\label{lem:gaussian_maxima}
	Let $(X_1, \dots,X_n)$ be independent and normally distributed.  We have 
	\begin{equation*}
	\E[\max_{1\le i\le n} |X_i|] \le 3 \sqrt{\log n} \max_{1\le i \le n} \sqrt{\E X_i^2}.
	\end{equation*}
\end{lemma}
\begin{proof}[Proof of Lemma \ref{lem:gaussian_maxima}]
	This is a well-known result of Gaussian maxima. We omit its proof and refer the reader to, for example, \cite{wainwright2009information} Lemma 9.
\end{proof}
\begin{lemma}\label{lem:normalincrease}
	Let $Y\sim \Normal(\mu, \sigma^2)$. Suppose $\mu\le \mu_0$, $\sigma\le \sigma_0$, and $\rho \ge \mu_0$, then the probability $\P(Y\ge \rho) \le \P(Z\ge \rho)$, where $Z\sim \Normal(\mu_0.\sigma_0)$.
\end{lemma}
\begin{proof}[Proof of Lemma \ref{lem:normalincrease}]
	By definition, we have
	$$
	\P(Y\ge \rho) = 1-\Phi\left(\frac{\rho - \mu}{\sigma}\right) = \Phi\left(\frac{\mu-\rho }{\sigma}\right).
	$$
	And similarly, we have
	$$
	\P(Z\ge \rho) = \Phi\left(\frac{\mu_0 -\rho}{\sigma_0}\right).
	$$
	By the assumption, we know that 
	$$
	\frac{\mu-\rho }{\sigma} \le \frac{\mu_0 -\rho}{\sigma_0},
	$$
	and thus the lemma follows from the fact that $\Phi$ is increasing.
\end{proof}

\begin{lemma} \label{lem:sign_consistence}
Under the working assumptions, the event $E = \{\zh_{S_0} = sign(\bbeta_{S_0})\}$ satisfies 
$$
\P(E^C) = o_{n}(1).
$$
\end{lemma}
\begin{proof}[Proof of Lemma \ref{lem:sign_consistence}]
This is the classic result of Theorem 3 in \cite{wainwright2009information}. For the Lasso problem considered in \Cref{eq:oracle_lasso_problem}, with the identity covariance matrix, all conditions (26a), (26b), (26c) therein are easily satisfied with $C_{\min} = C_{\max} = 1$. As long as we have the condition
\begin{equation}\label{eq:beta_condition_sign}
\min_{j = 1, .., s} \beta_j > g(\lambda) = c_3 \lambda + 20 \sqrt{\frac{\sigma^2\log s }{n}},
\end{equation}
for some constant $c_3>0$, we can guarantee that $sign(\bbeta_{S_0}) = sign(\hat{\bbeta}_{S_0})\equiv \zh_{S_0}$ with probability $1 - c_1\exp(-c_2 \min\{s, \log(k - s)\})$ for some constant $c_1, c_2 > 0$. 

To verify condition (\ref{eq:beta_condition_sign}), recall that $\lambda = n^b$, $s = O(n / \log p)$, and $\min_{j = 1, .., s} \beta_i = M(n)^{k - s + 1}$. Because $M(n) = n^a$ and $b< (k - s + 1)a - \frac{3}{2}$, we have
\begin{align*}
\min_{j = 1, .., s} \beta_j = M(n)^{k - s + 1} &= n^{(k - s + 1) a}\\
&\ge n^{b + \frac{3}{2}}\\
& = \lambda n^{\frac{3}{2}}\\
&\ge c_3 \lambda + 20 \sqrt{\frac{\sigma^2\log s }{n}} = g(\lambda).
\end{align*}
for sufficiently large $n$. This implies $P(E)\to 1$, or $P(E^C) = o_n(1)$.
\end{proof}

\subsection{Technical proofs for \Cref{sec:proofs}}\label{sec:additional_proofs}

\subsubsection{A property of FDP and TPP: any trade-off curve is strictly increasing}\label{sec:fundamental_prop_proof}

A natural belief on the pair of $(\TPP, \FDP)$ is that FDP (type-I error) should increase with TPP (power), which may be strengthened by our simulation plots. However, along a single Lasso path, this is in general not necessarily true. It is well-known that Lasso is not monotone \citep{donoho2009observed}, so it is possible that with more and more true variables entering the Lasso path, fewer and fewer noise variables retain in the Lasso path. In such a case, FDP is no longer a monotone function of TPP. This possibility complicates our analysis, yet the following lemma asserts that this possibility is impossible. We prove that the asymptotic FDP is strictly increasing with the asymptotic TPP. Formally speaking, as $\lambda$ varies, $\fdp_\lambda$ can be seen as a function of $\tpp_\lambda$, and $\fdp_\lambda$ is a strictly increasing function of $\tpp_\lambda$.

To be rigorous, in the following lemma---indeed throughout the paper---we consider the regime below the Donoho--Tanner phase transition. We refer interested readers to \cite{wang2020bridge} for results above this phase transition.


\begin{lemma}\label{lem:fundamental_tppfdp}
	Fix $\epsilon, \delta, \sigma$, and $\Pi\ne0$. We have that $\fdp_\lambda(\Pi)$ is a strictly increasing function of $\tpp_\lambda(\Pi)$. That is, $\fdp_\lambda(\tpp_\lambda)$ is a well-defined function, and $\textnormal{fdp}_\lambda^{\infty\prime}(\cdot)\bigr|_{\tpp_\lambda}>0$ for any valid value of $\tpp_\lambda$.\footnote{As we will see in Lemma \ref{lem:lemD4}, the valid range of $\tpp_\lambda$ is the range $(0, u^\star)$. In this paper, we only focus on the case where $u^\star = 1$.}
\end{lemma}

To prove this lemma, we need the following characterizations among $\alpha$, $\lambda$, $\fdp_\lambda$ and $\tpp_\lambda$. 
\begin{lemma}\label{lem:Foundation}
	Fix $\epsilon, \delta, \sigma$, and $\Pi\ne0$. Consider any $\alpha$, $\tau$, $\lambda$ that solve equations (\ref{basic}). We have the following facts
	\begin{eqnarray}
	&&\frac{\mathrm{d}\alpha}{\mathrm{d}\lambda}>0 \label{gen1}\\
	&&\frac{\mathrm{d} \tpp}{\mathrm{d}\alpha}<0 \label{gen2}\\
	&&\frac{\mathrm{d}\fdp}{\mathrm{d}\alpha}<0\label{gen2.5}\\
	&&\frac{\mathrm{d}\tpp}{\mathrm{d}\fdp}>0 \label{gen3}
	\end{eqnarray}
\end{lemma}

We note that the denotations of $\fdp$ and $\tpp$ stand for $\fdp_\alpha$ and $\tpp_\alpha$, where we treat $\alpha$ as the free parameter. We often suppress this dependence on $\alpha$ in the following proof when it is clear from the context, and use the denotations of $\fdp$ and $\tpp$ for simplicity. 
\begin{proof}[Proof of Lemma \ref{lem:Foundation}]
	
	The (\ref{gen1}) is a well-known result, and one can refer to, for example, Lemma 4.11 of \cite{ mousavi2018consistent} for a proof.
	
	To prove (\ref{gen2}), we note that $\tpp = \mathbb{P}(|\Pi^{\star}+\tau W|>\alpha\tau)$, where $\Pi^{\star}$ is the distribution of an entry of $\beta$ given it's not zero. For any $\Pi$ that is a proper distribution and satisfies (\ref{basic}), proving $\frac{\mathrm{d}}{\mathrm{d}\alpha} \mathbb{P}(|\Pi+\tau W|>\alpha\tau)<0$ will suffice. And this result follows from Lemma 4.10 of 	\cite{mousavi2018consistent}.
	Now we left to prove (\ref{gen2.5}) and (\ref{gen3}). We note that, however, (\ref{gen3}) follows directly from (\ref{gen2.5}), and therefore we only need to prove (\ref{gen2.5}). To see this fact, we note that $\text{tpp}^\infty$ is a strictly increasing function of $\alpha$, so $\alpha$ is also a function of $\text{tpp}^\infty$. By the chain rule, we have 
	$$\frac{\mathrm{d}\text{fdp}^\infty}{\mathrm{d}\text{tpp}^\infty} =\frac{\mathrm{d}\text{fdp}^\infty}{\mathrm{d}\alpha}\cdot\frac{\mathrm{d}\alpha}{\mathrm{d}\text{tpp}^\infty}.$$
	Note that we have already proven $\frac{\mathrm{d}\text{tpp}^\infty}{\mathrm{d}\alpha}<0$, which implies $\frac{\mathrm{d}\alpha}{\mathrm{d}\text{tpp}^\infty}<0$. Therefore, we only need to show (\ref{gen2.5}) is true to prove (\ref{gen3}).
	
	Now, to prove (\ref{gen2.5}), we observe that
	\begin{eqnarray}
	\text{fdp}^\infty(\alpha)&=&\frac{1}{1+\frac{\epsilon\mathbb{P}\left(\left|\frac{\Pi^{\star}}{\tau}+ W\right|>\alpha\right)}{2(1-\epsilon)\Phi(-\alpha)}}, \notag\\
	\frac{\mathrm{d}\text{fdp}^\infty}{\mathrm{d}\alpha}&=&\frac{\epsilon}{2(1-\epsilon)}\cdot\frac{\frac{\mathrm{d}\left(\frac{\mathbb{P}\left(\left|\frac{\Pi^{\star}}{\tau}+ W\right|>\alpha\right)}{\Phi(-\alpha)}\right)}{\mathrm{d}\alpha}} {\left(1+\frac{\epsilon\mathbb{P}\left(\left|\frac{\Pi^{\star}}{\tau}+ W\right|>\alpha\right)}{2(1-\epsilon)\Phi(-\alpha)}\right)^2}.
	\end{eqnarray}
	Since the denominator is positive, we only need to show the numerator is negative. For simplicity, we will abuse the notation a little bit by using $\Pi$ for $\Pi^\star$ in the rest of the proof. We need to show for all $\Pi\ne 0$, 
	$$
	\frac{\mathrm{d}\left(\frac{\mathbb{P}\left(\left|\frac{\Pi}{\tau}+ W\right|>\alpha\right)}{\Phi(-\alpha)}\right)}{\mathrm{d}\alpha}<0,
	$$
	or equivalently,
	\begin{equation}\label{dstar0}
	\frac{\mathrm{d}\mathbb{P}\left(\left|\frac{\Pi}{\tau}+ W\right|>\alpha\right)}{\mathrm{d}\alpha}\cdot \Phi(-\alpha)+\mathbb{P}\left(\left|\frac{\Pi}{\tau}+ W\right|>\alpha\right)\cdot \phi(\alpha)>0.
	\end{equation}
	We observe that
	\begin{eqnarray}
	&&\mathbb{P}\left(\left|\frac{\Pi}{\tau}+ W\right|>\alpha\right)\notag\\
	&=& \mathbb{E}\left[\mathbb{P}_W\left(\left.W>\alpha-\frac{\Pi}{\tau}\right|\Pi\right)+\mathbb{P}_W\left(\left.W<-\alpha-\frac{\Pi}{\tau}\right|\Pi\right)\right] \notag\\
	&=& \mathbb{E}\left[\Phi\left(\frac{\Pi}{\tau}-\alpha\right)+\Phi(-\frac{\Pi}{\tau}-\alpha)\right].\label{dstar1}
	\end{eqnarray}
	
	Substituting the expressions of (\ref{dstar1}) and (\ref{star8}) into (\ref{dstar0}), we obtain
	\begin{eqnarray}
	&&\frac{\mathrm{d}\mathbb{P}\left(\left|\frac{\Pi}{\tau}+ W\right|>\alpha\right)}{\mathrm{d}\alpha}\cdot \Phi(-\alpha)+\mathbb{P}\left(\left|\frac{\Pi}{\tau}+ W\right|>\alpha\right)\cdot \phi(\alpha)\notag\\
	&=&\frac{1}{\frac{\sigma^2\delta}{\tau^3}+\mathbb{E}_{\Pi,W}\left[\frac{\Pi^2}{\tau^3}1_{\{-\alpha<\frac{\Pi}{\tau}+W<\alpha\}}\right]}\left[
	\overset{\Omega_1}{\overbrace{\mathbb{E}_\Pi\left[\Phi\left(\frac{\Pi}{\tau}-\alpha\right)+\Phi(-\frac{\Pi}{\tau}-\alpha)\right]\cdot\frac{\sigma^2\delta}{\tau^3}\cdot\phi(\alpha)}} \right. \notag \\
	&&+ \overset{\Omega_2}{\overbrace{\mathbb{E}_\Pi\left[\Phi\left(\frac{\Pi}{\tau}-\alpha\right)+\Phi(-\frac{\Pi}{\tau}-\alpha)\right]\cdot \mathbb{E}_{\Pi,W}\left[\frac{\Pi^2}{\tau^3}1_{\{-\alpha<\frac{\Pi}{\tau}+W<\alpha\}}\right]\cdot\phi(\alpha)}} \notag \\
	&& \overset{\Omega_3}{\overbrace{-\mathbb{E}_\Pi\left[\phi\left(\alpha-\frac{\Pi}{\tau}\right)+\phi\left(\alpha+\frac{\Pi}{\tau}\right)\right]\cdot\frac{\sigma^2\delta}{\tau^3}\cdot\Phi(-\alpha)}} \notag \\
	&& \overset{\Omega_4}{\overbrace{-\mathbb{E}_\Pi\left[\phi\left(\alpha-\frac{\Pi}{\tau}\right)+\phi\left(\alpha+\frac{\Pi}{\tau}\right)\right]\cdot \mathbb{E}_{\Pi,W}\left[\frac{\Pi^2}{\tau^3}1_{\{-\alpha<\frac{\Pi}{\tau}+W<\alpha\}}\right]\cdot\Phi(-\alpha)}}\notag \\
	&&+\overset{\Omega_5}{\overbrace{\mathbb{E}_\Pi\left[\frac{\Pi}{\tau^2}\left(-\phi\left(\alpha-\frac{\Pi}{\tau}\right)+\phi\left(\alpha+\frac{\Pi}{\tau}\right)\right)\right]
			\cdot \mathbb{E}_\Pi\left[\alpha\left(\int^{\infty}_{\alpha-\frac{\Pi}{\tau}} \phi(w)\mathrm{d}w +
			\int_{-\infty}^{-\alpha-\frac{\Pi}{\tau}} \phi(w)\mathrm{d}w\right)\right]\cdot\Phi(-\alpha)}}\notag \\
	&&+\left.\overset{\Omega_6}{\overbrace{\mathbb{E}_\Pi\left[\frac{\Pi}{\tau^2}\left(-\phi\left(\alpha-\frac{\Pi}{\tau}\right)+\phi\left(\alpha+\frac{\Pi}{\tau}\right)\right)\right]
			\cdot \mathbb{E}_\Pi\left[\int^{\infty}_{\alpha-\frac{\Pi}{\tau}} -w \phi(w)\mathrm{d}w +
			\int_{-\infty}^{-\alpha-\frac{\Pi}{\tau}} w\phi(w)\mathrm{d}w\right]\cdot\Phi(-\alpha)}}\right]. \notag\\
	&&\label{dstar1.5}
	\end{eqnarray}
	Note the denominator is positive, so we only need to prove that the numerator is positive. 
	Let $g(u) = (\Phi(u-\alpha)+\Phi(-u-\alpha))\phi(\alpha)-(\phi(u-\alpha)+\phi(u+\alpha))\Phi(-\alpha)$. By Lemma \ref{lem3}, $g(u)>0$ for $u\ne0$, and therefore we have
	\begin{eqnarray}\label{dstar2}
	\Omega_1+\Omega_3=\frac{\sigma^2\delta}{\tau^3}\mathbb{E}_\Pi \left[g\left(\frac{\Pi}{\tau}\right)\right]>0.
	\end{eqnarray}
	For $\Omega_5$, we observe that if $\Pi>0$, then $-\phi\left(\alpha-\frac{\Pi}{\tau}\right)+\phi\left(\alpha+\frac{\Pi}{\tau}\right)<0$; and if $\Pi<0$, then $-\phi\left(\alpha-\frac{\Pi}{\tau}\right)+\phi\left(\alpha+\frac{\Pi}{\tau}\right)>0$. Therefore, we have
	$$
	\mathbb{E}_\Pi\left[\frac{\Pi}{\tau}\left(-\phi\left(\alpha-\frac{\Pi}{\tau}\right)+\phi\left(\alpha+\frac{\Pi}{\tau}\right)\right)\right]\leq0
	$$
	So, by the fact that $\Phi(-\alpha)~\leq~ \frac{\phi(\alpha)}{\alpha}$, the definition of $\Phi(x)$, and (\ref{star8lem}), we have
	\begin{align}
	&~~~~\Omega_5\notag\\
	&= \mathbb{E}_\Pi\left[\frac{\Pi}{\tau^2}\left(-\phi\left(\alpha-\frac{\Pi}{\tau}\right)+\phi\left(\alpha+\frac{\Pi}{\tau}\right)\right)\right]
	~\mathbb{E}_\Pi\left[\left(\int^{\infty}_{\alpha-\frac{\Pi}{\tau}} \phi(w)\mathrm{d}w +
	\int_{-\infty}^{-\alpha-\frac{\Pi}{\tau}} \phi(w)\mathrm{d}w\right)\right]\cdot\alpha\Phi(-\alpha)\notag\\
	&\geq \mathbb{E}_\Pi\left[\frac{\Pi}{\tau^2}\left(-\phi\left(\alpha-\frac{\Pi}{\tau}\right)+\phi\left(\alpha+\frac{\Pi}{\tau}\right)\right)\right]
	~\mathbb{E}_\Pi\left[\left(\int^{\infty}_{\alpha-\frac{\Pi}{\tau}} \phi(w)\mathrm{d}w +
	\int_{-\infty}^{-\alpha-\frac{\Pi}{\tau}} \phi(w)\mathrm{d}w\right)\right]\cdot\phi(\alpha)\notag\\
	&=\mathbb{E}_\Pi\left[\frac{\Pi}{\tau^2}\left(-\phi\left(\alpha-\frac{\Pi}{\tau}\right)+\phi\left(\alpha+\frac{\Pi}{\tau}\right)\right)\right]
	~\mathbb{E}_\Pi\left[\Phi\left(\frac{\Pi}{\tau}-\alpha\right)+\Phi\left(-\frac{\Pi}{\tau}-\alpha\right)\right]\cdot\phi(\alpha)\notag\\
	&=\mathbb{E}_\Pi\left[\frac{\Pi}{\tau^2}\int^{\alpha-\frac{\Pi}{\tau}}_{-\alpha-\frac{\Pi}{\tau}} w\phi(w)\mathrm{d}w\right] ~\mathbb{E}_\Pi\left[\Phi\left(\frac{\Pi}{\tau}-\alpha\right)+\Phi\left(-\frac{\Pi}{\tau}-\alpha\right)\right]\cdot\phi(\alpha).\label{dstar2.5}
	\end{align}
	Similarly, by (\ref{dstar1.5}) and (\ref{dstar2.5}), and then by the definition of $f(u)$ in Lemma \ref{lem2}, we obtain
	\begin{align}
	&~~~~\Omega_2+\Omega_5\notag\\
	&\geq \mathbb{E}_\Pi\left[\Phi\left(\frac{\Pi}{\tau}-\alpha\right)+\Phi\left(-\frac{\Pi}{\tau}-\alpha\right)\right]~\mathbb{E}_\Pi\left[\frac{\Pi}{\tau^2}
	\int^{\alpha-\frac{\Pi}{\tau}}_{-\alpha-\frac{\Pi}{\tau}} w\phi(w)\mathrm{d}w+\frac{\Pi^2}{\tau^3} \int^{\alpha-\frac{\Pi}{\tau}}_{-\alpha-\frac{\Pi}{\tau}}
	\phi(w)\mathrm{d}w\right]\cdot\phi(\alpha)\notag\\
	&=\mathbb{E}_\Pi\left[\Phi\left(\frac{\Pi}{\tau}-\alpha\right)+\Phi\left(-\frac{\Pi}{\tau}-\alpha\right)\right]~\mathbb{E}_\Pi\left[\frac{1}{\tau}~ f\left(\frac{\Pi}{\tau}\right)\right]\cdot\phi\left(\alpha\right).\label{dstar3}
	\end{align}
	Similar to the proof of equation (\ref{star9}), we have
	\begin{equation}
	\Omega_4+\Omega_6 = -\mathbb{E}_\Pi\left[\left(\alpha-\frac{\Pi}{\tau}\right)+\phi\left(\alpha+\frac{\Pi}{\tau}\right)\right]\mathbb{E}_\Pi\left[\frac{1}{\tau}f\left(\frac{\Pi}{\tau}\right)\right].\label{dstar4}
	\end{equation}
	Combining the last display, (\ref{dstar3}), (\ref{dstar4}), Lemma \ref{lem2} and the well-known result that $\Phi(-\alpha)~\leq~ \frac{\phi(\alpha)}{\alpha}$, we obtain
	\begin{align}\label{eq:A51}
	\Omega_2+\Omega_4+\Omega_5+\Omega_6~~\geq~~\mathbb{E}_\Pi\left[g\left(\frac{\Pi}{\tau}\right)\right]\cdot\mathbb{E}_\Pi\left[\frac{1}{\tau}f\left(\frac{\Pi}{\tau}\right)\right]~~>~0,
	\end{align}
	Put together (\ref{dstar2}), (\ref{dstar1.5}) and (\ref{eq:A51}), we have for all $\Pi\neq0$
	$$\frac{\mathrm{d}\mathbb{P}\left(\left|\frac{\Pi}{\tau}+ W\right|>\alpha\right)}{\mathrm{d}\alpha}\cdot \Phi(-\alpha)+\mathbb{P}\left(\left|\frac{\Pi}{\tau}+ W\right|>\alpha\right)\cdot \phi(\alpha)>0,$$
	which, by (\ref{dstar0}), amounts to (\ref{gen2.5}), or
	$$
	\frac{\mathrm{d}\fdp}{\mathrm{d}\alpha}<0. 
	$$
	Therefore, combining with (\ref{gen2}), we obtain that
	$$\frac{\mathrm{d}\tpp}{\mathrm{d}\fdp}>0. $$
\end{proof}

Summarizing the result we have proven, it is very easy to prove Lemma \ref{lem:fundamental_tppfdp}.

\begin{proof}[Proof of Lemma \ref{lem:fundamental_tppfdp}]
	Observe $\tpp(\alpha)$ is a strictly increasing function of $\alpha$, and thus $\tpp$ is a one-to-one function of $\alpha$. The inverse function therefore exists, so $\alpha$ is a strictly increasing function of $\tpp$. Similarly, $\fdp$ is also a strictly increasing function of $\alpha$. Therefore, we conclude that $\fdp = \fdp(\alpha) = \fdp(\alpha(\tpp)) = \fdp(\tpp)$ is a strictly increasing function of $\tpp$, and that $\frac{\mathrm{d}\fdp}{\mathrm{d}\tpp}>0$ holds for any valid value of $\tpp$.
\end{proof}

Now, we prove the lemmas that we have used in the proof of Lemma \ref{lem:Foundation}.

\begin{lemma}\label{lem2}
	Let $f(u) = u\int_{-\alpha-u}^{\alpha-u} (w+u)\phi(w)\mathrm{d}w$. We have $f(u)> 0$, for all $u\neq0 \in \mathbb{R}$.
\end{lemma}

\begin{proof}[Proof of Lemma \ref{lem2}]
	Observe that
	\begin{eqnarray*}
		f(u)&=& u\int_{-\alpha-u}^{\alpha-u} (w+u)\phi(w)\mathrm{d}w \\
		&\overset{w'=w+u}{=}& u\int_{-\alpha}^{\alpha} w'\phi(w'-u)\mathrm{d}w' \\
		&=&u\int_{0}^{\alpha} w'[\phi(w'-u)-\phi(-w'-u)]\mathrm{d}w'.
	\end{eqnarray*}
	So, if $u>0$, then $\phi(w'-u)-\phi(-w'-u)>0,$ for any $w'\in(0,\alpha]$, thus $f(u)>0$; and if $u<0$, then $\phi(w'-u)-\phi(-w'-u)<0,$ for any $w'\in(0,\alpha]$, thus $f(u)>0$.
\end{proof}

\begin{lemma}\label{lem3}
	For any fixed $\alpha>0$, let $g(u) = (\Phi(u-\alpha)+\Phi(-u-\alpha))\phi(\alpha)-(\phi(u-\alpha)+\phi(u+\alpha))\Phi(-\alpha)$. Then we have $g(u)\geq0$, and the equality $g(u)=0$ holds if and only if $u=0$.
\end{lemma}

\begin{proof}[Proof of Lemma \ref{lem3}]
	We observe that
	\begin{eqnarray*}
		g(u)&=& \phi(\alpha)\left(\int^{\infty}_{\alpha-u}\phi(w)\mathrm{d}w + \int^{\infty}_{\alpha+u}\phi(w)\mathrm{d}w -\int^{\infty}_{\alpha}\phi(w)\mathrm{d}w \left(\frac{\phi(u-\alpha)}{\phi(\alpha)}+\frac{\phi(u+\alpha)}{\phi(\alpha)}\right)\right) \notag \\
		&=& \phi(\alpha)\left(\int^{\infty}_{\alpha}\phi(w) \e^{wu}\cdot \e^{\frac{-u^2}{2}} \mathrm{d}w + \int^{\infty}_{\alpha}\phi(w) \e^{-wu}\cdot \e^{\frac{-u^2}{2}} \mathrm{d}w \right.\notag\\
		&& \left.-\int^{\infty}_{\alpha}\phi(w) (\e^{\alpha u}+\e^{-\alpha u})\cdot \e^{\frac{-u^2}{2}} \mathrm{d}w\right) \notag\\
		&=&\phi(\alpha)\e^{\frac{-u^2}{2}}\int^{\infty}_{\alpha}\phi(w) \left[(\e^{wu}+\e^{-wu}) -(\e^{\alpha u}+\e^{-\alpha u})\right] \mathrm{d}w.
	\end{eqnarray*}
	Since for any $w>\alpha>0$, we have
	$$\e^{wu}+\e^{-wu}> \e^{\alpha u}+\e^{-\alpha u}, ~\text{ for any } u\in \mathbb{R}.$$
	We obtain $g(u)\geq0$, and it is clear the equality holds if and only if $u=0$.
\end{proof}

\subsubsection{Miscellaneous proofs for \Cref{sec:q_upper}}\label{sec:upper_additional_proof}
In this section, we prove all the necessary lemmas for \Cref{thm:upper}. To start with, we state the following lemma that specifies the range of all valid $\tpp$'s.
\begin{lemma}\label{lem:lemD4}[Lemma C.1 and Lemma C.4 in \citet{su2017false}]
	Put
	\[
	u^{\star}(\delta, \epsilon):=\left\{\begin{array}{ll}
	1-\frac{(1-\delta)\left(\epsilon-\epsilon^{\star}\right)}{\epsilon\left(1-\epsilon^{\star}\right)}, & \delta<1 \text { and } \epsilon>\epsilon^{\star}(\delta), \\
	1, & \text { otherwise.}
	\end{array}\right.
	\]
	Then
	\[
	\operatorname{tpp}^{\infty}<u^{\star}(\delta, \epsilon).
	\]
	Moreover, any $u$ between $0$ and $u^\star$ can be uniquely realized as $\tpp$, by setting $\alpha = t^\star(u)$ which is the root to (\ref{eq:t_lower}).
\end{lemma}
From this lemma, we know for $\delta <1$ and large $\epsilon$, it is possible that the range of $\tpp$ is no longer $(0,1)$. In such a case, we are ``above the Donoho--Tanner phase transition (DTPT)''; and symmetrically, when $\tpp$ has the range $(0,1)$, we are ``below the DTPT''. The purpose of this lemma is mainly for the completeness of the theory. In the following, however, we will always assume the range of $\tpp$ is $(0,1)$ to avoid extra complicity. This assumption will simplify our argument, but the proofs of the theorems can be extended to the case when this assumption is not true.

Now, we prove that the upper curve can be achieved by any $(\epsilon, M)$-homogeneous prior (\ref{eq:least}). This implies that the homogeneous effect sizes are the least desired.
\begin{lemma}\label{lem:upper}
	Given $(\epsilon, \delta)$ and $\sigma=0$. Any $(\epsilon, M)$-homogeneous prior gives the same unique trade-off curve $q^\nabla$ on $(0, u^{\star})$. Furthermore, this curve has the expression specified in (\ref{eq:q_upper}).
\end{lemma}

\begin{proof}[Proof of Lemma \ref{lem:upper}]
	
	We start with the proof to show the curve $q^\nabla$ is \textit{unique} in the sense that any two different $(\epsilon, M)$-homogeneous priors give the same trade-off curve. Consider any two $(\epsilon, M)$-homogeneous priors $\Pi_1$ and $\Pi_2$. Let their nonzero conditional priors be $\Pi_1^{\star}\equiv M_1$ and $\Pi_2^{\star}\equiv M_2$. Treating $\alpha>\alpha_0$ as the free parameter, we denote the solution to $\tau$ in equation (\ref{basic}) with prior $\Pi_1$ by $\tau_1$. We have
	\begin{equation*}
	\delta = (1-\epsilon)\E(\eta_{\alpha}(W)^2) + \epsilon \E\left(\eta_{\alpha}\left(\frac{M_1}{\tau_1}+W\right)-\frac{M_1}{\tau_1}\right)^2.
	\end{equation*}
	It is clear from a simple calculation that $ \tau_2 = \tau_1\frac{M_1}{M_2}$, $\alpha$ and $\Pi_2$ also solve the first equation in (\ref{basic}), that is,
	\begin{equation}
	\delta = (1-\epsilon)\E(\eta_{\alpha}(W)^2) + \epsilon \E\left(\eta_{\alpha}\left(\frac{M_2}{\tau_2}+W\right)-\frac{M_2}{\tau_2}\right)^2,
	\end{equation}
	which implies $\tau_2$ is the solution to (\ref{basic}) given $\alpha$ and prior $\Pi_2$.
	Observe the relationships $\tau_2 = \tau_1\frac{M_1}{M_2}$, $\Pi_1^{\star}\equiv M_1$ and $\Pi_2^{\star}\equiv M_2$. We have
	\begin{align*}
	\P\left(\left|\frac{\Pi_1^\star}{\tau_1}+ W\right|>\alpha\right) &=  \P\left(\left|\frac{\Pi_2^\star}{\tau_2}+ W\right|>\alpha\right).
	\end{align*}
	Therefore, combining the equality above with (\ref{eq:fdp_tpp_infty}), we obtain
	$$
	\tpp_{\alpha}(\Pi_1) = \P\left(\left|\frac{\Pi_1^\star}{\tau_1}+ W\right|>\alpha\right) = \P\left(\left|\frac{\Pi_2^\star}{\tau_2}+ W\right|>\alpha\right) =\tpp_{\alpha}(\Pi_2),
	$$
	and
	$$
	\fdp_{\alpha}(\Pi_1) =\frac{2(1-\epsilon)\Phi(-\alpha)}{2(1-\epsilon)\Phi(-\alpha) + \epsilon\P\left(\left|\frac{\Pi_1^\star}{\tau_1}+ W\right|>\alpha\right) } =\frac{2(1-\epsilon)\Phi(-\alpha)}{2(1-\epsilon)\Phi(-\alpha) + \epsilon\P\left(\left|\frac{\Pi_2^\star}{\tau_2}+ W\right|>\alpha\right) } =\fdp_{\alpha}(\Pi_1).
	$$
	This means that any point on $q^{\Pi_1}(\cdot)$ is also on $q^{\Pi_2}(\cdot)$. Similarly, any point on $q^{\Pi_2}(\cdot)$ is also on $q^{\Pi_1}(\cdot)$. By Lemma \ref{lem:fundamental_tppfdp}, they are both strictly increasing function, and thus must equal everywhere on the entire domain $(0,u^\star)$.

	Now, we proceed to prove that this unique trade-off curve has the expression given by (\ref{eq:q_upper}). Fix some $(\epsilon, M)$-homogeneous prior $\Pi^\nabla$. Let $u$ be some point between $0$ and $u^\star = 1$, and set $\alpha$ such that $\tpp_\alpha(\Pi^\nabla)$ equals to $u$. We have	
	$$
	u = \tpp_\alpha= \PP(|\Pi^{\nabla\star} + \tau W| > \alpha \tau) = \Phi(-\alpha+\widetilde{M})+\Phi(-\alpha-\widetilde{M}),$$
	where $\widetilde{M} = \frac{M}{\tau}$.
	Let $\varsigma = -\alpha + \widetilde{M}$, then the equation above becomes
	\begin{equation}\label{eq:alpha_gamma_u}
	\Phi(\varsigma) + \Phi(-2\alpha -\varsigma) = u.
	\end{equation}
	According to (\ref{basic}), we have
	\begin{equation}\label{eq:A53}
	\delta = (1-\epsilon) \EE[\eta_\alpha(W)]^2 + \epsilon\EE[\eta_\alpha(\widetilde{M}+W) - \widetilde{M}]^2.
	\end{equation}	
	By a simple algebraic fact
	$$
	\EE[\eta_\alpha(W)]^2 = 2[(1+\alpha^2)\Phi(-\alpha)-\alpha\phi(\alpha)],
	$$
	and the fact
	\begin{align*}
	\EE[\eta_\alpha(\widetilde{M}+W) - \widetilde{M}]^2 =& -(\alpha+\widetilde{M})\phi(\alpha-\widetilde{M})-(\alpha-\widetilde{M})\phi(\alpha+\widetilde{M}) +\\
	& (1+\alpha^2)[\Phi(-\alpha+\widetilde{M})+\Phi(-\alpha-\widetilde{M})] + \widetilde{M}^2[\Phi(\alpha-\widetilde{M})-\Phi(-\alpha-\widetilde{M})],
	\end{align*}
	we can plug-in the definition of $\varsigma$ into (\ref{eq:A53}) and obtain
	\begin{align}
	\delta = & 2(1-\epsilon)[(1+\alpha^2)\Phi(-\alpha)-\alpha\phi(\alpha)] + \epsilon[-(2\alpha+\varsigma)\phi(\varsigma) + \varsigma \phi(2\alpha+\varsigma) + \notag \\
	& (1+\alpha^2)[\Phi(\varsigma)+\Phi(-2\alpha-\varsigma)] + (\varsigma+\alpha)^2[\Phi(-\varsigma) + \Phi(-2\alpha-\varsigma)]]. \label{eq:alpha_gamma}
	\end{align}
	So $\varsigma = \varsigma(\alpha;\epsilon,\delta)$ is the solution of the equation above. 
	And finally combining the last equation with (\ref{eq:alpha_gamma_u}), we get an equation in $\alpha$ 
	\begin{equation}\label{eq:alpha_varsigma}
	\Phi(\varsigma(\alpha)) + \Phi(-2\alpha-\varsigma(\alpha)) = u,
	\end{equation}
	Denote the solution of $\alpha$ to the equation above by $t^\nabla = t^\nabla(u; \epsilon, \delta)$.	We have
	\begin{equation}\label{eq:q_nabla}
	q^\nabla(u;\epsilon,\delta) = \fdp_{t^\nabla}(\Pi^\nabla) = \frac{2(1-\epsilon)\Phi(-t^\nabla(u))}{2(1-\epsilon)\Phi(-t^\nabla(u)) + \epsilon u}.
	\end{equation}
	Therefore, the expression for the upper boundary is just defined by (\ref{eq:q_nabla}), where $t^\nabla$ is solved from (\ref{eq:alpha_gamma}) and (\ref{eq:alpha_varsigma}).
\end{proof}
We comment about the existence of $\alpha$ in the proof above. Note that both equations (\ref{eq:alpha_gamma}) and (\ref{eq:alpha_varsigma}) are derived from the AMP equations, which for any $\alpha>\alpha_0$, have unique solution $\tau$. Note that $\varsigma$ is a function of $\tau$, and thus it is also unique. Therefore, the solution to (\ref{eq:alpha_gamma}) also uniquely exists by Lemma \ref{lem:lemD4}.

\subsubsection{Miscellaneous proofs for \Cref{sec:lower_curve}}
In this section, we prove all necessary lemmas needed for \Cref{thm:lower}, and then prove \Cref{thm:lower}. We start by giving the proof to Lemma \ref{lem:n_points_convergence}. 
\begin{proof}[Proof of Lemma \ref{lem:n_points_convergence}]
	
	We treat $\tau$ as the free parameter instead of $\lambda$. To explicitly express the limiting process of $\bM$, we consider a sequence of priors $\{\Pi^\Delta(\bM^{(t)},\bm{\gamma})\}_t$, where $M_1^{(t)}\to\infty$ and $M_{i+1}^{(t)}/M_i^{(t)} \to \infty$ as $t\to\infty$. From (\ref{eq:fdp_tpp_infty}), the asymptotic TPP of $\Pi^\Delta(\bM^{(t)},\bm{\gamma})$ at $\tau$ can be written as
	\begin{equation}\label{eq:tppinftym}
	\begin{aligned}
	\tpp_\tau(\Pi^{\Delta}(\bM^{(t)},\bm{\gamma}))
	=&  ~\P \left(|\Pi^{\Delta\star}(\bM^{(t)},\bm{\gamma})+\tau W| >\alpha \tau\right)\\
	=& \left[ \gamma_1 \P\left(\left|W+\frac{M_1^{(t)}}{\tau}\right|  > \alpha\right)+\gamma_2 \P\left(\left|W+\frac{M_2^{(t)}}{\tau}\right|  > \alpha\right)\right.\\
	&\left.+\dots +\gamma_m \P\left(\left|W+\frac{M_m^{(t)}}{\tau} \right|  > \alpha\right) \right],
	\end{aligned}
	\end{equation}
	where $\alpha$ is solved from (\ref{basic}). We denote the last display by $\tppt_\tau$ for convenience.
	Similarly, we denote 
	\begin{equation}\label{eq:fdpinftym}
	\fdpt_\tau \equiv \fdp_\tau(\Pi^{\Delta}(\bM^{(t)},\bm{\gamma})) = \frac{2(1-\epsilon)\Phi(-\alpha)}{2(1-\epsilon)\Phi(-\alpha) + \epsilon \tppt(\tau)}.
	\end{equation}
	
	In the following proof, we will choose an $m$-tuple of $(\tau_i^{(t)})_{i = 1}^m$ for each fixed $t$, such that as $t\to \infty$, $(\tppt_{\tau_i^{(t)}}, \fdpt_{\tau_i^{(t)}}) \to (u_i, q^\Delta(u_i))$ at $m$ different $u_i$'s. This implies exactly that the limit of trade-off curves $q^{\Pi^\Delta(\bM^{(t)},\bm{\gamma})}$ agrees with $q^\Delta$ at (at least) different $m$ points in $(0,1]$.
	A natural way to pick such an $m$-tuple $(\tau_i^{(t)})_{i = 1}^m$ is 
	$$\tau_i^{(t)} = \sqrt{M_i^{(t)}M_{i+1}^{(t)}}, 1\le i\le m-1,$$
	and $\tau_m^{(t)} = m \times M_{m}^{(t)}$ when $i = m$.\footnote{In fact, one can pick any $\tau_i^{(t)}$ such that (\ref{eq:tau_regime}) and (\ref{eq:tau_regime_l}) hold.}  Under the regime of $M_1^{(t)}\to\infty$ and $M_{i+1}^{(t)}/M_i^{(t)}\to\infty$ for all $i$, we know $\tau_i^{(t)}$ satisfies 
	\begin{equation}\label{eq:tau_regime}
	|M_i^{(t)}| = o(\tau_i^{(t)}),\text{  and  }\tau_i^{(t)} = o(|M_{i+1}^{(t)}|) , ~~~1\le i\le m-1
	\end{equation}
	and 
	\begin{equation}\label{eq:tau_regime_l}
	|M_m^{(t)}| = o(\tau_m^{(t)}),
	\end{equation}
	as $t\to \infty$.
	Moreover, for any $\alpha>\alpha_0$ and any $1\le i\le m$, we have
	\begin{equation}\label{eq:P_asymp}
	\P\left(\left|W+\frac{M_j^{(t)}}{\tau_i^{(t)}}\right|>\alpha\right) = 
	\begin{cases}
	\P(|W|>\alpha) + o_t(1)  & \text{for } j= 1,2,\cdots, i,\\
	1 + o_t(1)  & \text{for } j = i+1,\cdots, m,\\
	\end{cases}
	\end{equation}
	and 
	\begin{equation}\label{eq:E_asymp}
	\eta_{\alpha}\left(W+\frac{M_j^{(t)}}{\tau_i^{(t)}}\right) - \frac{M_j^{(t)}}{\tau_i^{(t)}} = 
	\begin{cases}
	\eta_{\alpha}{(W)}+o_t(1) & \text{for } j = 1,2,\cdots, j,\\
	W - \alpha +o_{\P,t}(1) & \text{for } j = i+1,\cdots, m.\\
	\end{cases}
	\end{equation}
	Let $\gamma^{(j)} = \sum_{i = j+1}^m \gamma_i$ and use $\alpha_i^{(t)}$ to denote the solution of $\alpha$ to (\ref{basic}) given $\tau_i^{(t)}$. We have
	\begin{equation}\label{eq:mpointbasic}
	\left(1-\frac{\sigma^2}{\tau_i^{(t)2}}\right)\delta =  \epsilon \cdot  \underbrace{\E\left(\eta_{\alpha}\left(\frac{\Pi^\Delta(\bM^{(t)},\bm{\gamma})}{\tau_i^{(t)}} + W\right) - \frac{\Pi^\Delta(\bM^{(t)},\bm{\gamma})}{\tau_i^{(t)}}\right)}_{(*)} + (1-\epsilon) \E(\eta_{\alpha}(W))^2.
	\end{equation}
	By (\ref{eq:E_asymp}), the $(*)$ part of the last display is
	\begin{align*}
	(*) =& \sum_{j =1}^i \gamma_{j}\E\left(\eta_{\alpha}\left(\frac{M_j^{(t)}}{\tau_i^{(t)}} + W\right)-\frac{M_i^{(t)}}{\tau^{(t)}}\right)^2 +
	\sum_{j = i+1}^m \gamma_{j} \E\left(\eta_{\alpha}\left(\frac{M_j^{(t)}}{\tau_i^{(t)}} + W\right)-\frac{M_j^{(t)}}{\tau_i^{(t)}}\right)^2\\
	= &~ (1-\gamma^{(i)})\E(\eta_{\alpha}(W)^2) + \gamma^{(i)} \E((W- {\alpha})^2) + o_t(1).
	\end{align*}
	Observe the fact that $\sigma$ is fixed and thus $\frac{\sigma^2}{\tau_i^{(t)2}} \to 0$. With some simple calculation, (\ref{eq:mpointbasic}) can be written as
	\begin{equation}\label{eq:limitlinking}
	\epsilon \gamma^{(i)}(1+\alpha^2)+2(1-\epsilon \gamma^{(i)})[(1+\alpha^2)\Phi(-\alpha)-2\alpha\phi(\alpha)]=\delta + o_t(1).
	\notag
	\end{equation}
	Therefore the solution $\alpha_i^{(t)}$ of the equation above has a limit \footnote{By the existence asserted by AMP theory, the equation (\ref{eq:linking}) has a unique solution, denote it by $\alpha^{(i)}$, we know the solution $\alpha_i^{(t)}$ of (\ref{eq:limitlinking}) must converge into it, since the left-hand side of (\ref{eq:limitlinking}) is continuous in $\alpha$.}
	
	\begin{equation}\label{eq:alphalimit}
	\alpha_i^{(t)} \to \alpha^{(i)}, ~~\text{   as   } t\to\infty,
	\end{equation} 
	which solves the equation
	\begin{equation}\label{eq:linking}
	\epsilon \gamma^{(i)}(1+\alpha^2)+2(1-\epsilon \gamma^{(i)})[(1+\alpha^2)\Phi(-\alpha)-2\alpha\phi(\alpha)]=\delta.
	\end{equation}
	Note $\alpha^{(i)}$ is independent of the choice of $\{M^{(t)}\}_{t=1}^\infty$ and $(\tau_i^{(t)})_{i = 1}^m$.
	Direct calculation can verify that each solution $\alpha^{(i)}$ also satisfies the equation (\ref{eq:t_lower}) with setting
	\begin{equation}\label{eq:ui}
	u = u_i = 2\Phi(-\alpha^{(i)})(1-\gamma^{(i)})+\gamma^{(i)}.
	\end{equation}
	This implies $\alpha^{(i)}$ is also the unique solution of (\ref{eq:t_lower}), so 
	\begin{equation}\label{alpha}
	\alpha^{(i)}=t^\Delta(u_i).
	\end{equation}
	Combining (\ref{eq:tppinftym}), (\ref{eq:fdpinftym}), (\ref{eq:P_asymp}) and (\ref{eq:alphalimit}), the limits of $\tppt_{\tau_i^{(t)}}$ and $\fdpt_{\tau_i^{(t)}}$ are
	\begin{equation}\label{tdgamma}
	\begin{cases}
	\tppt_{\tau_i^{(t)}} &\to \ 2\Phi(-\alpha^{(j)})(1-\gamma^{(i)})+ \gamma^{(i)},\\
	\fdpt_{\tau_i^{(t)}} &\to \frac{2(1-\epsilon)\Phi(-\alpha^{(i)})}{2(1-\epsilon)\Phi(-\alpha^{(i)})+\epsilon(2\Phi(-\alpha^{(i)})(1-\gamma^{(i)})+\gamma^{(i)})}.
	\end{cases}
	\end{equation}
	By $\alpha^{(i)}=t^\Delta(u_i)$ and (\ref{eq:ui}), we obtain from (\ref{eq:q_lower}) that
	\begin{align}\label{eq:qlowerm}
	q^\Delta(u_i;\delta,\epsilon) &= \frac{2(1-\epsilon)\Phi(-t^\Delta(u_i))}{2(1-\epsilon)\Phi(-t^\Delta(u_i))+\epsilon u_i}\notag\\ &=\frac{2(1-\epsilon)\Phi(-\alpha^{(i)})}{2(1-\epsilon)\Phi(-\alpha^{(i)})+\epsilon(2\Phi(-\alpha^{(i)})(1-\gamma^{(i)})+\gamma^{(i)})}.
	\end{align}
	Combining (\ref{eq:fdp}), (\ref{eq:fdp_tpp_infty}), (\ref{eq:ui}), (\ref{tdgamma}) and (\ref{eq:qlowerm}), we finally obtain as $t\to\infty$,
	\begin{equation*}
	\begin{cases}
	\textnormal{tpp}_t^{\infty}(\tau_i^{(t)})\to u_i,\\
	\textnormal{fdp}_t^{\infty}(\tau_i^{(t)})\to q^\Delta(u_i;\delta,\epsilon).
	\end{cases}
	\end{equation*}
	Therefore, the limiting function of the trade-off curves of priors $\Pi^\Delta(\bM^{(t)},\bm{\gamma})$ agrees with the lower boundary $q^\Delta(\cdot;\delta,\epsilon)$ at $(u_i, q^\Delta(u_i)),~\text{for}~i=1,2,\cdots, m$. 
	Since the $m$ different points are nonzero, there must be at least $m-1$ points in $(0,1)$.
\end{proof}
An important set of equations is (\ref{eq:ui}). Recall that  (\ref{alpha}) asserts $\alpha^{(i)}=t^\Delta(u_i)$ for all $i$, and therefore the equations (\ref{eq:ui}) are, for all $i$,
\begin{equation}\label{eq:ui2}
u_i = 2\Phi(-t^\Delta(u_i))(1-\gamma^{(i)})+\gamma^{(i)}.
\end{equation}
The last display allows us to quantify the exact points $u_i$ the limit of $\Pi^\Delta(\bM^{(t)},\bm{\gamma})$ agrees with $q^\nabla$. This fact allows us to set $\bm{\gamma}$ cleverly so that the distance between any two consecutive $u_i$'s are small enough. This is formalized in the following lemma.
\begin{lemma}\label{lem:delta_close}
	For any $\xi>0$, there is some $\bm{\gamma} = \{\gamma_1, ...,\gamma_m\}$, with $\sum_i \gamma_i = 1$, such that the $m$ points specified by (\ref{eq:ui2}), together with $u_0 = 0$ and $u_{m+1} = 1$ \footnote{Technically speaking, it should be $u_{m+1} = u^\star$, yet as discussed earlier, we will focus on the case when we are below the Donoho--Tanner phase transition, so always $u^\star = 1$.} satisfy the following
	\begin{equation}\label{eq:u_close}
	\max_{1\le j\le m+1} |u_j-u_{j-1}| < \frac{\xi}{2}.
	\end{equation}
	
\end{lemma}
\begin{proof}[Proof of Lemma \ref{lem:delta_close}]
	We first prove that the difference $u_{m+1}-u_m = 0$. Since $\gamma^{(m)} = 1$, and by (\ref{eq:ui}), we obtain that $u_m = 1$, and thus $u_{m+1} - u_m = 1-1 = 0$.\footnote{One might want to verify the existence of $\alpha^{(m)}$. Since we always assume that we are below the DTPT, then for any $u<u^\star = 1$, the $\alpha = t^\Delta(u_i)$ exists as the solution to (\ref{eq:t_lower}) by Lemma \ref{lem:lemD4}. By setting $\gamma^{(m)} = 1$, one can directly verify this corresponds to set $u = u_m = 1^-$, and by the continuity of $t$ in equation (\ref{eq:t_lower}), we know $\alpha^{(m)} = t^\Delta(u_m)$ exists and less than infinite. And since all other $\alpha^{(i)}$'s are less than $\alpha^{(m)}$, they also exist.}
	With this in mind, we only need to prove the following two quantities can be arbitrarily small to ensure (\ref{eq:u_close}),
	\begin{equation}\label{eq:j=1}
	u_1-0 = \gamma_m + 2(1-\gamma_m)\Phi(-\alpha^{(m)}),
	\end{equation}
	and for all $m\ge j \ge2$,
	\begin{equation}\label{eq:j>1}
	u_j-u_{j-1} = \gamma_{m-j+1}+2(1-\gamma^{(j)})\Phi(-\alpha^{(j)})- 2(1-\gamma^{(j-1)})\Phi(-\alpha^{(j-1)}),
	\end{equation}
	where we remind the reader that by definition, $\gamma^{(m)} = \gamma_m$ and $\gamma^{(j)} - \gamma^{(j-1)} = \gamma_{m-j+1}$.
	
	For the expression in (\ref{eq:j>1}), we observe that 	
	\begin{align}
	u_j-u_{j-1} \le& ~\gamma_{m-j+1}|1+2(\Phi(-\alpha^{(j)})-\Phi(-\alpha^{(j-1)}))|\notag\\
	& +2(1-\gamma^{(j)})|\Phi(-\alpha^{(j)})-\Phi(-\alpha^{(j-1)})| \notag\\
	\le& 5\gamma_{m-j+1} + 2|\Phi(-\alpha^{(j)})-\Phi(-\alpha^{(j)}-\gamma_{m-j+1} ))|.\label{eq:a105}
	\end{align}
	We observe that in equation (\ref{eq:ui}) or (\ref{eq:ui2}), the dependence of $\alpha^{(j)}$ on $\gamma^{(j)}$ is only through linear functional of $\gamma^{(j)}$. Therefore, it is not hard to realize that $\alpha^{(j)}$ is continuous in $\gamma^{(j)}$. When 
	all $\{\gamma_s\}_{s>m-j+1}$ are fixed, the $\alpha^{(j)}$ is a continuous function in $\gamma_{m-j+1}$, and so is the expression in (\ref{eq:a105}). So we can pick $\gamma_{m-j+1}$ sufficiently small to ensure the (\ref{eq:a105}) is less than $\frac{\xi}{2}$.
	
	For the expression in (\ref{eq:j=1}), we pick some $M$ sufficiently large such that $\Phi(-M)<\frac{\xi}{8}$. By Lemma \ref{lem:cali_of_alpha}, we can pick $\gamma_m<\frac{\xi}{4}$ such that the solution to (\ref{eq:t_lower}) with $u$ being (\ref{eq:j=1}) satisfies that  $\alpha^{(m)}>M$, or $\Phi(-\alpha^{(m)}) <\Phi(-M)<\frac{\xi}{8}$. Therefore 
	$$
	\gamma_m + 2(1-\gamma^{(m)})\Phi(-\alpha(\gamma^{(m)})) <\frac{\xi}{4} + 2\cdot1\cdot\frac{\xi}{8} = \frac{\xi}{2}.	
	$$ 
\end{proof}

\begin{lemma}\label{lem:cali_of_alpha}
	For any fixed $\delta,~\epsilon$, $\xi>0$ and $M>0$. There exists $\gamma<\frac{\xi}{4}$, such that the solution $\alpha$ to (\ref{eq:t_lower}) with $u =  \gamma + 2(1-\gamma)\Phi(-\alpha)$ satisfies $\alpha >M$.
\end{lemma}

\begin{proof}[Proof of Lemma \ref{lem:cali_of_alpha}]
	We will use the following fact: there exists $\gamma<\frac{\xi}{4}$ and large $M''>M'>M$, such that:
	\begin{equation}\label{eq:a106}
	\begin{split}
	&(1-\epsilon\gamma)A(M')+\epsilon\gamma(1+M'^2)<\delta,\\
	&(1-\epsilon\gamma)A(M'')+\epsilon\gamma(1+M''^2)>\delta.
	\end{split}
	\end{equation}
	where $A(M) =2[(1+M^2)\Phi(-M)-M\phi(M)]$.
	
	Taken this as given for the moment, we set $u =  \gamma + 2(1-\gamma)\Phi(-\alpha)$ with $\gamma$ such that (\ref{eq:a106}) holds. Then (\ref{eq:t_lower}) becomes\footnote{Since the solution to $t$ is $\alpha$ here, we can plug in $u = \gamma + 2(1-\gamma)\Phi(-t)$.}
	$$
	\begin{aligned}
	\frac{2(1 - \epsilon)\left[ (1+t^2)\Phi(-t) - t\phi(t) \right] + \epsilon(1 + t^2) - \delta}{\epsilon\left[ (1+t^2)(1-2\Phi(-t)) + 2t\phi(t) \right]} = \frac{1 - \gamma - 2(1-\gamma)\Phi(-t)}{1 - 2\Phi(-t)}.
	\end{aligned}
	$$
	Observe the right hand side of the last display is just $1-\gamma$, so it is equivalent to
	\begin{equation*}
	2(1 - \epsilon)\left[ (1+t^2)\Phi(-t) - t\phi(t) \right] + \epsilon(1 + t^2) - \delta = (1-\gamma)\epsilon\left[ (1+t^2)(1-2\Phi(-t)) + 2t\phi(t) \right],
	\end{equation*}
	or,
	\begin{equation}\label{eq:a107}
	2(1 - \epsilon\gamma)\left[ (1+t^2)\Phi(-t) - t\phi(t) \right] + \epsilon\gamma (1 + t^2) - \delta = 0.
	\end{equation}
	By relationship (\ref{eq:a106}) and the fact that there exists unique solution $\alpha = t^\nabla(u)$ to (\ref{eq:t_lower}) and thus to (\ref{eq:a107}), we know the solution $\alpha\in (M', M'')$, and thus especially $\alpha > M$.
	
	Now, to prove (\ref{eq:a106}), we first note that it is direct to verify $A(t) = E[\eta_t(W)^2]$, and thus it is decreasing in $t$. And as $t\to\infty$, $A(t)\to 0$. Therefore for any $\delta>0$, we can pick $M'$ large enough such that $A(M') <\frac{\delta}{2}$, and now pick $\gamma\frac{\xi}{4}$ small enough such that $\epsilon \gamma (1+M'^2) < \frac{\delta}{2}$. Therefore, the left-hand side of the first equation of (\ref{eq:a106}) is bounded by $\delta$.
	For the second equation in (\ref{eq:a106}), pick $M''$ large enough so that the term $\epsilon \gamma (1+M''^2)>\delta$, and since $(1-\epsilon\gamma)A(M'')>0$, the second line also holds.
\end{proof}

The agreeing points asserted by Lemma \ref{lem:n_points_convergence} are close to each other in their x-coordinate distances. Therefore, by the uniform continuity of the lower curve $q^\nabla$ and Cantor's diagonalization argument, we can extend the result from Lemma \ref{lem:n_points_convergence} to uniform convergence.

\begin{lemma}\label{lem:uniform_convergence}
	There exist a sequence of prior of $\Pi^{(t)} = \Pi^\Delta(\bM^{(t)},\bm{\gamma}^{(t)})$, such that their trade-off curve $q^{\Pi^{(t)} }$ converge uniformly to $q^\Delta$ on any compact interval in $(0,1)$.
\end{lemma}

\begin{proof}[Proof of Lemma \ref{lem:uniform_convergence}]
	Fix any compact interval $\mathcal{I} = [a, b]$ in $(0,1)$. 
	As in Lemma \ref{lem:n_points_convergence}, we first consider prior $\Pi^{(t)} = \Pi^\Delta(\bM^{(t)},\bm{\gamma}^{(t)})$ with $\bm{\gamma}^{(t)} = \bm{\gamma}$ being some fixed $m$-tuple. By Lemma \ref{lem:fundamental_tppfdp}, we know that both $q^{\Pi^{(t)}}(u)$ and $q^\Delta(u)$ are continuous and strictly increasing. Consider any two adjacent agreeing points $u_j, u_{j+1}$ specified in (\ref{eq:ui}), such that $q^\Delta(u_j) = \lim_{t\to\infty}q^{\Pi^{(t)}}(u_j) $ and $q^\Delta(u_{j+1}) = \lim_{t\to\infty}q^{\Pi^{({t})}}(u_{j+1}) $. Since in the interval $(u_j,u_{j+1})$ the difference is controlled by
	\begin{equation*}
	q^{\Pi^{({t})}}(u)-q^\Delta(u) \le q^\Delta(u_{j+1})-q^\Delta(u_{j}), ~\text{for any} ~u\in (u_j,u_{j+1})
	\end{equation*}
	by the monotonicity of $q^{\Pi^{({t})}}$. This difference will be small as long as the gap $q^\Delta(u_{j+1})-q^\Delta(u_{j})$ is small, so we proceed to prove we can select $\Pi^{(t)}$ to ensure the gaps $q^\Delta(u_{j+1})-q^\Delta(u_{j})$ are small for all $i$.
	
	Fix any $\theta >0$. Since $q^\Delta$ is uniformly continuous on the compact set $\mathcal{I}$, there exists $\xi >0$ such that for any $u, v\in \mathcal{I}$,
	\begin{equation}\label{eq:4.38}
	|q^\Delta(u)-q^\Delta(v)|<\frac{\theta}{2}, ~\text{ as long as }|u-v|<\xi
	\end{equation} 
	By the proof of Lemma \ref{lem:n_points_convergence}, we can construct $\bm{\gamma}^{(t)} = \bm{\gamma}_\theta$ to be specified later, and $\bM^{(t)}_{\bm{\gamma}_\theta}$, such that the limit of $q^{\Pi^{(t)}}$ agrees with $q^\Delta$ at $m$ points $u_1, \cdots, u_m$. This implies there exists some $T_\theta$ such that for all $t \ge T_\theta$, 
	\begin{equation}\label{eq:4.39}
	\max_j \left|q^{\Pi^{({t})}}(u_j) - q^{\Delta}(u_{j})\right| <\frac{\theta}{2}.
	\end{equation}
	To specify the choice of  $\bm{\gamma}_\theta$, note that by Lemma \ref{lem:delta_close}, we can choose $\bm{\gamma}_\theta$ such that $u_1, .. ,u_m$ satisfies $u_1 - a < u_1 - 0 < \frac{\xi}{2}$, $u_m -b <\frac{\xi}{2}$ and
	\begin{equation*}
	\max_{2\le j\le m} |u_j-u_{j-1}| < \frac{\xi}{2}.
	\end{equation*}
	With this choice of $\bm{\gamma}_\theta$ together with (\ref{eq:4.38}) and (\ref{eq:4.39}), we obtain
	$$
	\sup_{t \ge T_\theta, u\in \mathcal{I}} \left|q^{\Pi^{({t})}}(u) - q^{\Delta}(u)\right| < \theta.
	$$
	Specifically, we have the equation above holds for $t = T_\theta$,
	$$
	\sup_{ u\in \mathcal{I}} \left|q^{\Pi^{({ T_\theta})}}(u) - q^{\Delta}(u)\right| < \theta.
	$$
	Since $\Pi^{(T_\theta)} = {\Pi^\Delta\left(\bM^{(T_\theta)}_{\bm{\gamma}_\theta}, \bm{\gamma}_\theta\right)}$, the inequality above is simply 
	\begin{equation}\label{eq:A79}
	\sup_{ u\in \mathcal{I}} \left|q^{\Pi^\Delta(\bM^{(T_\theta)}_{\bm{\gamma}_\theta}, \bm{\gamma}_\theta)}(u) - q^{\Delta}(u)\right| < \theta.
	\end{equation}
	Now we can apply Cantor's diagonalization trick since the last display is true for any $\theta>0$. We set $\theta_\zeta = \frac{1}{\zeta}\to 0, \zeta \ge 1$, and choose the priors $$\Pi^{(\zeta)} = \left\{\Pi^\Delta\left(\bM_{\gamma_{\theta_\zeta}}^{(T_{\theta_\zeta})},\bm{\gamma_{\theta_\zeta}}\right)\right\}_\zeta.
	$$
	Then we know from (\ref{eq:A79}) that $q^{\Pi^{(\zeta)}}$ converges to $q^\Delta$ uniformly on $\mathcal{I}$ as $\zeta\to\infty$.
\end{proof}

Now, we can proceed to prove \Cref{thm:lower}, whose proof is very similar to that of \Cref{thm:upper} in \Cref{pf:thm_upper}. 
\begin{proof}[Proof of \Cref{thm:lower} (a)]\label{pf:thm_lower}
	Consider any non-constant prior $\Pi$, we first prove that there exists some $\Pi^\Delta$ and $\nu > 0$ such that for all $c < \lambda, \lambda' < C$,
	\begin{equation}\label{eq:not_hold2}
	\tpp_{\lambda}(\Pi^\Delta) < \tpp_{\lambda'}(\Pi) + \nu \text{ and } \fdp_{\lambda}(\Pi^{\Delta}) > \fdp_{\lambda'}(\Pi) - \nu
	\end{equation}
	cannot hold simultaneously. To see this fact, first find $0 < u_1 < u_2 < 1$ such that the asymptotic powers $\tpp_{\lambda}(\Pi^\Delta), \tpp_{\lambda'}(\Pi)$ are always between $u_1$ and $u_2$ for $c < \lambda, \lambda' < C$. Next, we know from Lemma \ref{lem:lower_strict} that $q^\Delta$ is the strictly below any trade-off curve. So, for any prior $\Pi$, we have $q^\Delta(u)<q^\Pi(u)$ for any $u\in\mathcal{I} = [u_1, u_2]$. Note that both $q^\Delta$ and $q^\Pi$ are uniformly continuous on $\mathcal{I}$ and thus one can set $\nu' > 0$ to be
	\begin{equation}\label{eq:4.42}
	\nu' := \inf_{u_1 \le u_\le u_2} \left(q^\Pi(u) - q^\Delta(u)\right) > 0.
	\end{equation}
	Since $q^\Delta$ is a continuous function on the closed interval $[0, 1]$, we can make use of its uniform continuity, which ensures
	\begin{equation}\label{eq:4.43}
	\left| q^\Delta(u) - q^\Delta(u') \right| < \frac{\nu'}{4}
	\end{equation}
	as long as $|u - u'| \le \nu''$ for some $\nu'' > 0$. 
	By the assertion of Lemma \ref{lem:uniform_convergence}, we can choose a prior $\Pi^\Delta$ such that it is $\frac{\nu'}{4}$-close to $q^\Delta$ on $\mathcal{I}$,
	\begin{equation}\label{eq:4.44}
	\sup_{u_1\le u\le u_2} (q^{\Pi^\Delta}(u)-q^\Delta(u)) < \frac{\nu'}{4}.
	\end{equation}
	Now we can prove \eqref{eq:not_hold} cannot hold simultaneously with our choice of $\Pi^\Delta$ and $\nu = \min\{\nu'/2, \nu'' \}$. To see this, suppose we already have $\tpp_{\lambda}(\Pi^\Delta) < \tpp_{\lambda'}(\Pi) + \nu$. Now observe that
	\[
	\begin{aligned}
	\fdp_{\lambda}(\Pi^{\Delta}) &= q^{\Pi^\Delta}(\tpp_{\lambda}(\Pi^{\Delta})) \\
	&\le q^\Delta\left( \tpp_{\lambda}(\Pi^{\Delta}) \right) + \frac{\nu'}{4}\\
	&< q^\Delta\left( \tpp_{\lambda'}(\Pi^{\Delta}) + \nu \right) + \frac{\nu'}{4}\\
	& \le q^\Delta\left( \tpp_{\lambda'}(\Pi)\right) + \frac{\nu'}{2}\\
	& \le q^\Pi\left( \tpp_{\lambda'}(\Pi)\right) - \nu' + \frac{\nu'}{2}\\
	& = q^\Pi\left( \tpp_{\lambda'}(\Pi)\right) - \frac{\nu'}{2}\\
	& \le q^\Pi\left( \tpp_{\lambda'}(\Pi)\right) - \nu\\
	& = \fdp_{\lambda'}(\Pi) - \nu,
	\end{aligned}
	\]
	where the first inequality follows from (\ref{eq:4.44}); the second inequality follows from the fact that $\fdp(\tpp)$ is strictly increasing, and $\tpp_{\lambda}(\Pi^\Delta) < \tpp_{\lambda'}(\Pi) + \nu$; the third inequality is by (\ref{eq:4.43}); the fourth inequality is by (\ref{eq:4.42}); the last inequality is by the definition of $\nu$. 
	
	As such, the first inequality in \eqref{eq:not_hold2} leads to the violation of the second inequality.
	Having shown \eqref{eq:not_hold2}, it is easy to prove \Cref{thm:lower}. Lemma~\ref{lem:fdp_tpp_fix_lambda} ensures that both the following four terms
	\begin{equation}\label{eq:very_close2}
	\left| \TPP_\lambda(\Pi^\Delta) -  \tpp_\lambda(\Pi^\Delta) \right|, \left| \FDP_\lambda(\Pi^\Delta) -  \fdp_\lambda(\Pi^\Delta) \right|, \left| \TPP_{\lambda'}(\Pi) -  \tpp_{\lambda'}(\Pi) \right|, \left| \FDP_{\lambda'}(\Pi) -  \fdp_{\lambda'}(\Pi) \right|
	\end{equation}
	are all smaller than $\nu/2$ for all $c < \lambda, \lambda' < C$, with probability tending to one as $n, p \goto \infty$. On this event, it is easy to check that $\TPP_{\lambda}(\Pi^\Delta) \le \TPP_{\lambda'}(\Pi)$ implies $\tpp_{\lambda}(\Pi^\Delta) < \tpp_{\lambda'}(\Pi) + \nu$, and $\FDP_{\lambda}(\Pi^{\Delta}) \ge \FDP_{\lambda'}(\Pi)$ implies $\fdp_{\lambda}(\Pi^{\Delta}) > \fdp_{\lambda'}(\Pi) - \nu$. As such, in the event \eqref{eq:very_close2}, the impossibility of \eqref{eq:not_hold2} uniformly for all $c < \lambda, \lambda' < C$ implies the impossibility of 
	\[
	\TPP_{\lambda}(\Pi^\Delta) \le \TPP_{\lambda'}(\Pi) \text{ and } \FDP_{\lambda}(\Pi^{\Delta}) \ge \FDP_{\lambda'}(\Pi)
	\]
	for all $c < \lambda, \lambda' < C$.
\end{proof}
It is the same as the comment after the proof of \Cref{thm:upper}, we can prove part (b) of \Cref{thm:lower} similarly, and we omit for simplicity.

In closing, we present the following lemma to be self-contained. It shows that the lower boundary is strictly below any trade-off curve, on which the proof of \Cref{thm:lower} relies.
\begin{lemma}[Lemma C.3 in \cite{su2017false}]
	\label{lem:lower_strict}
	Consider any $\epsilon$-sparse prior $\Pi$. The lower boundary $q^\Delta$ is strictly below the trade-off curve $q^\Pi(\cdot)$, that is, $q^\Delta(u) < q^\Pi(u)$ for any $u$.
\end{lemma}
\begin{proof}[Proof of Lemma \ref{lem:lower_strict}]
	This is just a re-statement of Lemma C.3 in \cite{su2017false}. They have proved that for any $\tpp = u$, $\fdp > q^\Delta(u)$, which implies $q^\Pi(u) > q^\Delta(u)$.
\end{proof}


\end{document}